\documentclass[11
pt,a4paper]{article}
\usepackage[utf8]{inputenc}
\usepackage[francais]{babel}
\usepackage[T1]{fontenc}
\usepackage{amsmath}
\usepackage{amsfonts}
\usepackage{amssymb}
\usepackage{amsthm} 
\usepackage[all]{xy} 
\usepackage[a4paper]{geometry}
\usepackage{mathrsfs}
\usepackage{frcursive}
\usepackage{color}
\usepackage{tikz-cd}
\usepackage{mathtools}
\usepackage{graphicx}

 \theoremstyle{plain}
 \newtheorem{thm}{Théor\`{e}me}[section]
 \newtheorem{defi}[thm]{D\'{e}finition}
 \newtheorem{prop}[thm]{Proposition}
 \newtheorem{lem}[thm]{Lemme}
 \newtheorem{coro}[thm]{Corollaire}

 \theoremstyle{remark}
 \newtheorem{rem}[thm]{Remarque}
 
 \newtheorem{ex}[thm]{Exemple}

\author{Sylvain Gaulhiac}
\title{Reconstruction anabélienne du squelette des courbes analytiques}

\date{}

 \font\bf= cmbx10 at 10pt

 \newcommand{\ds}{\displaystyle}

 \newcommand{\A}{\mathbb A}
 \newcommand{\C}{\mathbb{C}}
 \newcommand{\F}{\mathbb{F}}
 \newcommand{\G}{\mathbb{G}}
 \newcommand{\Gm}{\mathbb{G}_{\mathrm {m}}}

 \newcommand{\N}{\mathbb{N}}
 \newcommand{\PP}{\mathbb{P}}
 \renewcommand{\P}{\PP}
 \newcommand{\Q}{\mathbb{Q}}
 \newcommand{\R}{\mathbb{R}}
 \newcommand{\Z}{\mathbb{Z}}

\newcommand{\OO}{\mathscr{O}}

 \newcommand{\Aut}{\mathrm{Aut}}

\newcommand{\Gal}{\mathrm{Gal}}
 
 \newcommand{\Hom}{\mathrm{Hom}}
 \newcommand{\im}{\mathrm{im}}

 \newcommand{\Spec}{\mathrm{Spec}}

\newcommand{\car}{\mathrm{car}}

\newcommand{\tp}{\mathrm{temp}}
\newcommand{\ob}{\mathrm{Ob}}
\newcommand{\cov}{\mathrm{cov}}

\newcommand{\an}{\mathrm{an}}
\newcommand{\br}{\mathrm{br}}
\newcommand{\Cov}{\mathrm{Cov}}
\newcommand{\Covt}{\mathrm{Cov}^\mathrm{t}}

\newcommand{\ttt}{\mathrm{t}}
\newcommand{\fet}{\mathrm{f\acute{e}t}}
\newcommand{\pp}{\, (p')}
\newcommand{\St}{\mathrm{St}}
\newcommand{\ta}{\mathrm{tame}}
\newcommand{\HH}{\mathrm{H}}
\newcommand{\et}{\mathrm{\acute{e}t}}
\newcommand{\Kum}{\mathsf{Kum}}
\newcommand{\sep}{\mathrm{sep}}
\newcommand{\s}{\mathrm{s}}
\newcommand{\Mod}{\mathrm{Mod}}

\newcommand{\HC}{\mathscr{H}}

 \newcommand{\EE}{\mathcal{E}}
 \newcommand{\VV}{\mathcal{V}}
 \newcommand{\BB}{\mathcal{B}}
 \newcommand{\GG}{\mathcal{G}}
 \newcommand{\XX}{\mathcal{X}}
 \newcommand{\YY}{\mathcal{Y}}
 \newcommand{\TT}{\mathcal{T}}
 \newcommand{\HHH}{\mathcal{H}}
 \newcommand{\CC}{\mathcal{C}}
 \newcommand{\UU}{\mathcal{U}}
 \newcommand{\DD}{\mathcal{D}}
  \newcommand{\E}{\mathsf{E}}
 
 \newcommand{\CX}{\mathcal{C}(X,S)}
 \newcommand{\DX}{\mathcal{D}(X,S)}

 \newcommand{\iso}{\xrightarrow{\sim}}

\begin{document}

\maketitle

\renewcommand{\abstractname}{Résumé}
\begin{abstract} Ce travail met en évidence certains comportements \emph{anabéliens} des courbes analytiques au sens de la géométrie de Berkovich. Nous montrons que la connaissance du \emph{groupe fondamental tempéré} de certaines courbes appelées \emph{analytiquement anabéliennes} détermine leurs \emph{squelettes analytiques} en tant que graphes. Le fameux demi-plan de Drinfeld en est un exemple. Le groupe fondamental tempéré d'un espace de Berkovich, introduit par André, a permis à Mochizuki de démontrer le premier résultat anabélien en géométrie de Berkovich, concernant les analytifiées de \emph{courbes hyperboliques} sur $\overline{\mathbb{Q}}_p$. \`{A} cette fin, Mochizuki développe le langage des \emph{semi-graphes d'anabélioïdes} et des \emph{tempéroïdes}. Cet article consiste à associer un semi-graphe d'anabélioïdes à une courbe analytique munie d'une \emph{triangulation minimale}, puis s'inspirer des résultats de Mochizuki afin de retrouver le squelette analytique de la courbe. La nouveauté de ces résultats est que les courbes que nous étudions ne sont plus supposées de \emph{nature algébrique}.

\bigskip
\textbf{Mots-clés :} Géométrie anabélienne, espaces de Berkovich, groupe fondamental tempéré, courbes analytiques, semi-graphe d'anabélioïdes.
\end{abstract}\

\renewcommand{\abstractname}{Abstract}
\begin{abstract} \textbf{Reconstruction of the skeleton of analytic curves.} 
This work brings to light some \emph{anabelian behaviours} of analytic curves in the context of Berkovich geometry. We show that the knowledge of the \emph{tempered fundamental group} of some curves called \emph{analytically anabelian} determines their \emph{analytic skeletons} as graphs. The famous Drinfeld half-plane is an example of such a curve. The tempered fundamental group of a Berkovich space, introduced by André, enabled Mochizuki to prove the first anabelian result in Berkovich geometry, dealing with analytifications of \emph{hyperbolic curves} over $\overline{\Q}_p$. To that end, Mochizuki develops the language of \emph{semi-graphs of anabelioids} and \emph{temperoids}. This article consists in associating a semi-graph of anabelioids to a Berkovich curve equipped with a \emph{minimal triangulation} and in adapting the results of Mochizuki in order to recover the analytic skeleton. The novelty here is that the curves we are interested in are not supposed anymore to be of \emph{algebraic nature}.

\end{abstract}\

\newpage
\tableofcontents\

\newpage
\section*{Introduction}\

Dans la lignée des développements de la géométrie anabélienne dans le cas des courbes hyperboliques, il est loisible de se poser des questions d'inspiration anabélienne dans les espaces des Berkovich, à savoir :\begin{center}
\textit{Dans quelle mesure un espace $k$-analytique est-il déterminé par son \og\emph{groupe fondamental}\fg ?}
\end{center}
 
Il est possible de définir en théorie de Berkovich plusieurs \og groupes fondamentaux\fg\, selon que ceux-ci classent les revêtements \emph{topologiques}, \emph{étales finis}, ou \emph{étales}. Néanmoins le groupe fondamental \emph{topologique} est souvent \og trop petit\fg\,tandis que le groupe fondamental \emph{étale} a tendance à être \og trop grand\fg\,. 
Le groupe qui semble au mieux capturer les comportements anabéliens des espaces analytiques sur les corps non archimédiens est le groupe fondamental \emph{tempéré}, introduit par Yves André dans \cite{And}. Ce groupe classe les revêtements dits \emph{tempérés}, définis comme les revêtements étales devenant topologiques après un changement de base fini étale, ce qui équivaut à demander que ces revêtements soient quotient d'un revêtement topologique d'un revêtement étale fini. Les revêtements topologiques tout autant que les revêtements étales finis sont des exemples de revêtements tempérés. \\

C'est Yves André, dans son article \cite{And1}, qui obtient pour la première fois des résultats de nature anabélienne en utilisant les espaces de Berkovich. À l'aide du groupe fondamental tempéré, il donne une description de nature géométrique des sous-groupes de Galois locaux $G_{\Q_p}$ du groupe de Galois absolu $G_\Q$ dans l'esprit de l'\textit{Esquisse d'un programme} \cite{Gro} de Grothendieck,
puis introduit un avatar $p$-adique du groupe de Grothendieck-Teichmüller $\widehat{GT}$. À la suite de cet article, ce sont Shinichi Mochizuki puis Emmanuel Lepage qui obtiennent les résultats les plus significatifs de géométrie anabélienne dans les espaces de Berkovich. Leurs résultats relient le groupe fondamental tempéré de l'analytifiée d'une courbe algébrique hyperbolique au graphe dual de sa réduction stable. Si $X$ est une courbe hyperbolique définie sur un corps ultramétrique complet $k$ algébriquement clos, nous savons que le type d'homotopie de l'analytifiée $X^\an$ (au sens des espaces de Berkovich) peut être décrit à partir du modèle stable $\mathscr{X}$ de $X$. Plus précisément, si $\mathscr{X}_s$ désigne la fibre spéciale de $\mathscr{X}$, le \emph{graphe dual de la réduction stable}, noté $\G_X$, est le graphe fini dont les sommets correspondent aux composantes irréductibles de $\mathscr{X}_s$, et dont les arêtes entre les sommets correspondent aux nœuds (i.e. aux  singularités en les points doubles ordinaires, qui sont les seules singularités de $\mathscr{X}_s$ par hypothèse de stabilité) entre les composantes irréductibles (une arête boucle sur le même sommet lorsqu'il s'agit d'un point d'auto-intersection). Notons $\overline{X}$ la compactification normale de $X$. En définissant les \emph{cusps} de $X$ comme l'ensemble des points de $\overline{X}$ n'appartenant pas à $X$, il est possible d'enrichir $\G_X$ en un semi-graphe $\G_X^\mathtt{c}$, en ajoutant une arête ouverte (i.e. reliée seulement à un sommet en un bout, ouverte en l'autre bout) pour chaque cusp. Il existe alors un plongement topologique canonique $\G_X^\mathtt{c}\hookrightarrow X^\an$ qui admet une rétraction par déformation $X^\an\twoheadrightarrow \G_X^\mathtt{c}$ topologiquement propre (\cite{Ber4}). L'espace analytique $X^\an$ a donc le même type d'homotopie que $\G_X$ (et que $\G_X^\mathtt{c}$). Le graphe $\G_X$ ne garde pas la trace des cusps, tandis que ceux-ci jouent un rôle dans le modèle stable. Dans le cas où $X$ est la courbe hyperbolique obtenue en retirant $3$ points à $\P_k^1$, $\G_X$ est réduit à un singleton qui sera un point de type $2$ dans $X^\an$, alors que $\G_X^\mathtt{c}$ est constitué de trois arêtes ouvertes reliées à ce même sommet. 

\subsubsection*{Résultats de Mochizuki et  Lepage}

Grâce au langage des \emph{anabélioïdes}, des \emph{semi-graphes d'anabélioïdes} et des \emph{tempéroïdes} qu'il introduit en toute généralité dans \cite{M2} et \cite{M3} afin de donner une saveur géométrique aux objets définis catégoriquement telles les catégories galoisiennes, Mochizuki prouve dans \cite{M3} que le groupe fondamental tempéré de l'analytifiée d'une courbe hyperbolique détermine le graphe dual de sa réduction stable : 

\begin{thm}[\cite{M3}, Corollary $3.11$]
Si $X_1$ et $X_2$ sont deux courbes hyperboliques sur $\overline{\Q}_p$, alors tout isomorphisme extérieur $\varphi : \pi_1^\tp(X_{1, \C_p}^\an) \simeq \pi_1^\tp(X_{2, \C_p}^\an)$ détermine, fonctoriellement en $\varphi$, un unique isomorphisme de graphes entre les graphes duaux étendus de leurs réductions stables : $\overline{\varphi} : \G_{X_1}^\mathtt{c} \simeq  \G_{X_2}^\mathtt{c}$.
\end{thm}

Mochizuki montre plus précisément qu'il est possible de reconstruire le graphe $\G_X$ de la réduction stable d'une courbe hyperbolique $X$ à partir d'une $(p')$-version $\pi_1^{\tp, \,\pp}(X^\an)$ du groupe fondamental tempéré qui se construit à partir de celui-ci en prenant la limite projective des quotients ayant un sous-groupe distingué sans torsion d'indice fini premier à $p$. Ce groupe classe l'ensemble des revêtements (tempérés) qui deviennent topologiques après un changement de base étale fini de degré premier à $p$. La reconstruction du graphe à partir de ce groupe est la suivante :

\begin{itemize}
\item[•]les sommets de $\G_X$ correspondent aux classes de conjugaison des sous-groupes compacts maximaux de $\pi_1^{\tp, \, \pp}(X^\an)$;
\item[•]les arêtes de $\G_X$ correspondent aux classes de conjugaison des intersections non triviales de deux sous-groupes compacts maximaux de $\pi_1^{\tp, \, \pp}(X^\an)$.
\end{itemize}\

Quelques années plus tard, en $2010$, Emmanuel Lepage précise le résultat de Mochizuki en montrant que si le groupe fondamental tempéré détermine en effet le graphe dual de la réduction stable d'une courbe hyperbolique, il le détermine également dans certains cas en tant que graphe métrique, c'est-à-dire que la connaissance du groupe fondamental tempéré permet de retrouver la métrique naturelle qui existe sur le graphe dual $\G_X$ de la réduction stable. Cette métrique est telle que la longueur d'une arête correspondant à un nœud est la largeur de la couronne correspondant à la fibre générique de la complétion formelle en ce nœud. Une condition suffisante assurant un tel résultat de reconstruction métrique est de se restreindre aux courbes dites \emph{de Mumford}, c'est-à-dire aux courbes algébriques propres dont les composantes irréductibles normalisées de la réduction stable sont isomorphes à $\P^1$ (ce qui en théorie de Berkovich équivaut à dire que l'analytifié $X^\an$ est localement isomorphe à un ouvert de $\P^{1, \an}$, ou encore que le genre topologique de $X^\an$ est le genre de $X$, ou encore que $X^\an$ admet une uniformisation par un ouvert de $\P^{1, \an}$, ou encore qu'il n'a pas de point de type $2$ de genre $>0$ ...).

\begin{thm}[\cite{Lep$1$}, Th. $4.12$; \cite{Lep$2$}, Th. $12$]\label{théorème lepage}
Soient $X_1$ et $X_2$ deux courbes hyperboliques de Mumford sur $\overline{\Q}_p$, et $\varphi : \pi_1^\tp(X_{1}^\an) \simeq \pi_1^\tp(X_{2}^\an)$ un isomorphisme. Alors l'isomorphisme entre les graphes duaux des réductions stables : $\overline{\varphi} : \G_{X_1} \simeq  \G_{X_2}$ est un isomorphisme de \emph{graphes métriques}.
\end{thm}

Lepage va plus loin dans son article \cite{Lep$3$}. En utilisant des techniques de résolution des non-singularités et en prouvant qu'une courbe de Mumford sur $\overline{\Q}_p$ vérifie une telle résolution, il montre que le groupe tempéré de l'analytifiée d'une courbe de Mumford détermine son espace topologique sous-jacent. Il utilise pour cela un homéomorphisme entre l'analytifiée de toute $\overline{\Q}_p$-courbe propre et lisse et la limite projective des graphes duaux de toutes les réductions semi-stables de la courbe.

\begin{thm}[\cite{Lep$3$}, Theorem  $3.9$]
Soient $X_1$ et $X_2$ deux courbes hyperboliques de Mumford sur $\overline{\Q}_p$, et $\varphi : \pi_1^\tp(X_{1}^\an) \simeq \pi_1^\tp(X_{2}^\an)$ un isomorphisme. Alors $\varphi$ induit un homéomorphisme entre les analytifiées : $$X_{1}^\an\iso X_{2}^\an.$$
\end{thm}

Dans le cas des courbes de Tate, il obtient même mieux en montrant, toujours grâce à la résolution des non-singularités, que le groupe fondamental tempéré d'une courbe de Tate épointée détermine l'espace analytique de la courbe de Tate.

\begin{thm}
Soient $q_1, q_2\in \overline{\Q}_p$, tels que $\vert q_1\vert <1$ et $\vert q_2\vert <1$. Pour $i\in \lbrace 1,2\rbrace$, posons $E_i=\G_m^\an/q_i^\Z$ et $X_i=E_i\setminus \lbrace 1\rbrace$. Si les groupes tempérés $\pi^\tp(X_1)$ et $\pi^\tp(X_2)$ sont isomorphes, alors il existe $\sigma\in G_{\Q_p}=\Gal(\overline{\Q}_p/\Q_p)$ tel que $\sigma(q_1)=q_2$, et $E_1$ et $E_2$ sont isomorphes en tant qu'espaces $\Q_p$-analytiques. 
\end{thm}

\subsection*{Contenu de l'article}

Ces résultats concernent les courbes analytiques de \emph{nature algébrique}, c'est-à-dire les analytifiées de courbes algébriques. Or, la théorie des espaces analytiques de Berkovich est suffisamment riche pour faire apparaître toutes sortes de courbes qui sont essentiellement de \emph{nature analytique}, sans provenir directement de courbes algébriques. Les exemples les plus simples mais non moins importants de telles courbes sont les \emph{couronnes} et les \emph{disques} analytiques. Si $k$ est un corps ultramétrique complet, les couronnes (respectivement les disques) $k$-analytiques sont les espaces $k$-analytiques isomorphes au domaine $k$-analytique de $\P_k^{1, \an}$ défini par la condition $\{\vert T \vert \in I \}$ pour un certain intervalle non vide $I$ de $\R_+^\times$ (resp. de $\R_+$ et contenant $0$ et rencontrant $\R_+^*$). Dans le sillage des résultats de Mochizuki et de Lepage, il est naturel de se demander si des résultats anabéliens similaires existent pour les courbes $k$-analytiques sans imposer que celles-ci soient de nature algébrique. 
\bigskip

Les résultats évoqués ci-dessus faisant apparaître le graphe dual de la réduction stable $\G_X$, il convient d'en trouver pour les courbes $k$-analytiques un analogue défini uniquement dans le langage de la géométrie analytique sans requérir une quelconque nature algébrique. C'est le \emph{squelette analytique} qui jouera ce rôle : si $X$ est une courbe $k$-analytique quasi-lisse, le \emph{squelette analytique} de $X$, noté $S^\an(X)$, est défini comme l'ensemble des points de $X$ qui n'appartiennent à aucun disque analytique ouvert. Le squelette $S^\an(X)$ est alors un semi-graphe localement fini. Remarquons  que $S^\an(X)=\emptyset$ lorsque $X=\P_k^{1, \an}$ et que $\P_k^{1, \an}$ est l'unique courbe $k$-analytique quasi-lisse, connexe et compacte de squelette analytique vide lorsque $k$ est algébriquement clos. Par ailleurs, le squelette analytique est homéomorphe à un cercle pour une \emph{courbe de Tate}, c'est-à-dire l'analytifiée d'une courbe elliptique à mauvaise réduction. Dès que l'application naturelle $\pi_0(S^\an(X))\rightarrow \pi_0(X)$ est surjective, il existe une rétraction par déformation $X\twoheadrightarrow S^\an(X)$ topologiquement propre. En particulier, dans ce cas, $X$ a le type d'homotopie de $S^\an(X)$. \\

Lorsque $k$ est algébriquement clos, nous appellerons \emph{courbe $k$-analytiquement hyperbolique} (\ref{définition des courbes marquées anabéliennes}) toute courbe $k$-analytique $X$, quasi-lisse et connexe, telle que l'ensemble des nœuds (\ref{noeuds}) de son squelette analytique soit non vide et constitué uniquement de \emph{noeuds hyperboliques} (\ref{noeuds hyperboliques}). Une courbe \emph{$k$-analytiquement anabélienne} sera alors définie comme une courbe $k$-analytiquement hyperbolique vérifiant une condition technique, dite d'\emph{ascendance vicinale} (\ref{définition des courbes analytiquement anabéliennes}), permettant la reconstruction des arêtes ouvertes. L'analytifiée $\mathscr{X}^\an$ d'une $k$-courbe algébrique $\mathscr{X}$ est $k$-analytiquement hyperbolique si et seulement si $\mathscr{X}$ est une courbe hyperbolique, auquel cas $\mathscr{X}^\an$ est $k$-analytiquement anabélienne. Nous montrons dans ce papier le résultat suivant lorsque $k$ est un corps non archimédien complet, non trivialement valué et algébriquement clos : 

\bigskip
\textbf{Théorème A : }
\textit{Soient $X_1$ et $X_2$ deux courbes $k$-analytiquement anabéliennes. Tout isomorphisme entre les groupes tempérés : $\varphi : \pi_1^\tp(X_1)\iso \pi_2^\tp(X_2)$ induit (fonctoriellement en $\varphi$) un isomorphisme de semi-graphes entre les squelettes analytiques associés : $\overline{\varphi}: S^\an (X_1)\iso S^\an (X_2).$}\\

Pour certaines courbes anabéliennes, dites \emph{à arêtes relativement compactes} (cf. \ref{définition courbes hyperboliques à aretes relativement compactes}), l'isomorphisme $\varphi$ du théorème ci-dessus peut être remplacé par l'hypothèse moins forte d'un isomorphisme entre les versions premières à $p$, $\pi_1^{\tp, \pp}(X_i)$ (cf. \ref{p' version du groupe fondamental}), des groupes tempérés $\pi_1^\tp(X_i)$.\\

Un exemple intéressant de courbe $k$-analytiquement anabélienne (à arêtes relativement compactes) est celui du \emph{demi-plan de Drinfeld}, d'un intérêt majeur en théorie des nombres et analogue non archimédien du demi-plan de Poincaré.
Soit $k_0$ un corps non-archimédien \emph{local} de caractéristique positive (i.e. $k_0$ est soit une extension finie de $\Q_p$, soit du type $\F_{q}((T))$) contenu dans le corps ultramétrique complet et algébriquement clos $k$ dont la valuation prolonge celle de $k_0$. Le \emph{demi-plan de Drinfeld}, noté $\mathfrak{D}^1_{k/k_0}$, est la courbe $k$-analytique définie en tant qu'ouvert de $ \P_{k}^{1, \an}$ par :
 
 $$\mathfrak{D}^1_{k/k_0}= \P_{k}^{1, \an}\setminus  \P_{k}^{1, \an}(k_0).  $$
 
 Il s'agit d'une courbe $k$-analytiquement anabélienne, non compacte, dont le squelette analytique n'est autre que l'\emph{immeuble de Bruhat-Tits} de $\mathrm{SL}_2/k_0$, noté $\mathfrak{B}^1_{k_0}$. L'anabélianité du demi-plan de Drinfeld nous permet d'énoncer le résultat suivant : 
 
\bigskip
\textbf{Proposition B} :
\textit{Pour $i\in \lbrace 1,2 \rbrace$, soit $p_i$ un nombre premier et $k_i$ un corps de nombres $p_i$-adique, c'est-à-dire une extension finie de $\Q_{p_i}$. Supposons les groupes suivants isomorphes : $${\pi_1^{\tp, (p_1')}(\mathfrak{D}^{1}_{\C_{p_1}/k_1})\simeq \pi_1^{\tp, (p_2')}(\mathfrak{D}^{1}_{\C_{p_2}/k_2})},$$
alors $p_1=p_2$, et les extensions $k_1/\Q_{p_1}$ et $k_2/\Q_{p_2}$ ont même degré d'inertie : $[\widetilde{k_1} : \F_{p_1}]=[\widetilde{k_2} : \F_{p_2}]$, d'où un isomorphisme $\widetilde{k_1}\simeq \widetilde{k_2}$.}\\

Certaines courbes $k$-analytiquement anabéliennes peuvent par ailleurs être hyperboliques tout en ayant dans leur squelette analytique des arêtes ouvertes (appelées cusps). La quête de la propriété d'anabélianité parmi les courbes $k$-analytiquement hyperboliques nous amène à définir deux types de cusps privilégiés : les cusps \emph{coronaires finis} et \emph{discaux épointés} (\ref{définition cusps coronaires finis et discaux épointés}). En traitant les cusps coronaires finis à l'aide d'une condition de total déploiement de certains $\mu_p$-torseurs (\ref{corollaire permettant de trouver un revêtement d'une couronne totalement décomposé en un de ses bouts}) et les cusps discaux épointés par un résultat de Mochizuki (\ref{anabélianité des analytifiées des courbes hyperboliques}), nous obtenons finalement le théorème suivant :

\bigskip
\textbf{Théorème C}
 \textit{Soit $k$ un corps ultramétrique complet algébriquement clos et non trivialement valué. Alors les trois types de courbes ci-dessous sont $k$-analytiquement anabéliennes : }

 \begin{itemize}
 \item[•]\textit{les courbes $k$-analytiquement hyperboliques à arêtes relativement compactes; }
 \item[•]\textit{si $k$ est de caractéristique mixte : les courbes $k$-analytiquement hyperboliques dont tous les cusps sont coronaires finis;}
 \item[•]\textit{si $k$ est de caractéristique mixte : les courbes $k$-analytiquement hyperboliques dont le squelette analytique est un semi-graphe fini et dont les cusps sont coronaires finis ou discaux épointés. }
 \end{itemize}\
 
Ces cas ne traitent néanmoins pas toutes les courbes $k$-analytiques hyperboliques, et nous ne savons pas pour le moment si toutes les courbes $k$-analytiquement hyperboliques sont $k$-analytiquement anabéliennes. 
\bigskip 
\bigskip

Nous présentons dans $\S 1$ le langage des semi-graphes d'anabélioïdes et des tempéroïdes que Mochizuki introduit dans \cite{M3}, afin d'énoncer le théorème \ref{reconsmochi} (\cite{M3}, Corollary $3.9$) selon lequel, sous certaines hypothèses que nous qualifierons de \emph{mochizukiennes}, un graphe d'anabélioïdes $\GG$ est déterminé par le tempéroïde connexe $\BB^\tp(\GG)$ qui lui est associé. 

La partie $\S 2$ se concentre sur les courbes $k$-analytiques. Après avoir introduit dans \ref{courbes analytiques} quelques notions de base sur les courbes $k$-analytiques, nous définissons dans \ref{STN} la notion de triangulation, de triangulation généralisée, de squelette et de nœuds d'une courbe $k$-analytique, en nous inspirant très largement de l'étude systématique des courbes analytiques entreprise par A. Ducros dans \cite{Duc}. Cela nous permet d'associer à une courbe $k$-analytique quasi-lisse munie d'une triangulation généralisée $S$ des graphe et semi-graphe d'anabélioïdes $\GG^\natural(X,S)$ et $\GG(X,S)$ que l'on définit à partir des revêtements \emph{modérés} autour des points de la triangulation. La partie \ref{groupe fondamental tempéré} est l'occasion de définir les revêtements tempérés ainsi que le groupe fondamental tempéré d'un espace strictement $k$-analytique. 

C'est dans ce nouveau cadre analytique que nous transposons, dans $\S 3$, les résultats anabéliens de Mochizuki de reconstruction d'un graphe d'anabélioïdes à partir du tempéroïde connexe associé. Nous entreprenons dans \ref{section description} la description du semi-graphe d'anabéloïdes $\GG(X,S)$. La description des composantes coronaires découle d'un résultat de Berkovich sur l'aspect kummérien des revêtements modérés des couronnes analytiques, tandis que la description des composantes sommitales, plus délicate, repose sur la trivialité des revêtements modérés des disques (\cite{Ber4}), mais également sur la description des revêtements autour d'un point. Nous introduisons pour cela la \emph{courbe résiduelle} $\mathscr{C}_x$ associée à un point $x$ de type $2$, ainsi que des résultats assez fins de \cite{Duc} sur les valuations associées aux branches et l'application bijective de l'ensemble des branches issues de $x$ sur l'ensemble des points fermés d'un ouvert de Zariski de la courbe résiduelle  $\mathscr{C}_x$. La partie \ref{Compatibilité} consiste à comprendre sous quelles contraintes sur la courbe $X$ et la triangulation généralisée $S$ le graphe d'anabélioïdes $\GG^\natural(X,S)$ vérifie l'ensemble des conditions mochizukiennes. Ces hypothèses étant malaisées à vérifier directement, nous nous servirons d'une remarque de Mochizuki (\cite{M3}, Example $2.10$) selon laquelle l'ensemble des hypothèses mochizukiennes sont vérifiées sous des conditions plus simples, et notre travail consiste alors à montrer sous quelles contraintes sur $X$ et $S$ ces conditions sont vérifiées par $\GG^\natural(X,S)$. Ces considérations nous amènent à définir l'ensemble des courbes \emph{$k$-analytiquement hyperboliques} (\ref{définition des courbes marquées anabéliennes}), puis \emph{$k$-analytiquement anabéliennes} (\ref{définition des courbes analytiquement anabéliennes}), avant d'énoncer et de prouver nos résultats.\\

\def\restriction#1#2{\mathchoice
              {\setbox1\hbox{${\displaystyle #1}_{\scriptstyle #2}$}
              \restrictionaux{#1}{#2}}
              {\setbox1\hbox{${\textstyle #1}_{\scriptstyle #2}$}
              \restrictionaux{#1}{#2}}
              {\setbox1\hbox{${\scriptstyle #1}_{\scriptscriptstyle #2}$}
              \restrictionaux{#1}{#2}}
              {\setbox1\hbox{${\scriptscriptstyle #1}_{\scriptscriptstyle #2}$}
              \restrictionaux{#1}{#2}}}
\def\restrictionaux#1#2{{#1\,\smash{\vrule height .8\ht1 depth .85\dp1}}_{\,#2}}

\bigskip
\begin{center}
 \large
 R\normalsize EMERCIEMENTS
 \end{center} 
 
  Je tiens à exprimer toute ma reconnaissance à Antoine Ducros pour sa disponibilité lors de nos échanges mathématiques et pour ses nombreuses idées et commentaires dont la pertinence a souvent éclairé des points qui m'étaient obscurs. Je tiens également à remercier Emmanuel Lepage dont vient l'idée de départ de cet article et des méthodes mises en œuvre, ainsi que Shinichi Mochizuki qui a pris le temps de répondre consciencieusement à certaines de mes questions.
  \\
 
\section{Résultats de Mochizuki}

Nous allons présenter le langage des semi-graphes d'anabélioïdes et de tempéroïdes utilisé par Mochizuki dans \cite{M1}, \cite{M2} et \cite{M3} .\\

Commençons par définir un \emph{semi-graphe} $\G$ qui consiste en une collection de données suivante : 
\begin{center}
\begin{itemize}
\item[•] un ensemble $\mathcal{V}(\G)$ dont les éléments sont appelés \emph{sommets};
\item[•] un ensemble $\mathcal{E(\G)}$ dont les éléments sont appelés \emph{arêtes} et consistent en des ensembles de cardinal au plus $2$ dont les éléments sont les \emph{branches}, et tels que si $e, e'\in \mathcal{E}(\G)$ :  $e\neq  e' \Rightarrow  e\cap e'=\emptyset$;
\item[•] une collection $\zeta (\G)$ de flèches $\left( \zeta_e\right)_{e \in\mathcal{E}(\G)}$ avec $\zeta_e : e\rightarrow \mathcal{V}(\G)$
\end{itemize}
\end{center}

Si $v=\zeta_e(b)$ pour une arête $e$, une branche $b$ de $e$ et un sommet $v$, on dira que $e$ (ou $b$) \emph{aboutit} à $v$. Un semi-graphe $\G$ sera dit \emph{fini} (resp. \emph{dénombrable}) si les deux ensembles $\mathcal{V}(\G)$ et $\mathcal{E}(\G)$ sont finis (resp. dénombrables), et \emph{localement fini} si pour chaque sommet $v\in \mathcal{V}(\G)$ l'ensemble des arêtes aboutissant à $v$ est fini. On parlera de \emph{graphe} lorsque toutes les arêtes contiennent deux branches.\\

Un morphisme entre deux semi-graphes $$\G=\lbrace\mathcal{V}(\G), \mathcal{E}(\G), \zeta (\G)\rbrace\rightarrow \G'=\lbrace\mathcal{V}(\G'), \mathcal{E}(\G'), \zeta (\G')\rbrace$$
est une collection d'applications $\mathcal{V}(\G)\rightarrow \mathcal{V}(\G')$, $\mathcal{E}(\G)\rightarrow \mathcal{E}(\G')$ qui rend compatibles entre elles les collections $\zeta (\G)$ et $\zeta (\G')$ et telle que si l'arête $e$ de $\G$ est envoyée sur l'arête $e'$ de $\G'$, alors on a une injection $e\hookrightarrow e'$. Une arête de $\G$ ne peut donc pas être envoyée sur une arête de $\G'$ qui a un nombre strictement inférieur de branches.

\subsection{Le monde des anabélioïdes}\label{le monde des anabélioïdes}

Si $G$ est un groupe profini, nous noterons $\mathcal{B}(G)$ la catégorie des ensembles finis sur lesquels $G$ opère continûment (i.e. de manière à ce que les stabilisateurs des points soient des sous-groupes ouverts de $G$). Soit $\mathfrak{Ens}^\mathrm{f}$ la catégorie des ensembles finis, et $$\beta_G^* : \mathcal{B}(G)\rightarrow \mathfrak{Ens}^\mathrm{f}$$ le foncteur « oubli » défini en oubliant l'action de $G$. Le foncteur $\beta_G^*$ est \emph{exact} (i.e. commute aux colimites finies et limites finies) et \emph{conservatif} (i.e. un morphisme $\alpha$ de $\mathcal{B}(G)$ est un isomorphisme si et seulement si $\beta_G^*(\alpha)$ est un isomorphisme de $\mathfrak{Ens}^\mathrm{f}$).\\

Rappelons (voir \cite{Cad} ou \cite{SGA1}, Exposé $5$ pour une exposition détaillée) qu'une catégorie $\XX$ est appelée \emph{catégorie galoisienne} si elle est équivalente à une catégorie de la forme $\mathcal{B}(G)$ avec $G$ un groupe profini. \\

Grothendieck en donne également une caractérisation ne faisant appel à aucun groupe profini. Plus précisément, en vertu de \cite{SGA1} (Exposé $5$, Th. $4.1$), une catégorie $\XX$ est galoisienne si elle seulement si elle vérifie les conditions suivantes : 
\begin{itemize}
\item[•]les limites finies et les colimites finies existent dans $\XX$,
\item[•]il existe un foncteur covariant $\beta^* : \XX \rightarrow\mathfrak{Ens}^\mathrm{f}$ qui est \emph{exact} et \emph{conservatif}, 
\item[•]tout morphisme $u : X\rightarrow Y$ de $\XX$ se factorise par $X \xrightarrow{u'} Y' \xrightarrow{u''} Y$ où $u'$ est un épimorphisme strict et $u''$ un monomorphisme qui induit un isomorphisme sur un sommande direct de $Y$.
\end{itemize}

 Un tel foncteur $\beta^* : \XX \rightarrow\mathfrak{Ens}^\mathrm{f}$ exact et conservatif est appelé un \emph{foncteur fondamental} de la catégorie galoisienne $\XX$. La donnée d'un foncteur fondamental $\beta^*$ permet de construire une équivalence de catégorie entre $\XX$ et $\mathcal{B}(\pi_1(\XX, \beta^*))$ où le groupe profini $\pi_1(\XX, \beta^*)$ est défini par :$$\pi_1(\XX, \beta^*) :=\Aut(\beta^*),$$ équivalence par laquelle le foncteur $\beta^*$ est identifié à $\beta^*_{\pi_1(\XX, \beta^*)}$.\\

 Le groupe $\pi_1(\XX, \beta^*)$ est appelé le \emph{groupe fondamental} de la catégorie galoisienne $\XX$ munie du foncteur fondamental $\beta^*$. Deux foncteurs fondamentaux sont toujours isomorphes, et un tel isomorphisme induit fonctoriellement un isomorphisme de groupes fondamentaux. Par ailleurs, tout automorphisme d'un foncteur fondamental induit un automorphisme intérieur du groupe fondamental. De là, la classe d'isomorphismes de $\pi_1(\XX, \beta^*)$ ne dépend pas de $\beta^*$, et nous pouvons  considérer le \emph{groupe fondamental  $\pi_1(\XX)$ de $\XX$}, défini à unique isomorphisme extérieur près.  C'est un groupe profini. Remarquons qu'une catégorie galoisienne est connexe puisqu'elle possède un objet final, la limite vide.\\

Afin de développer une intuition plus géométrique et d'être cohérent avec la philosophie anabélienne de Mochizuki nous conserverons la terminologie de \cite{M2} en appelant \emph{anabélioïde connexe} toute catégorie $\mathcal{X}$ équivalente à une catégorie de la forme $\mathcal{B}(G)$ où $G$ est un groupe profini. Ainsi un anabélioïde connexe n'est autre qu'une catégorie galoisienne, mais dont les morphismes sont définis dans l'autre sens.

\begin{defi}\
Si $\XX$ et $\YY$ sont des anabélioïdes connexes, un \emph{morphisme d'anabélioïdes} ${\varphi : \XX\rightarrow\YY}$ est défini comme étant un foncteur exact $\varphi^* : \YY\rightarrow\XX$. Un \emph{isomorphisme d'anabélioïdes} entre deux anabélioïdes connexes est un morphisme dont le foncteur correspondant est une équivalence de catégories.
\end{defi}

\begin{ex}

Si $X$ et $Y$ sont soit deux $k$-schémas localement noethériens connexes non vides, soit deux bons espaces $k$-analytiques connexes, la catégorie $\mathrm{Cov}^\mathrm{alg}(X)$ (resp $\mathrm{Cov}^\mathrm{alg}(Y)$) des revêtements étales finis de $X$ (resp. de $Y$) est une catégorie galoisienne, donc un anabélioïde connexe. Un morphisme de $k$-schémas (ou d'espaces $k$-analytiques) $\varphi : X\rightarrow Y$ induit par changement de base un foncteur exact $\varphi^* : \mathrm{Cov}^\mathrm{alg}(Y)\rightarrow \mathrm{Cov}^\mathrm{alg}(X)$, i.e. un morphisme d'anabélioïdes $\varphi : \mathrm{Cov}^\mathrm{alg}(X)\rightarrow \mathrm{Cov}^\mathrm{alg}(Y)$, ce qui explique le choix de ce sens des flèches. Le groupe fondamental de l'anabélioïde connexe $\mathrm{Cov}^\mathrm{alg}(X)$ n'est autre que le \emph{groupe fondamental algébrique} $\pi_1^\mathrm{alg}(X)$. Ainsi : $$\mathrm{Cov}^\mathrm{alg}(X)\simeq \BB\left(\pi_1^\mathrm{alg}(X)  \right).$$

\end{ex}

\begin{rem}
Les catégories galoisiennes sont des exemples de topoi élémentaires (qui ne sont cependant pas des topoi de Grothendieck), et le sens des flèches dans la définition des morphismes d'anabélioïdes connexes coïncide avec la convention pour les morphismes géométriques de topoi élémentaires.
\end{rem}

\begin{rem}
La catégorie $\mathfrak{Ens}^\mathrm{f}$ est l'exemple le plus simple d'anabélioïde connexe (obtenu en appliquant $\BB(\--)$ au groupe trivial). Cela nous permet de considérer les foncteurs fondamentaux entre catégories galoisiennes comme des « points fondamentaux » des anabélioïdes connexes pensés comme objets géométriques.
\end{rem}

\begin{defi}
Un \emph{point de base} d'un anabélioïde connexe $\XX$ n'est autre qu'un morphisme d'anabélioïdes $\beta : \mathfrak{Ens}^\mathrm{f}\rightarrow\XX$.  
\end{defi}

\begin{rem}
Si $\mathcal{X}$ est une catégorie galoisienne, un foncteur exact de $\XX$ dans  $\mathfrak{Ens}^\mathrm{f}$ est conservatif, c'est-à-dire que c'est un foncteur fondamental. Par conséquent il existe une correspondance bijective entre les foncteurs fondamentaux de la catégorie galoisienne $\XX$ et les points de base de l'anabélioïde connexe $\XX$. Si $\beta$ est un point de base de $\XX$ (considéré en tant qu'anabélioïde) et $\beta^*$ est le foncteur associé, on définit le \emph{groupe fondamental} de l'anabélioïde $\XX$ muni de $\beta$, noté $\pi_1(\XX, \beta)$, comme égal au groupe $\pi_1(\XX, \beta^*)$. La classe d'isomorphisme de ce groupe étant indépendante du choix du point de base, on pourra considérer  $\pi_{1}(\mathcal{X})$ qui sera défini à unique isomorphisme extérieur près.
\end{rem}

\begin{rem}
Un foncteur qui définit une équivalence de catégories est nécessairement exact, par conséquent un isomorphisme d'anabélioïdes connexes $\varphi : \XX\rightarrow\YY$ n'est autre qu'une équivalence de catégories $\varphi^* : \YY\rightarrow\XX$.
\end{rem}

\begin{defi}\label{definition anabelioide general} 
Une catégorie $\XX$ est un \emph{anabélioïde} si elle est équivalente à une catégorie produit $$\XX_I:=\prod _{i\in I} \XX_i, $$ où $I$ est un ensemble fini et $\XX_i$ est un anabélioïde connexe pour tout $i\in I$. En d'autres termes, dans le langage de $\cite{SGA1}$, un anabélioïde est une \og \emph{catégorie multi-galoisienne} \fg. Un morphisme d'anabélioïdes sera défini comme dans le cas connexe, c'est-à-dire comme un foncteur exact dans le sens opposé.
\end{defi}

\begin{rem}
Les anabéloïdes forment une \emph{$2$-catégorie}. Si $\varphi, \psi\in \mathrm{Mor}(\mathcal{X},\mathcal{Y})$ sont des $(1)$-morphismes entre deux anabélioïdes $\mathcal{X}$ et $\mathcal{Y}$, un $(2)$-morphisme $\varphi\to \psi$ est défini comme une transformation naturelle entre les foncteurs $\varphi^*$ et $\psi^*$.

\end{rem}

\begin{defi} \label{definition d'une decomposition d'un objet dans une categorie} 
Soit $\mathcal{C}$ une catégorie admettant les colimites finies, possédant ainsi un objet initial $\mathbf{0}_\mathcal{C}$. 
\begin{enumerate}
\item Un objet $X\in \ob(\mathcal{C})$ différent de $\mathbf{0}_\mathcal{C}$ est dit \emph{connexe} lorsqu'il n'existe pas dans $\mathcal{C}$ d'écriture sous forme de coproduit $X=X_1\coprod X_2$ avec $X_i\neq \mathbf{0}_\mathcal{C}$ pour $i=1,2$.
\item Pour tout objet $X\neq \mathbf{0}_\mathcal{C}$ de $\mathcal{C}$, une \emph{décomposition de $X$} est une écriture de la forme : $$ X= \coprod_{\alpha\in A} X_\alpha $$ où $A$ est un ensemble fini, $X_\alpha\in \ob(\mathcal{C})$ et $X_\alpha\neq \mathbf{0}_{\mathcal{C}}$. En particulier l'objet $X$ est \emph{connexe} si et seulement si l'ensemble d'indexation de toute décomposition de $X$ est de cardinal $1$. 
\item Lorsque $\mathcal{C}$ est une catégorie galoisienne, alors $($cf. $\cite{Cad}, 3.2.1)$ tout objet $X\neq \mathbf{0}_\mathcal{C}$ de $\mathcal{C}$ peut être écrit : $$X=\coprod_{i=1}^r X_i,$$ où les $X_i\in \ob(\mathcal{C})$ sont des objets \emph{connexes} tous distincts de $\mathbf{0}_\mathcal{C}$. Cette décomposition est unique $($à permutations près$)$, et les $X_i$ sont appelées les \emph{composantes connexes de $X$}.
\end{enumerate}
\end{defi}

\begin{prop}\label{unicité de la décomposition d'un objet en composantes connexes dans une multi-galoisienne } Si $\XX_I:=\prod _{i\in I} \XX_i$ est l'anabélioïde défini en \ref{definition anabelioide general}, alors tout objet  $X\neq \mathbf{0}_\mathcal{\XX_I}$ de $\XX_I$ admet, à permutations près, une unique décomposition de la forme : $$X=\coprod_{j=1}^s X_j,$$ où les $X_j\in \ob(\mathcal{\XX_I})$ sont des objets \emph{connexes}, tous distincts de $\mathbf{0}_\mathcal{\XX_I}$, et appelés les \emph{composantes connexes de $X$}.
\end{prop}

\begin{proof}
En tant que catégorie multi-galoisienne, $\XX_I$ admet les limites et colimites finies, d'où un objet final (resp. initial) noté $\mathbf{1}_{\XX_I}$ (resp. $\mathbf{0}_{\XX_I}$), unique à unique isomorphisme près, et défini comme la limite (resp. colimite) vide.
 Pour $i\in I$, notons $\mathbf{0}_i$ et $\mathbf{1}_i$ les objets initiaux et finals de l'anabélioïde connexe $\XX_i$. Si $X_i\in \ob(\XX_i)$, notons $\widehat{X_i}\in \ob(\XX_I)$ défini en prenant le produit de $X_i$ avec les $\mathbf{0}_j$ pour $j\in I\setminus \lbrace i\rbrace$. Si $X\in \ob(\XX_I)$, notons $X[i]\in\ob(\XX_i)$ la composante de $X$ sur $\XX_i$ selon la décomposition $\XX_I:=\prod _{i\in I} \XX_i,$ de sorte que $X=\left(X[i]\right)_{i\in I}$. On a alors : $$X=\coprod_{i\in I} \widehat{X[i]}.$$
 Par conséquent un objet $X$ de $\XX_I$ n'est connexe que s'il est de la forme $\widehat{X_i}$ pour un certain $i\in I$ et $X_i \in\ob(\XX_i)$ objet connexe. En couplant ces considérations au résultat d'unicité de la décomposition d'un objet non vide (i.e. non initial) d'une catégorie galoisienne en somme finie d'objets connexes (cf. définition \ref{definition d'une decomposition d'un objet dans une categorie}, ou \cite{Cad}, $3.2.1$ pour la preuve), on obtient l'unicité voulue : soit $X\in \ob(\XX_I)$, $X\neq \mathbf{0}_{\XX_I}$, et $I'\subset I$ formé des éléments $i$ tels que $X[i]\neq \mathbf{0}_i$. Pour chaque $i\in I'$, considérons la décomposition (unique à permutation près des facteurs) de $X[i]$ en ces composantes connexes :$$X[i]=\coprod_{j\in J_i}X_{i,j}$$ où $J_i$ est un ensemble fini et $X_{i,j}\in\ob(\XX_i)$ est un objet connexe, pour tout $j\in J_i$. Il est alors clair que $X$ s'écrit : $$X=\coprod_{i\in I'}\coprod_{j\in J_i}\widehat{X_{i,j}},$$ que  tous les $\widehat{X_{i,j}}$ sont des objets connexes de $\XX_I$, et que cette décomposition est unique à permutations près.
\end{proof}

\begin{coro}\label{unicité des composantes connexes d'un anabélioide} 
Si $\XX_I$ est l'anabélioïde défini en \ref{definition anabelioide general}, alors l'ensemble fini $I$ et les anabélioïdes connexes $\XX_i$ (pour $i\in I$) sont canoniquement déterminés par l'anabélioïde abstrait $\XX_I$. Les $\XX_i$ sont appelées les \emph{composantes connexes de l'anabélioïde $\XX_I$}.
\end{coro}

\begin{proof}
Notons $\varepsilon_i=\widehat{\mathbf{1}_i}\in \ob(\XX_I)$ (rappelons qu'il est défini en prenant le produit de $\mathbf{1}_i$ avec les $\mathbf{0}_j$ pour $j\in I\setminus \lbrace i\rbrace$). Les $\varepsilon_i$ sont des objets connexes de $\XX_I$, et l'on obtient la décomposition : $$\mathbf{1}_{\XX_I}=\coprod_{i\in I}\varepsilon_i. $$
 Par la proposition \ref{unicité de la décomposition d'un objet en composantes connexes dans une multi-galoisienne } précédente, une telle décomposition est unique (à permutations des facteurs près). De là les $\varepsilon_i$ sont entièrement déterminés par $\XX_I$ en tant que composantes connexes de l'élément final $\mathbf{1}_{\XX_I}$, de même que l'ensemble $I$. Ainsi, pour chaque $i\in I$, l'anabélioïde connexe $\XX_i$ est déterminé tout aussi canoniquement par $\XX_I$ en tant que sous-catégorie pleine de $\XX_I$ formée des objets \emph{au-dessus de $\varepsilon_i$}.
\end{proof}

Soit $\XX$ un anabélioïde et $S\in \ob(\XX)$. Désignons par $\XX_S$ la catégorie des objets de $\XX$ au-dessus de $S$, c'est-à-dire la catégorie dont les objets sont les flèches $X\rightarrow S$ de $\XX$, et dont les flèches entre $X\rightarrow S$ et $X'\rightarrow S$ sont les $S$-morphismes $X\rightarrow X'$. Soit $i_S^*: \XX\rightarrow\XX_S$ le foncteur défini en prenant le produit fibré avec $S$.

\begin{prop}[\cite{M2}, Prop. $1.2.1$]
La catégorie $\XX_S$ est un anabélioïde, et le foncteur $i_S^*$ est exact, c'est-à-dire qu'il définit un morphisme d'anabélioïdes ${i_S : \XX_S\rightarrow\XX}$.
\end{prop}

\begin{rem}
Les composantes connexes de l'anabélioïde $\XX_S$ sont ici en bijection avec les composantes connexes de $S$ en tant qu'objet de $\XX$.  Par conséquent $\XX_S$ est un \emph{anabélioïde connexe} si et seulement si $S$ est lui-même un objet connexe. 
\end{rem}

\begin{defi}
Un morphisme d'anabélioïdes $\varphi : \YY\rightarrow\XX$ est un \emph{revêtement fini étale} s'il se factorise comme la composée d'un isomorphisme $\alpha : \YY\iso\XX_S$ avec le morphisme $i_S : \XX_S\rightarrow\XX$ où $S\in \ob(\XX)$.
\end{defi}

\begin{rem}
Un morphisme continu $\theta : H\rightarrow G$ entre deux groupes profinis $H$ et $G$ induit naturellement un foncteur exact $\theta^* : \BB(G)\rightarrow\BB(H)$, c'est-à-dire un morphisme d'anabélioïdes $\theta : \BB(H)\rightarrow\BB(G)$. Un tel morphisme est un revêtement fini étale si et seulement si $\theta$ est une injection qui permet d'identifier $H$ à un sous-groupe ouvert de $G$. Par ailleurs, tout revêtement fini étale d'anabélioïdes connexes peut être écrit sous cette forme sous un choix approprié des points de base à la source et au but.
\end{rem}

Nous allons maintenant définir la notion de \emph{semi-graphe d'anabélioïdes} introduite par Mochizuki dans \cite{M3}, qui est la transcription dans le langage anabélioïdique des \emph{semi-graphes de groupes profinis} introduits par le même auteur dans \cite{M1}. Nous introduirons ensuite la catégorie des \emph{recouvrements finis étales} d'un tel semi-graphe qui permettra de définir ses \emph{revêtements finis étales}.

\begin{defi}\label{définition semi-graphe d'anabélioïdes}

Un \emph{semi-graphe d'anabélioïdes (connexes)} $\GG$ consiste en une collection de données suivante : 
\begin{itemize}
 \item[•]un semi-graphe $\lbrace\G, \VV, \EE\rbrace$;
 \item[•]pour chaque sommet $v\in \VV$, un anabélioïde connexe $\GG_v$;
 \item[•]pour chaque arête $e\in \EE$, un anabélioïde connexe $\GG_e$;
 \item[•]pour chaque branche $b$ d'une arête $e$ aboutissant à un sommet $v$, un morphisme d'anabélioïdes $b : \GG_e\rightarrow\GG_v$, c'est-à-dire un foncteur exact $b^* : \GG_v\rightarrow\GG_e$.
 \end{itemize} \
 
 Nous appellerons \emph{anabélioïdes constituants} les différents anabélioïdes connexes $\GG_e$ et $\GG_v$. Un semi-graphe d'anabélioïdes $\GG$ est dit \emph{connexe} (resp. \emph{fini}, \emph{dénombrable}, \emph{localement fini}) lorsque le semi-graphe sous-jacent $\G$ est connexe (resp. fini, dénombrable, localement fini). Nous parlerons de \emph{graphe d'anabélioïdes} lorsque $\G$ sera un graphe.\\
 
 Il y a une notion naturelle de \emph{morphisme de semi-graphes d'anabélioïdes} : un morphisme entre deux semi-graphes d'anabélioïdes $\GG$ et $\GG'$ est la donnée :
 \begin{itemize}
 \item[•] d'un morphisme entre les semi-graphes associés : $\G\rightarrow\G'$,
 \item[•] pour chaque sommet $v$ (resp. arête $e$) de $\G$ envoyé sur un sommet $v'$ (resp. $e'$) de $\G'$, des morphismes d'anabélioïdes $\GG_v\rightarrow\GG'_{v'}$ et $\GG_e\rightarrow\GG'_{e'}$,
 \item[•] un isomorphisme $\psi_b$ entre les morphismes d'anabélioïdes composés ${\GG_e\rightarrow\GG_v\rightarrow\GG_{v'}}$ et $\GG_e\rightarrow\GG_{e'}\rightarrow\GG_{v'}$ (i.e un isomorphisme entre les foncteurs associés et définis dans l'autre sens) dès que $b\in e$ est une branche aboutissant à $v$.
 \end{itemize}
\end{defi}

\begin{rem}
Les $(1)$-morphismes entre $\GG$ et $\GG'$ forment de manière naturelle une catégorie (\cite{M3}, $2.4.2$), de sorte que l'on pourra parler de la \emph{$2$-catégorie des semi-graphes d'anabélioïdes}, notée $\mathfrak{Sanab}$. Soit $\vert \mathfrak{Sanab} \vert$ la $1$-catégorie dont les objets coïncident avec ceux de $\mathfrak{Sanab}$ et dont les morphismes entre deux semi-graphes d'anabélioïdes $\GG$ et $\GG'$ correspondent aux classes d'isomorphismes des $1$-morphismes $\GG\to \GG'$ dans $\mathfrak{Sanab}$. On appellera $\vert \mathfrak{Sanab} \vert$ la \emph{catégorie grossière des semi-graphes d'anabélioïdes}.\

D'après \cite{M3} $2.4.2$, si $\varphi : \GG\to \GG'$ est un $1$-morphisme \emph{localement ouvert} entre deux semi-graphes d'anabélioïdes \emph{totalement détachés} et \emph{sommitalement minces} (ces notions seront définies dans la section \ref{section hypothèses mochizukiennes}, voir aussi \cite{M3} $2.2$ pour \emph{localement ouvert}), $\varphi$ n'a aucun $2$-automorphisme non trivial dans $\mathfrak{Sanab}$. De là, si l'on se restreint aux morphismes localement ouverts entre semi-graphes d'anabélioïdes totalement détachés et sommitalement minces, ces morphismes peuvent être considérés comme de simples morphismes de la catégorie $\vert \mathfrak{Sanab}\vert$ plutôt que comme des $1$-morphismes de la $2$-catégorie $\mathfrak{Sanab}$. De là, dans ce cadre, le produit fibré $\mathcal{Y}\times_{\mathcal{X}}\mathcal{Z}$ de deux semi-graphes d'anabélioïdes $\mathcal{Y}$ et $\mathcal{Z}$ au-dessus d'un troisième $\mathcal{X}$ est bien défini et existe dans la catégorie grossière $\vert \mathfrak{Sanab}\vert$. Hors de ce cadre, le produit fibré est à considérer au sens du $2$-produit fibré dans la $2$-catégorie $\mathfrak{Sanab}$. Remarquons néanmoins que les semi-graphes d'anabélioïdes que nous étudierons par la suite vérifieront les \emph{hypothèses mochizukiennes} (voir section \ref{section hypothèses mochizukiennes}), en particulier les semi-graphes seront totalement détachés, sommitalement minces, et les morphismes localement ouverts, ce qui nous permettra en pratique de considérer le produit fibré (ou le changement de base) dans $\vert \mathfrak{Sanab}\vert$.

\end{rem}

\begin{defi}

Soit $\GG$ un semi-graphe d'anabélioïdes non vide ayant un nombre fini de composantes connexes.\

 Définissons la catégorie des \emph{recouvrements finis de $\GG$}, notée $\BB(\GG),$ comme la catégorie dont les objets sont déterminés par la donnée de $\lbrace S_c, \varphi_b  \rbrace$ où : 

\begin{itemize}
\item[•]pour chaque sommet ou arête $c$, $S_c\in \ob(\GG_c)$,
\item[•]pour chaque branche $b$ d'une arête $e\in \EE$ aboutissant au sommet $v\in \VV$, ${\varphi _b : b^* S_{v}\iso S_{e}}$ est un isomorphisme dans $\GG_e$,
\end{itemize}
et dont les morphismes sont obtenus de manière naturelle à partir des morphismes au sein des différents anabélioïdes constituants $\GG_c$ de sorte que ces morphismes soient compatibles avec les différents $\varphi_b$.

\end{defi}

\begin{prop}[\cite{M3}]\
La catégorie $\BB(\GG)$ des recouvrements finis de $\GG$ est un anabélioïde, qui est connexe lorsque le semi-graphe d'anabélioïdes $\GG$ est lui-même connexe. 
\end{prop}

Lorsque $\GG$ est connexe, $\BB(\GG)$ étant un anabélioïde connexe, il a un groupe fondamental.

\begin{defi}
Si $\GG$ est un semi-graphe d'anabélioïdes connexe, son \emph{groupe fondamental}, noté $\pi_1(\GG)$, est défini comme égal au groupe fondamental $\pi_1(\BB(\GG))$ de l'anabélioïde connexe $\BB(\GG)$. Il est défini à unique isomorphisme extérieur près.
\end{defi}

Considérons maintenant un semi-graphe d'anabélioïdes $\GG$. Prenons $S\in \ob(\BB(\GG))$ et considérons le revêtement fini étale d'anabélioïdes $i_S : \BB(\GG)_S\rightarrow\BB(\GG)$. Alors les définitions précédentes impliquent que $\BB(\GG)_S$ est exactement la catégorie des recouvrements finis (i.e. un $\BB(\--)$) d'un certain semi-graphe d'anabélioïdes $\GG'$ équipé d'un morphisme $\GG'\rightarrow\GG$ de semi-graphes d'anabélioïdes au-dessus d'un morphisme de semi-graphes $\G'\rightarrow\G$ qui est \emph{propre}, c'est-à-dire qui préserve le nombre de branches des arêtes. Remarquons que $\GG'$ est un semi-graphe connexe dès que $\GG$ et $S$ sont connexes. \\

Plus précisément, si $v$ est par exemple un sommet de $\G$, l'ensemble des sommets $v'$ de $\G'$ au-dessus de $v$ correspond à l'ensemble des composantes connexes de $S_v$, et l'anabélioïde connexe $\GG'_{v'}$ est isomorphe à $(\GG_v)_{S_{v'}}$ où $S_{v'}$ est la composante connexe de $S_v$ correspondant à $v'$.

\begin{defi}\

\begin{enumerate}
\item Si $\GG$ est un semi-graphe d'anabélioïdes, un \emph{revêtement fini étale} de $\GG$ est un morphisme de semi-graphes d'anabélioïdes $\GG'\rightarrow\GG$ induit par le revêtement fini étale d'anabélioïdes $i_S : \BB(\GG)_S\rightarrow\BB(\GG)$ pour un certain $S\in \ob(\BB(\GG))$. 
\item Un morphisme de semi-graphes d'anabélioïdes  $\GG'\rightarrow\GG$ est dit \emph{localement trivial} (resp. \emph{localement fini étale}) si et seulement si les morphismes induits entre chaque anabélioïde constituant sont des isomorphismes (resp. sont finis étales).
\end{enumerate}

\end{defi}

\subsection{Le monde des tempéroïdes}\label{section le monde des tempéroïdes}

Nous allons ici décrire le langage des tempéroïdes introduit par Mochizuki dans \cite{M3}, qui permettra de définir le \emph{groupe fondamental tempéré d'un semi-graphe d'anabélioïdes}, analogue dans un contexte beaucoup plus catégorique du \emph{groupe fondamental tempéré} d'espaces analytiques que nous introduirons dans la suite. \\

Si $\Pi$ est un groupe topologique, notons $\BB^\tp(\Pi)$ la catégorie dont les objets sont les ensembles dénombrables discrets munis d'une action continue (i.e. à stabilisateurs ouverts) de $\Pi$, et dont les morphismes sont les morphismes de $\Pi$-ensembles. Un objet $T$ de $ \BB^\tp(\Pi)$ est connexe (au sens des catégories) si et seulement si l'action de $\Pi$ est transitive sur $T$. Soit $\mathfrak{Ens}^{\mathrm{dd}}$ la catégorie des ensembles dénombrables discrets.

\begin{defi}\
\begin{enumerate}
\item Un groupe topologique $\Pi$ est dit \emph{tempéré} lorsqu'il peut être écrit comme une limite projective de groupes topologiques discrets dénombrables. 
\item Un \emph{tempéroïde connexe} est défini comme un couple $(\TT, \Upsilon_{\TT})$ où : 
\begin{itemize}
\item[•]$\TT$ est une catégorie équivalente à $\BB^\tp(\Pi)$ pour un certain groupe tempéré $\Pi$,
\item[•]$\Upsilon_{\TT}$ est une classe d'isomorphismes de foncteurs $\TT  \to \mathfrak{Ens}^{\mathrm{dd}}$ conservatifs préservant les limites finies et les colimites dénombrables.
\end{itemize}
Les foncteurs $\TT\to \mathfrak{Ens}^\mathrm{dd}$ dont la classe d'isomorphisme est $\Upsilon_{\TT}$ sont appelés les \emph{foncteurs fibres} de $\TT$.
\item Un tempéroïde est défini comme un couple $(\TT,\Upsilon_\TT)$ où : 
\begin{itemize}
\item[•]$\TT$ est une catégorie équivalente à une catégorie produit $\ds{\prod _{j\in J} \TT_j},$ où $J$ est un ensemble dénombrable et $(\TT_j, \Upsilon_{\TT_j})$ est un tempéroïde connexe pour tout $j\in J$. 
\item[•]$\Upsilon_\TT$ est la donnée des classes d'isomorphismes des foncteurs $\TT\to \mathfrak{Ens}^\mathrm{dd}$ contenant les foncteurs se factorisant par les foncteurs fibres des $\TT_j$, c'est-à-dire les composés : $$\TT\simeq \prod _{j\in J} \TT_j \xrightarrow{\pi_{j_0}} \TT_{j_0} \xrightarrow{\beta_{j_0}} \mathfrak{Ens}^{\mathrm{dd}}$$ pour $j_0\in J$ et $\beta_{j_0}\in \Upsilon_{\TT_{j_0}}$.
\end{itemize}
Les foncteurs $\TT\to \mathfrak{Ens}^\mathrm{dd}$ dont la classe d'isomorphisme est dans $\Upsilon_{\TT}$ sont appelés les \emph{foncteurs fibres} de $\TT$.

\item Soient $(\TT_1,\Upsilon_{\TT_1})$ et $(\TT_2,\Upsilon_{T_2})$ des tempéroïdes. Un \emph{morphisme de tempéroïdes} $\varphi : (\TT_1,\Upsilon_{\TT_1})\to (\TT_2,\Upsilon_{T_2})$ est défini comme un foncteur $\varphi^* : \TT_2\rightarrow\TT_1$ dans le sens opposé qui préserve les colimites dénombrables et les limites finies, et tel que $$\Upsilon_{\TT_1}\circ \varphi^*\subseteq \Upsilon_{\TT_2},$$ c'est-à-dire qui préserve les foncteurs fibres.
\end{enumerate}
\end{defi}

\begin{rem}
Soit $J$ un ensemble dénombrable et $(\TT_j, \Upsilon_{\TT_j})$ un tempéroïde connexe, pour tout $j\in J$. Considérons $$\TT_J:=\prod_{j\in J}\TT_j.$$
Alors l'ensemble $J$ et les différents $\TT_j$ sont entièrement déterminés par $\TT_J$ tout entier. Les différents $\TT_j$ seront appelés les \emph{composantes connexes de $\TT_J$}. Ce résultat se montre comme dans \ref{unicité de la décomposition d'un objet en composantes connexes dans une multi-galoisienne } et \ref{unicité des composantes connexes d'un anabélioide}, et découle de l'unicité de la décomposition en objets connexes de l'objet initial $\mathbf{1}_{\TT_J}$ de $\TT_J$. 
\end{rem}

\begin{rem}
 Tout groupe profini est tempéré. Si $\Pi$ est un groupe tempéré, alors sa topologie  admet une base de sous-groupes ouverts distingués, et tout sous-groupe ouvert de $\Pi$ est fermé et d'indice dénombrable. Si $H\subseteq \Pi$ est un sous-groupe, alors l'ensemble $\Pi/H$ muni de l'action à gauche de $\Pi$ définit un objet de $\BB^\tp(\Pi)$ si et seulement si $H$ est ouvert dans $\Pi$. Par ailleurs $\BB^\tp(\Pi)$ est naturellement un tempéroïde connexe dont les foncteurs fibres sont définis comme l'ensemble des foncteurs isomorphes au foncteur oubli $$\beta_\Pi^\mathrm{oubli} : \BB^\tp(\Pi)\to \mathfrak{Ens}^\mathrm{dd}.$$ Soit $\Upsilon_\Pi^\mathrm{oubli}$ la classe d'isomorphisme de ces foncteurs fibres. 
\end{rem}

\bigskip
La définition que donne Mochizuki dans son article initial \cite{M3} d'un tempéroïde connexe $\TT$ ne fait pas intervenir la classe d'isomorphisme de foncteurs fibres $\Upsilon_\TT$. Néanmoins, une telle définition est insuffisante pour définir correctement le groupe fondamental d'un tempéroïde connexe. Afin de pallier cet écueil Mochizuki introduit dans les corrections \cite{M3'} la notion de \emph{Galois-dénombrable} (Galois-countable), qui pour les groupes topologiques considérés est une condition d'existence d'une base dénombrable. Lorsque tous les groupes, tempéroïdes et semi-graphes d'anabélioïdes sont Galois-dénombrables, tous les résultats structuraux $3.2$ jusqu'à $3.9$ de \cite{M3} restent vrai et le groupe fondamental tempéré est correctement défini. Néanmoins, dans notre cadre des courbes analytiques développé dans la partie \ref{section courbes analytiques}, il n'est pas clair que les groupes tempérés (au sens de \cite{And}, voir \ref{groupe fondamental tempéré}) soient Galois-dénombrables. Cette condition de Galois-dénombrabilité peut cependant être contournée par la donnée d'une classe de foncteurs fibres, et notre cadre des semi-graphes d'anabélioïdes des courbes analytiques munies d'une triangulation fera naturellement apparaître une classe de foncteurs fibres géométriques (cf. \ref{foncteurs fibre géométriques}).

\bigskip
En pensant aux tempéroïdes connexes comme à des analogues d'anabélioïdes connexes, stables par coproduits dénombrables et où les groupes prodiscret remplacent les groupes profinis, nous aimerions être capables de définir le \og groupe fondamental\fg{} d'un tempéroïde connexe. La proposition suivante le permet :

\begin{prop} [\cite{M3}, Prop. 3.2]
Soient $\Pi_1$ et $\Pi_2$ deux groupes tempérés. Alors il y a une bijection naturelle entre les \emph{classes d'isomorphismes de morphismes de tempéroïdes connexes} $(\BB^\tp(\Pi_1),\Upsilon_{\Pi_1}^\mathrm{oubli})\rightarrow (\BB^\tp(\Pi_2),\Upsilon_{\Pi_2}^\mathrm{oubli})$ et l'ensemble des \emph{homorphismes continus extérieurs} $\Pi_1\rightarrow\Pi_2$ (i.e. l'ensemble des orbites de l'action de $\Pi_2$ par conjugaison à droite sur l'ensemble $\Hom(\Pi_1\rightarrow\Pi_2)$ des homomorphismes continus de $\Pi_1$ dans $\Pi_2$).
\end{prop} 

\begin{rem}
Deux morphismes de tempéroïdes $\varphi, \psi : (\TT_1,\Upsilon_{\TT_1})\rightarrow(\TT_2,\Upsilon_{\TT_2})$ sont dits \emph{isomorphes} si les foncteurs associés $\varphi^*, \psi ^*: \TT_2\rightarrow\TT_1$ sont isomorphes, c'est-à-dire s'il existe une transformation naturelle entre $\varphi$ et $\psi$ dont tous les morphismes induits sont des isomorphismes. 
\end{rem}

\begin{defi}[Groupe fondamental tempéré d'un tempéroïde connexe]La proposition précédente permet d'affirmer que si $(\TT,\Upsilon_\TT)$ est un tempéroïde connexe, alors tous les groupes tempérés $\Pi$ tels que $(\TT,\Upsilon_\TT)$ soit isomorphe à  $(\BB^\tp(\Pi),\Upsilon_\Pi^\mathrm{oubli})$ sont isomorphes entre eux (de façon canonique modulo les automorphismes intérieurs). Cela permet ainsi de définir le \emph{groupe fondamental tempéré} de $(\TT,\Upsilon_\TT)$, noté $$\pi_1^\tp(\TT),$$ comme étant égal au groupe $\Pi$ défini à unique isomorphisme extérieur près. Il s'agit donc d'un groupe tempéré, et l'on a $\Pi\simeq \Aut(\upsilon_\TT)$ dès que $\upsilon_\TT$ est un foncteur fibre de $\TT$. 
\end{defi}

Nous allons maintenant introduire la notion de \emph{revêtement étale d'un tempéroïde}, qui se définit de manière très analogue à celle de \emph{revêtement étale fini d'un anabélioïde}.

\begin{prop}
Soit $(\TT,\Upsilon_\TT)$ un tempéroïde et $T\in \ob(\TT)$. Soit $\TT_T$ la catégorie des objets de $\TT$ au-dessus de $T$, et $\Upsilon_{\TT_T}$ l'ensemble des classes d'isomorphisme des foncteurs $\TT_T\to \mathfrak{Ens}^\mathrm{dd}$ contenant ceux qui se factorisent par les foncteurs fibres de $\TT$ le long du foncteur canonique $\TT_T\to \TT$ obtenu en oubliant la structure au-dessus de $T$. Alors $(\TT_T,\Upsilon_{\TT_T})$ est un tempéroïde dont les composantes connexes correspondent aux composantes connexes de $T$ en tant qu'objet de $\TT$. De plus, le foncteur $i_T^* : \TT\rightarrow \TT_T$ défini en prenant le produit avec $T$ commute avec les limites finies et les colimites dénombrables et définit donc un morphisme de tempéroïdes  : $$i_T : (\TT_T,\Upsilon_{\TT_T})\rightarrow (\TT,\Upsilon_\TT)$$
\end{prop}

\begin{defi}
Un morphisme de tempéroïdes connexes $\varphi : (\mathcal{U},\Upsilon_{\mathcal{U}})\rightarrow(\TT,\Upsilon_\TT)$ est un \emph{revêtement étale} s'il se factorise comme le composé d'un isomorphisme $\alpha : (\mathcal{U},\Upsilon_{\mathcal{U}})\iso(\TT_T,\Upsilon_{\TT_T})$ avec le morphisme $i_T : (\TT_T,\Upsilon_{\TT_T})\rightarrow(\TT,\Upsilon_\TT)$ où $T\in \ob(\TT)$.
\end{defi}

Si $\XX$ est une catégorie, nous noterons $\XX^\top$ la catégorie formée à partir de $\XX$ dont les objets sont tous les coproduits dénombrables \og formels \fg{}  arbitraires d'objets de $\XX$. Si $I$ et $J$ sont deux ensembles dénombrables d'indices, les morphismes sont définis de manière naturelle par : 
$$\Hom_{\XX^\top}(\coprod_{i\in I} A_i, \coprod_{j\in J} B_j) := \prod_{i\in I}\coprod_{j\in J} \Hom_\XX(A_i, B_j)$$

\begin{rem}
Soit $\XX$ est un anabélioïde connexe muni d'un point de base $\beta$ défini par un foncteur fondamental $\beta^* : \XX\to \mathfrak{Ens}^\mathrm{f}$. Ce foncteur se prolonge naturellement en un foncteur : $\beta^{*\top} : \XX^\top\to \mathfrak{Ens}^\mathrm{dd}$. Alors $\XX^\top$ est naturellement un tempéroïde connexe dont les foncteurs fibres sont définis comme les foncteurs de $\XX^\top$ à $\mathfrak{Ens}^\mathrm{dd}$ isomorphes à $\beta^{*\top}$. On appellera ce tempéroïde la \emph{tempérification} de $\XX$. En particulier, si $G$ est un groupe profini, on a une équivalence de catégories $\BB(G)^\top\simeq \BB^\tp(G)$.
\end{rem}

\begin{defi}

Soit $\GG$ un semi-graphe non vide d'anabélioïdes. \

Définissons la catégorie des \emph{recouvrements de $\GG$}, notée $$\BB^\cov(\GG),$$ comme la catégorie dont les objets sont déterminés par la donnée de $\lbrace S_c, \varphi_b  \rbrace$ où : 

\begin{itemize}
\item[•]pour chaque sommet ou arête $c$,  $S_c\in \ob(\GG_c^\top)$,
\item[•]pour chaque branche $b$ d'une arête $e\in \EE$ aboutissant au sommet $v\in \VV$, ${\varphi _b : b^* S_{v}\iso S_e}$ est un isomorphisme dans $\GG_e^\top$,
\end{itemize}
et dont les morphismes sont obtenus de manière naturelle à partir des morphismes au sein des $\GG_c^\top$ de sorte que ces morphismes soient compatibles avec les différents $\varphi_b$.

\end{defi}

\begin{rem}
Nous avons implicitement utilisés dans la définition précédente le fait que les foncteurs $b^* : \GG_{v}\rightarrow \GG_e$ se prolongent tout naturellement en des foncteurs : $b^{*\top} : \GG_{v}^\top\rightarrow \GG_e^\top$.
\end{rem}

Considérons maintenant un semi-graphe d'anabélioïdes $\GG$. Si $S\in \ob(\BB^\cov(\GG))$, les définitions précédentes impliquent que la catégorie $\BB^\cov(\GG)_S$ au-dessus de $S$ est exactement la catégorie des recouvrements (i.e. un $\BB^\cov(\--))$ d'un certain semi-graphe d'anabélioïdes $\GG'$ équipé d'un morphisme $\GG'\rightarrow\GG$ de semi-graphes d'anabélioïdes. Remarquons que $\GG'$ est connexe dès que $S$ est connexe.

\begin{defi}
Si $\GG$ est un semi-graphe d'anabélioïdes, un \emph{revêtement étale} de $\GG$ est un morphisme de semi-graphes d'anabélioïdes $\GG'\rightarrow\GG$ induit comme dans ce qui précède par un certain $S\in \ob(\BB^\cov(\GG))$.  
\end{defi}

Par analogie avec le monde des revêtements finis de semi-graphes d'anabélioïdes nous aimerions pouvoir dire que $\BB^\cov(\GG)$ est un tempéroïde connexe dès que $\GG$ est connexe. Or cela est faux généralement, et afin d'avoir des résultats similaires nous allons avoir besoin d'introduire une certaine catégorie de revêtements étales qui contiendra à la fois les revêtements finis et les revêtements étales triviaux (isomorphes à un coproduit dénombrable de copies de la base), et à laquelle il sera suffisant de se restreindre.\\

Si $\GG$ est un semi-graphe d'anabélioïdes, nous allons définir la catégorie des \emph{recouvrements tempérés} $\BB^\tp(\GG)$ de $\GG$ qui sera une sous-catégorie pleine de $\BB^\cov(\GG)$ contenant les recouvrements finis $\BB(\GG)$. On aura ainsi des injections pleines : 
$$\BB(\GG)\hookrightarrow\BB^\tp(\GG)\hookrightarrow\BB^\cov(\GG)$$

\begin{defi}[Revêtements tempérés]\
Un recouvrement étale $\lbrace S_v, \varphi_e  \rbrace\in \BB^\cov(\GG)$ est dit \emph{tempéré} si le revêtement étale de semi-graphes d'anabélioïdes $\GG'\rightarrow\GG$ auquel il donne lieu est tel qu'il existe un revêtement étale fini $\GG''\rightarrow\GG$ avec la propriété que pour chaque sommet ou arête $c$ de $\G$, la restriction de $\GG''\rightarrow\GG$ à $\GG_c$ trivialise la restriction de $\GG'\rightarrow\GG$ à $\GG_c$. Cela revient à dire que par le changement de base $\GG''_c\rightarrow\GG_c$ le revêtement $\GG'_c\rightarrow\GG_c$  est trivial, c'est-à-dire isomorphe au coproduit d'un nombre dénombrable de copies de $\GG''_c$.\

Un tel revêtement étale $\GG'\rightarrow\GG$ induit par un recouvrement tempéré sera appelé \emph{revêtement tempéré} du semi-graphe d'anabélioïdes $\GG$. La catégorie des recouvrements tempérés sera notée $$\BB^\tp(\GG).$$
\end{defi}\

Soit $\GG$ est un semi-graphe d'anabélioïdes. Si $c$ est une composante (sommet ou arête) du semi-graphe sous-jacent $\G$, tout foncteur fondamental $\beta^*_c$ de $\GG_c$ se prolonge en un foncteur $\beta_c^{*\top} : \GG_c^\top\to \mathfrak{Ens}^{\mathrm{dd}}$ conservatif et préservant les limites finies et colimites dénombrables

\begin{defi}\label{foncteurs fibre géométriques}
Un \emph{foncteur fibre géométrique} de $\BB^\tp(\GG)$ dans $\mathfrak{Ens}^\mathrm{dd}$ est défini comme un foncteur composé 
$$\BB^\tp(\GG)\to \GG_c^\top\xrightarrow{\beta_c^{*\top}}\mathfrak{Ens}^\mathrm{dd},$$ où $c$ est une composante de $\G$ et $\beta_c$ un point de base de $\GG_c$. Lorsque $\GG$ est connexe, tous les foncteurs fibres géométriques sont isomorphes.\

 Soit $\Upsilon^\mathrm{g\acute{e}o}_{\GG}$ l'ensemble des foncteurs de $\BB^\tp(\GG)$ dans $\mathfrak{Ens}^\mathrm{dd}$ isomorphes à un foncteur fibre géométrique. 
\end{defi}

\bigskip
Le prochain résultat, de \cite{M3}, nous montre que les revêtements tempérés sont suffisants dans la théorie des tempéroïdes, et assez riches pour pouvoir définir sous certaines hypothèses le groupe fondamental tempéré d'un semi-graphe d'anabélioïdes. 

\begin{prop}[Groupe tempéré d'un semi-graphe d'anabélioïdes]\label{groupe tempéré d'un semi-graphe d'anabélioïdes} \
Soient $\GG$ et $\GG'$ des semi-graphes d'anabélioïdes connexes, quasi-cohérents, totalement élevés, totalement détachés et sommitalement minces (ces notions seront définies dans la suite) . Alors : 

\begin{enumerate}
\item La catégorie $\BB^\tp(\GG)$ munie de la classe d'isomorphisme de foncteurs fibres $\Upsilon_{\GG}^\mathrm{g\acute{e}o}$ est un tempéroïde connexe, ce qui permet de définir le \emph{groupe fondamental tempéré} du semi-graphe $\GG$, noté $$\pi_1^\tp(\GG),$$ comme égal au groupe fondamental tempéré $\pi_1^\tp(\BB^\tp(\GG))$ du tempéroïde connexe $(\BB^\tp(\GG),\Upsilon_{\GG}^\mathrm{g\acute{e}o})$. On a donc une équivalence de catégories : $$\BB^\tp(\pi_1^\tp(\GG))\approx \BB^\tp(\GG).$$
\item L'injection pleine $\BB(\GG)\hookrightarrow\BB^\tp(\GG)$ induit une injection de groupes topologiques $${\pi_1^\tp(\GG)\hookrightarrow\pi_1(\GG)}.$$
\item Tout morphisme de semi-graphes d'anabélioïdes $\GG'\rightarrow\GG$ induit un morphisme de tempéroïdes connexes $\BB^\tp(\GG')\rightarrow\BB^\tp(\GG)$. Par ailleurs si $\GG'\rightarrow\GG$ est un revêtement tempéré et si $\GG$ est cohérent, alors le morphisme induit $\BB^\tp(\GG')\rightarrow\BB^\tp(\GG)$ est un revêtement étale.
\end{enumerate}

\end{prop}

\begin{rem}\label{description du pi1 temp d'un graphe d'anabélioïdes en terme de limite} 
Soit $\GG$ un semi-graphe d'anabélioïdes connexe, quasi-cohérent, totalement élevé, totalement détaché et sommitalement mince. Supposons également que chaque arête aboutit au moins à un sommet. Soit $(\GG_i\rightarrow \GG)_{i\in I}$ une famille de revêtements étales finis galoisiens de $\GG$ qui est cofinale (c'est-à-dire telle que si pour chaque $i\in I $ $\GG_i\rightarrow \GG$ est induit par un certain $S_i\in \ob(\BB(\GG))$, alors la famille $(S_i)_{i\in I}$ est cofinale dans $\BB(\GG)$). Nous savons alors que le groupe profini $\pi_1(\GG)$ peut être construit comme la limite projective : 
$$\pi_1(\GG)= \varprojlim_{i\in I} \Gal(\GG_i/\GG). $$
Pour tout $i\in I$ notons $\widetilde{\GG}_i\rightarrow \GG_i$ le \emph{\og revêtement universel \fg{} } de $\GG_i$, i.e. le revêtement étale localement trivial de $\GG_i$ dont le semi-graphe sous-jacent correspond au revêtement universel du semi-graphe $\G_i$ correspondant à $\GG_i$. Le groupe $T_i:=\Aut_{\GG_i}(\widetilde{\GG}_i)$ est alors un groupe libre discret. Par ailleurs $\widetilde{\GG}_i\rightarrow \GG$ est un revêtement tempéré galoisien, dont le groupe de Galois s'insère dans une suite exacte : 
$$\{1\} \longrightarrow T_i\longrightarrow \Gal(\widetilde{\GG}_i/\GG)\longrightarrow \Gal(\GG_i/\GG)\longrightarrow\{1\}.$$

La famille $(\Gal(\widetilde{\GG}_i/\GG))_{i\in I}$ est un système projectif, et le groupe tempéré de $\GG$ peut être décrit de la manière suivante : 

$$\pi_1^\tp(\GG)= \varprojlim_{i\in I} \Gal(\widetilde{\GG}_i/\GG). $$

\end{rem}

\subsection{Hypothèses mochizukiennes}\label{section hypothèses mochizukiennes}

Avant d'énoncer le résultat fondamental sur les semi-graphes d'anabélïoides dont nous nous servirons dans la suite, nous allons tout d'abord avoir besoin de définir certaines notions que l'on qualifiera de \emph{mochizukiennes}. \\

Si $\GG$ est un semi-graphe d'anabélïoides, notons $\pi _v$ et $\pi_b$ (ou $\pi_e$) les groupes fondamentaux des anabélioïdes connexes $\GG_v$ et $\GG_e$ dont on choisit des points de bases, pour chaque sommet $v$ et chaque branche $b$ d'une arête $e$ aboutissant à $v$. La branche $b$ définit un morphisme de groupes (extérieur) $b_* : \pi_b \rightarrow \pi_v$. On dira que  $\GG$ est de \emph{type injectif} lorsque tous les morphismes de ce type sont injectifs. Si $\GG$ est de type injectif, nous continuerons à  noter $\pi_b$ l'image de $\pi_b$ dans $\pi_v$ par $b_*$, celle-ci étant définie à conjugaison près dans $\pi_v$.

\begin{defi}
Soit $\GG$ un semi-graphe d'anabélioïdes de type injectif. Un \emph{approximateur de $\GG$} sera défini comme un morphisme de semi-graphes d'anabélioïdes ${\GG\rightarrow \GG'}$ au-dessus d'un même semi-graphe $\G$ tel que $\GG'$ est de type injectif et : 
\begin{enumerate}
\item il existe un entier $N\geqslant 1$ tel que pour tout sommet $v$ de $\G$ le groupe $\pi'_v$ soit fini de cardinal $\leqslant N$;
\item si $c$ est une arête ou un sommet de $\G$, alors le morphisme de groupes $\pi_c \rightarrow \pi_c'$ induit par le morphisme d'anabélioïdes $\GG_c\rightarrow \GG'_c$ est surjectif.
\end{enumerate}
\end{defi}

\begin{defi}
Un groupe profini $G$ est dit \emph{mince} lorsque le centralisateur $Z_G(H)$ de tout sous-groupe ouvert $H$ de $G$ est trivial.
\end{defi}

\begin{rem}
Un groupe profini $G$ est mince si et seulement si le centre de tout sous-groupe ouvert est trivial. En effet la nécessité de cette condition est évidente, tandis que l'autre sens se montre comme suit : si $H\subseteq G$ est un sous-groupe ouvert et $h\in Z_G(H)$, alors le sous-groupe $H'$ engendré par $H$ et $h$ est ouvert et $h$ est dans son centre, donc trivial par hypothèse.
\end{rem}

\begin{defi}\label{contraintes locales mochizukiennes}
Soit $\GG$ un semi-graphe d'anabélioïdes de type injectif de semi-graphe sous-jacent $\G$.
\begin{enumerate}
\item $\GG$ est dit \emph{totalement détaché} lorsque, pour tout sommet $v$ auquel une branche $b$ aboutit, $\pi_v$ est d'ordre infini et pour tout $g\in \pi_v$: 
\begin{itemize}
\item[i)] l'intersection $\pi_b\cap g.\pi_b.g^{-1}$ est triviale dans $\pi_v$ si $g \notin \pi_b$
\item[ii)] l'intersection $\pi_b\cap g.\pi_{b'}.g^{-1}$ est triviale dans $\pi_v$ si $b'$ est une branche différente de $b$ aboutissant à $v$.
\end{itemize}
\item $\GG$ est dit \emph{totalement élevé} si pour tout sommet $v$ de $\G$ et pour tout entier $M\geqslant 1$ il existe un approximateur ${\GG\rightarrow \GG'}$ et un sous-groupe $N_M \subseteq \pi'_v$ d'ordre $\geqslant M$ tel que dans $\pi'_v$ l'intersection $N_M\cap \pi'_b$ soit triviale pour toute branche $b$ aboutissant à $v$.
\item $\GG$ est dit \emph{sommitalement mince} si pour tout sommet $v$ de $\G$ l'anabélioïde connexe $\GG_v$ est \emph{mince}, c'est-à-dire si $\pi_v$ est un groupe mince.
\end{enumerate}
\end{defi}

Les définitions précédentes étant de nature locale, nous allons donner une dernière définition de nature plus globale, celle que \emph{quasi-cohérence}.

\begin{defi}
Un semi-graphe d'anabélioïdes $\GG$ de type injectif sera dit \emph{quasi-cohérent} lorsque pour tout entier $r\geqslant 1$ et toute collection de revêtements finis étales $(\HHH_c\rightarrow \GG_c)_c$ de degrés $\leqslant r$ (où $c$ décrit toutes les composantes du graphe $\G$, ie. l'ensemble des arêtes et des sommets), il existe un approximateur ${\GG\rightarrow \GG'}$ tel que, pour toute composante $c$ de $\G$, le tiré en arrière à $\GG_c$ du \og revêtement universel \fg{} $\HHH'_c\rightarrow \GG'_c$ de $\GG'_c$ (ie. le revêtement fini étale déterminé par le sous-groupe trivial de $\GG'_c$) trivialise $\HHH_c\rightarrow \GG_c$.

$\GG$ sera dit \emph{cohérent} s'il est \emph{quasi-cohérent} et que, pour chaque composante $c$ de $\G$, le groupe profini $\pi_c$ est topologiquement finiment engendré.
\end{defi}

\subsection{Reconstruction du semi-graphe sous-jacent d'un tempéroïde}\label{section reconstruction du semi-graphe sous-jacent d'un tempéroïde}

 Nous allons énoncer un résultat de Mochizuki, fondamental dans la suite de ce texte, qui justifie l'usage que l'on fait des semi-graphes d'anabélioïdes. Il faut le voir comme une version faible de \cite{M3} , Corollary 3.9. \\
Mochizuki définit les morphismes \emph{quasi-géométriques} de tempéroïdes puis énonce un résultat général important qui donne sous certaines hypothèses une correspondance entre les morphismes localement ouverts de semi-graphes d'anabéloïdes  et les morphismes quasi-géométriques des tempéroïdes correspondants. Il n'est cependant pas nécessaire ici de définir les morphismes quasi-géométriques car nous utiliserons une version plus faible de son résultat. En effet, en remarquant que tout isomorphisme de tempéroïdes est quasi-géométrique, on obtient le résultat suivant : 

\begin{thm}[Reconstruction du semi-graphe d'anabélioïdes sous-jacent, \cite{M3}  (Corollary 3.9)]\ \label{reconsmochi} 

Soient $\GG$ et $\HHH$ des graphes d'anabélioïdes connexes, de type injectif, quasi-cohérents, totalement élevés, totalement détachés et sommitalement minces. Alors en appliquant $"\BB^\tp(\--)"$ on obtient une bijection naturelle entre les \emph{isomorphismes de semi-graphes d'anabélioïdes} $\GG\iso\HHH$ et les \emph{isomorphismes de tempéroïdes connexes} $$(\BB^\tp(\GG),\Upsilon^\mathrm{g\acute{e}o}_{\GG})\iso(\BB^\tp(\HHH),\Upsilon^\mathrm{g\acute{e}o}_{\HHH}).$$
\end{thm}

\section{Courbes analytiques et semi-graphes d'anabélioïdes associés}\label{section courbes analytiques}

\subsection{Courbes $k$-analytiques}\label{courbes analytiques}

Soit $k$ un corps non archimédien complet \emph{algébriquement clos}. Le corps résiduel $\widetilde{k}$ peut être \emph{a priori} de caractéristique nulle. Désignons par $p$ l'exposant caractéristique de $\widetilde{k}$, i.e. sa caractéristique si elle est non nulle et $1$ sinon. Nous travaillerons dans la suite avec la notion d'espace analytique \emph{au sens de Berkovich} introduite dans le livre fondateur \cite{Ber1}. En particulier un espace de Berkovich sera \emph{localement compact} et \emph{localement connexe par arcs}. \\

Nous désignerons par \emph{courbe $k$-analytique} tout espace $k$-analytique séparé dont toutes les composantes connexes sont de dimension $k$-analytique égale à $1$. D'après \cite{Ber4} nous savons que les courbes $k$-analytiques ont localement l'allure d'\emph{ \og arbres réel \fg{} } et sont \emph{localement simplement connexes}, ce qui nous permettra d'y appliquer la théorie du groupe fondamental topologique et du revêtement universel. Si $X$ est une courbe $k$-analytique, en reprenant les définitions de \cite{Ber1} ou \cite{Duc}, nous parlerons des points de type $1$, $2$, $3$ ou $4$. Pour $i\in \lbrace 1, 2, 3, 4\rbrace$ notons $X_{[i]}$ l'ensemble des points de $x$ de type $i$. Alors $X_{[1]}, X_{[2]}, X_{[3]}$ et $X_{[4]}$ forment une partition de $X$. Notons également $X_{[0]}$ l'ensemble des points \emph{rigides} de $X$. Nous avons $X_{[1]}=X_{[0]}$ puisque $k$ est supposé algébriquement clos. 
\\

\begin{enumerate}
\item\textbf{Droites $k$-analytiques affine et projective} \\

\begin{itemize}
\item[•]La \emph{droite affine $k$-analytique}, $\A_k^{1, \an}$, est la courbe $k$-analytique connexe, sans bord et lisse dont les points sont les semi-normes multiplicatives sur l'anneau de polynômes $k[T]$ qui étendent la norme de $k$. Elle est recouverte par les domaines affinoïdes de la forme : $$\mathrm{B}(0, r)= \mathcal{M}\left(k\lbrace r^{-1}T\rbrace \right)   $$ définis par la condition $\vert T \vert \leqslant r$ (voir \cite{Ber1} pour la définition de \emph{domaine affinoïde}, de $k\lbrace r^{-1}T\rbrace$ et du spectre analytique $\mathcal{M}(\mathcal{A})$ d'une $k$-algèbre de Banach $\mathcal{A}$).\\

\item[•]La \emph{droite projective $k$-analytique}, $\P_k^{1, \an}$, est la courbe $k$-analytique propre (compacte et sans bord), connexe et lisse dont les points sont les semi-normes multiplicatives $x$ sur l'anneau des polynômes $k[T_0, T_1]$ vérifiant $\vert T_i(x)\vert \neq 0$ pour au moins un $i\in \lbrace 0,1 \rbrace$, modulo la relation d'équivalence qui identifie deux semi-normes $x_1$ et $x_2$ si et seulement s'il existe $\delta\in \R_+^*$ tel que l'on ait $\vert P(x_1)\vert = \delta^n\vert P(x_2)\vert $ pour tout polynôme homogène $P\in k[T_0, T_1]$ de degré $n$.  \

Soit $\lbrace \infty \rbrace$ le point de $\P_k^{1, \an}$ défini par la relation : $\vert T_0(x)\vert =0$. On a un isomorphisme d'espaces $k$-analytiques $$\rho : \P_k^{1, \an}\setminus \lbrace \infty \rbrace \iso \A_k^{1, \an}$$ donné par : $\vert T(\rho(x)) \vert= \vert T_1(x) \vert / \vert T_2(x)\vert$ pour tout $x\in \P_k^{1, \an}\setminus \lbrace \infty \rbrace$.\

En reprenant les notations de \cite{Duc}, pour tout couple $(a, r)\in k\times \R_+$, notons $\eta_{a, r}$ le point de $\P_k^{1, \an}\setminus \lbrace \infty \rbrace$ (identifié à $\A_k^{1, \an}$) défini par : $$\vert P (\eta_{a, r}) \vert=\mathrm{max}_{0\leqslant i \leqslant r}(\vert \alpha_i \vert r^i) \;\;\text{dès que} \; \;P=\sum_{i=0}^r \alpha_i (T-a)^i$$

Comme $k$ est algébriquement clos, points rigides et points de type $1$ coïncident, et correspondent aux « évaluations » en $a\in k$, $\eta_{a, 0}$, données par $\vert P(\eta_{a, 0})\vert= \vert P(a) \vert_{k}$. Si $r\neq 0$, $\eta_{a, r}$ est de type $2$ si et seulement si $r\in \vert k^*\vert$, dans le cas contraire il est de type $3$. Il existe également des points de type $4$ dès que $k$ n'est pas sphériquement clos, ce qui est le cas par exemple pour $k=\C_p$.
\end{itemize}\

Il existe une métrique naturelle sur l'ensemble des points de type $2$ et $3$ de $\P_k^{1, \an}$, invariante sous le groupe des automorphismes ($k$-analytiques) de $\P_k^{1, \an}$. Elle est définie par la formule suivante : 
$$
d(\eta_{a,r}, \eta_{a',r'})   = \left\{
    \begin{array}{ll}
        \log_p\left(\frac{\vert a-a' \vert}{r} \right)+ \log_p\left(\frac{\vert a-a' \vert}{r'} \right)& \hbox{si } \vert a-a' \vert \geqslant \max(r,r'); \\
        \vert \log_p\left( \frac{r}{r'}\right) \vert & \hbox{si } \vert a-a' \vert \leqslant \max(r,r')
    \end{array}
\right.
$$

\item\textbf{Disques et couronnes $k$-analytiques} \\

\end{enumerate}

Si $I=[a,b]$ est un intervalle compact de $\R_+^*$, notons $\CC(I)$ la courbe $k$-analytique associée à l'algèbre affinoïde $$k\lbrace b^{-1}T, aU\rbrace / (TU = 1).$$
Si $I\subset J$ sont des intervalles compacts de  $\R_+^*$, il y a un morphisme naturel $\CC(I)\rightarrow \CC(J)$ qui identifie $\CC(I)$ à un domaine analytique de $\CC(J)$. Si $I$ est un intervalle arbitraire de $\R_+^*$, on définit $$\CC(I) = \varinjlim_{J\subset I} \CC(J)\subset \G_m^\an$$ où $J$ décrit les intervalles compacts de $\R_+^*$. Il aurait été possible et équivalent de définir $\CC(I)$ comme le domaine analytique de $\A_k^{1, \an}$ défini par la condition $\mid T\mid \in I$. Une \emph{couronne} désignera une courbe $k$-analytique isomorphe à $\CC(I)$ pour un certain intervalle $I$ de $\R_+^*$. Les couronnes sont en particulier des courbes quasi-lisses (nous utilisons ici la terminologie de \cite{Duc18}), lisses dans le cas ouvert.\\

\begin{defi}[Module d'une couronne $k$-analytique]
\item Si $I$ est un intervalle de $\R_+^*$, le \emph{module} de la couronne $\CC(I)$ est défini comme la quantité : $$\Mod(\CC(I))=\frac{\sup I}{\inf I}.$$
Le \emph{module} d'une couronne $k$-analytique $\CC$, noté $\Mod(\CC)$, est défini comme le module de $\CC(I)$ pour tout intervalle $I$ de $\R_+^*$ tel que l'on ait un isomorphisme $\CC \simeq \CC(I)$. Par cohérence avec la définition donnée de la métrique sur $\P^{1,\an}_{k,[2,3]} $, la \emph{longueur} d'une couronne $\CC$ sera définie comme : $\ell(\CC)=\log_p\left(\mathrm{Mod}(\CC) \right)$. 
\item Si $r, r'\in \R_+^*$ et $a, a'\in k$ vérifient  $\vert a-a' \vert \leqslant \max(r,r')$, alors $d(\eta_{a,r}, \eta_{a',r'})=  \vert \log_p\left( \frac{r}{r'}\right) \vert $ correspond au module de l'unique couronne ouverte dont le bord est constitué des deux points $\eta_{a,r}$ et $\eta_{a',r'}=\eta_{a,r'}$.
\end{defi}\

\begin{rem}
Si $I$ et $J$ sont deux intervalles de $\R_+^*$, les couronnes $\CC(I)$ et $\CC(J)$ sont isomorphes si et seulement si $J\in \sqrt{\vert k^\times \vert}I^{\pm1}$. Par conséquent la définition ci-dessus du module d'une couronne ne dépend pas \textit{a posteriori} du choix de l'intervalle $I$. 
\end{rem}\

Nous appellerons \emph{disque} toute courbe $k$-analytique isomorphe au domaine analytique de $\A_k^{1, \an}$ définit par les conditions $\mid T\mid < r$ ou $\mid T\mid \leqslant r$ pour un certain $r\in \R_+^*$. Les disques sont en particulier des courbes quasi-lisses, lisses dans le cas ouvert. \\

\begin{defi}[Revêtement étale modéré d'un espace analytique et groupe fondamental modéré] \

\begin{itemize}
\item[•] 
Soit $\varphi : Y\rightarrow X$ un morphisme étale entre deux espaces analytiques. Si $y\in Y$, $\varphi$ est dit \emph{modéré} en $y$ si l'entier ${[L : \HC(\varphi(y))]}$ est premier à $p$, où $L$ est une clôture galoisienne de $\HC(y)$. On dit que $\varphi$ est \emph{modéré} s'il l'est en tout point $y\in Y$.

Un tel morphisme modéré est un \emph{revêtement modéré de X} lorsque $X$ est recouvert par des ouverts $U$ avec $\varphi^{-1}(U)=\coprod Y_j$ et chaque $Y_j\to U$ est fini étale (voir la notion de \emph{covering space} de \emph{\cite{DJg}}). 

Nous noterons $\Covt(X)$ la catégorie des revêtements modérés finis de $X$. 
\item[•]Lorsque $X$ est connexe $\Covt(X)$ est une catégorie galoisienne (et a fortiori un anabélioïde connexe) dont nous noterons $\pi_1^\ttt(X)$ le groupe fondamental et qui sera appelé le \emph{groupe fondamental modéré de $X$}. En particulier la catégorie $\Covt(X)$ sera équivalente à la catégorie $\BB(\pi_1^\ttt(X))$ des ensembles finis munis d'une action continue de $\pi_1^\ttt(X)$.  

\end{itemize}
\end{defi}

\begin{rem}
Cette définition de morphisme modéré est légèrement différente de la notion de \og \emph{tame morphism} \fg{}  que Berkovich introduit dans \cite{Ber2} en imposant seulement que $p$ soit premier avec ${[\HC(y) : \HC(\varphi(y))]}$. Notre nouvelle définition a l'avantage d'être plus maniable. Attention à ne pas penser qu'être modéré pour $\varphi$ en $y\in Y$ signifierait que l'extension $\HC(y)/\HC(\varphi(y))$ devrait être modérément ramifiée. Une telle définition serait moins restrictive (puisque $\HC(y)/\HC(\varphi(y))$ est automatiquement modérément ramifiée dès que ${[\HC(y) : \HC(\varphi(y))]}$ est premier avec $p$, ce qui est en particulier le cas lorsque $p$ est premier avec ${[L : \HC(\varphi(y))]}$). 
\end{rem}\

Si $x$ est un point d'une courbe analytique $X$, l'ensemble des \emph{branches de $X$ en $x$} est défini par : 
$$\br(X,x) = \varprojlim_U \pi_0(U\setminus \lbrace x\rbrace)$$ où $U$ décrit l'ensemble des voisinages ouverts de $x$. Si $V$ est un voisinage ouvert de $x$, $b(V)$ est défini comme l'image de $b$ dans $\pi_0(V\setminus \lbrace x\rbrace)$. Une \emph{section} de $b$ est un $b(V)$ pour un certain voisinage ouvert $V$ de $x$.\\

\begin{defi}[Degré d'un morphisme fini et plat]\

Soit $\varphi : Y\rightarrow X$ un morphisme \emph{fini} et \emph{plat} entre espaces $k$-analytiques, et $x\in X$. \

\begin{itemize}
\item[•] Si $U$ est un domaine affinoïde de $X$ contenant $x$, son image réciproque $\varphi^{-1}(U)$ est un domaine affinoïde de $Y$ dont l'algèbre des fonctions est un $\mathscr{O}_U(U)$-module fini et localement libre; si $U$ est connexe le rang de ce module est bien défini et ne dépend pas du choix de $U$. On l'appellera le \emph{degré de $\varphi$ au-dessus de $x$}, et on le notera $\deg_x\varphi$. La fonction $x\rightarrow \deg_x\varphi$ est localement constante sur $X$, ce qui permet sans ambiguïté de définir $\deg \varphi$ dès que $X$ est supposé \emph{connexe} et non vide.

\item[•] Supposons que $X$ et $Y$ sont des courbes $k$-analytiques, et soit $a\in \br(X,x)$. Comme $\varphi$ induit une application naturelle $\varphi : \br(Y,y)\rightarrow \br(X,x)$, il est loisible de considérer l'ensemble $\varphi^{-1}(a)\subset \br(Y,y)$. Si $b\in \varphi^{-1}(a)$, considérons un voisinage ouvert $V$ de $x$ dans $X$ qui est un arbre tel que $\varphi^{-1}(V)$ soit une réunion disjointe (finie) d'arbres qui sépare les antécédents de $x$. Soit $W$ la composante connexe de $\varphi^{-1}(a(V))$ correspondant à $b$. Le morphisme $W\rightarrow a(V)$ est alors fini et plat, et son degré ne dépend pas du choix de $U$. On le note alors : $\deg(b\rightarrow a)$.
\end{itemize}
\end{defi}\

\begin{rem}
Lorsque $X$ est un \emph{bon espace} (i.e. que chaque point de $X$ possède un voisinnage $k$-affinoïde), alors $Y$ l'est aussi, et dans ce cas $\deg_x\varphi$ n'est rien d'autre que le rang du $\mathscr{O}_{X,x}$-module libre $\displaystyle{\prod_{y\in \varphi^{-1}(x)}\mathscr{O}_{Y,y}}$ . C'est le cas en particulier dès que $X$ est une courbe $k$-analytique en vertu de la proposition $3.3.7$ de \cite{Duc} selon laquelle toute courbe $k$-analytique est un bon espace.
\end{rem}\

Un morphisme $\varphi : Y\rightarrow X$ entre deux espaces $k$-analytiques est dit \emph{totalement déployé} en $x\in \im(\varphi)$ si pour tout $y\in \varphi^{-1}(\lbrace x\rbrace)$, le morphisme $\HC(x)\rightarrow\HC(y)$ est un isomorphisme. Lorsque $\varphi$ est de degré $n$, $\varphi$ est totalement déployé en $x$ si et seulement si la fibre $\varphi^{-1}(\lbrace x\rbrace)$ est de cardinal $n$, ce qui revient à dire que $\varphi$ est localement au voisinage de $x$ un revêtement topologique (cf. \cite{And}, III.$1.2.1$).

\subsection{Squelette, triangulation et noeuds d'une courbe $k$-analytique}\label{STN} \

\begin{defi} Soit $X$ une courbe $k$-analytique quasi-lisse :
\begin{enumerate}
\item Le \emph{squelette analytique} de $X$, noté $S^\an(X)$, est défini comme l'ensemble des points de $X$ qui n'appartiennent à aucun disque analytique ouvert. 
\item Un voisinage \emph{totalement découpé autour de} $x\in X_{[2,3]}$ est un voisinage ouvert $U$ de $x$ tel que l'application $\br(X,x)\rightarrow \pi_0(U\setminus \lbrace x\rbrace)$ soit bijective et que toute composante connexe de $U\setminus \lbrace x\rbrace$ soit une couronne ou un disque. 
\end{enumerate}

\end{defi}

Si $\pi_0(S^\an)\rightarrow \pi_0(X)$ est surjective, alors d'après les théorèmes $1.6.13$ et $5.1.11$ de \cite{Duc} nous savons qu'il existe une rétraction par déformation $r_X : X\rightarrow S^\an(X)$. En particulier $X$ a le type d'homotopie de $S^\an(X)$.\\

\begin{rem}  
Un point $x$ de type $2$ ou $3$ d'une courbe $k$-analytique quasi-lisse $X$ admet toujours un voisinage totalement découpé :
\begin{itemize}
\item[•]Si $x\in X_{[3]}$, $x$ admet un voisinage ouvert totalement découpé $Z$ qui est une couronne telle que $x\in S^\an(Z)$ et dont l'adhérence dans $X$ est un arbre compact. En outre $Z$ est de type $]*,*[ $ si $x\notin \partial^\an X$, de type $]*,*]$ (avec $\partial^\an Z=\lbrace x\rbrace$) si $x\in \partial^\an X$ et $x$ n'est  pas un point isolé, et $Z=\lbrace x \rbrace$ lorsque $x$ est un point isolé de $X$ (cf. \cite{Duc}, Théorème $4.3.5$).
\item[•]Si $x\in X_{[2]}$, un voisinage $U$ de $x$ totalement découpé peut être choisi tel que son adhérence est un arbre compact et que $U\setminus \lbrace x\rbrace$ est réunion disjointe de disques et d'un nombre \emph{fini} de couronnes. Cela peut être vu comme une conséquence observée par Berkovich du théorème de réduction semi-stable. Néanmoins, nous nous servirons plutôt de \cite{Duc} où l'auteur propose un cheminement \og dans l'autre sens \fg{} , c'est-à-dire en partant de l'étude locale pour arriver au théorème de réduction stable. La remarque ci-dessus découle alors de  \cite{Duc}, Théorème 4.5.4.
\end{itemize}
\end{rem} \

\begin{rem}\label{distinction graphe/semi-graphe}
Dans \cite{Duc} l'auteur utilise le terme de \og graphe \fg{} même lorsqu'une arête est \og ouverte \fg{} au sens où elle ne contient qu'une branche. Néanmoins, dans ce texte, afin d'être cohérent avec le langage de Mochizuki, nous différencierons un graphe (toutes les arêtes ont $2$ branches) d'un semi-graphe (les arêtes peuvent avoir également $0$ ou $1$ branche). D'après \cite{Duc}, Théorème 5.1.11, le squelette analytique $S^\an(X)$ d'une courbe analytique quasi-lisse $X$ est un semi-graphe localement fini contenu dans $X_{[2,3]}$ et qui contient le bord analytique de $X$. Si $S^\an(X)$ rencontre chaque composante connexe de $X$, alors c'est un semi-graphe \emph{analytiquement admissible} dans le sens de \cite{Duc}, c'est-à-dire que c'est un sous-semi-graphe fermé de $X$ dont toute composante connexe du complémentaire est un disque relativement compact dans $X$.
\end{rem}

\begin{defi}[Squelette tronqué]
Soit $X$ est une courbe $k$-analytique quasi-lisse. Nous savons que $S^\an(X)$ est un semi-graphe (cf. remarque \ref{distinction graphe/semi-graphe} ci-dessus). Notons $S^\an(X)^\natural$, appelé le \emph{squelette analytique tronqué} de $X$, le plus grand sous-graphe de $S^\an(X)$. Il est obtenu à partir de $S^\an(X)$ en retirant les arêtes qui ont strictement moins de deux branches aboutissant à des sommets. 
\end{defi}

Une interprétation du squelette tronqué en terme de composantes relativement connexe sera donnée dans la suite de ce texte. Cette définition est cohérente avec celle de \emph{squelette tronqué d'une triangulation} que nous donnons en \ref{définition triangulation généralisée et des squelettes}.

\begin{defi}\label{définition d'une sous-couronne} 
Si $X$ est une couronne $k$-analytique, un domaine analytique $Y$ de $X$ est une \emph{sous-couronne} de $X$ s'il peut être défini comme l'image réciproque selon la rétraction $r_X : X\rightarrow S^\an(X)$ d'un sous-intervalle non vide $J$ de $S^\an(X)$. On a dans ce cas $S^\an(Y)=J$.
\end{defi}

\begin{rem}\label{chapeautée d'une couronne} 
Il y a un homéomorphisme naturel $S^\an(\CC(I))\approx I$ dès que $I$ est un intervalle de $\R_+^*$. En particulier si $\CC$ est une couronne $S^\an(\CC)$ admet une compactification topologique canonique $\widehat{S}^\an(\CC)$ dont les éléments n'appartenant pas à $\CC$ sont appelés les \emph{bouts de $\CC$}.
\end{rem}\

\begin{defi}[Triangulation]\label{définition du squelette d'une triangulation} 
Soit $X$ une courbe $k$-analytique quasi-lisse. Nous appellerons \emph{triangulation de $X$} tout sous-ensemble $S\subset X_{[2,3]}$ discret et fermé dans $X$ tel que toute composante connexe de $X\setminus S$ soit un disque relativement compact ou une couronne $k$-analytique.

\end{defi}\

\begin{rem}
Attention, notre définition de \emph{triangulation} diffère de celle que Ducros donne dans \cite{Duc} dans la mesure où celui-ci exige que les composantes connexes de $X\setminus S$ soient \emph{relativement compactes}. L'écart que nous prenons vis-à-vis de cette définition en autorisant les couronnes qui apparaissent parmi ces composantes à ne pas être relativement compactes nous permettra notamment de traiter le cas des courbes non compactes. Ducros autorise également des points rigides dans sa définition de triangulation, tandis que nous nous restreignons ici aux triangulations constituées uniquement de points de type $2$ ou $3$. Cette restriction permettra une bonne description des composantes sommitales du graphe d'anabélioïdes que l'on associera à une triangulation (voir \ref{cmportement des composantes sommitales}). Remarquons que dans le cas où $X$ est une courbe $k$-analytique \emph{compacte}, squelette et squelette tronqué coïncident.
\end{rem}

\begin{rem}
Si $X$ est \emph{strictement $k$-analytique} et si la valuation sur $k$ n'est pas triviale, alors $X$ possède au moins une triangulation constituée uniquement de points de type $2$ (\cite{Duc}, Théorème 5.1.14).
\end{rem}\

\begin{defi} [Nœuds d'un sous-semi-graphe]\label{noeuds}
Soit $\Lambda$ un sous-semi-graphe localement fini de $X$. Si $x\in \Lambda$, le \emph{genre de $x$}, noté $g(x)$, est défini comme égal à $0$ si $x$ est de type $1, 3$ ou $4$, et comme étant le genre de la courbe résiduelle $\mathscr{C}_x$ (que nous définirons dans la section \ref{courbe résiduelle}) si $x$ est de type $2$. En reprenant les notations de \cite{Duc} nous dirons qu'un point $x\in \Lambda$ est un \emph{nœud de $\Lambda$} s'il satisfait l'une au moins des conditions suivantes : 
 \begin{enumerate}
 \item $x$ est un sommet topologique de $\Lambda$ de valence $\geqslant 3$ 
 \item $x\in \partial^\an X$
 \item $g(x)>0$ 
 \end{enumerate}
 \end{defi}\

\begin{rem}
Nous savons d'après \cite{Duc} ($4.2.11.2$) que si $x\in X_{[3]}$, il y a au plus deux branches de $X$ issues de $x$, avec égalité si et seulement si $x$ est intérieur. Il est possible également que $\br(X,x)$ soit vide, cela se produit si et seulement si $x$ est un point isolé de $X$. Si $x\in \Lambda$, on en déduit que $x$ est un nœud de $\Lambda$ si et seulement si $x\in \partial^\an X$.
\end{rem}\
 
 \begin{rem}\label{remarque sur l'ensemble des noeuds qui peut ne pas être une triangulation} 
 Soit $X$ une courbe $k$-analytique connexe, quasi-lisse, non vide et $S^\an(X)$ son squelette dont on note $\Sigma_X$ l'ensemble des nœuds. Lorsque $\Sigma_X$ est non vide, les composantes connexes relativement compactes de $X\setminus \Sigma_X$ sont des couronnes ou des disques (ce n'est plus le cas cependant lorsque $\Sigma_X$ est vide, par exemple pour une courbe de Tate). Si les composantes connexes non relativement compactes de $X\setminus \Sigma_X$ sont des couronnes, alors $\Sigma_X$ est une triangulation de $X$. En revanche il se peut que certaines composantes non relativement compactes ne soient pas des couronnes, auquel cas $\Sigma_X$ ne constitue pas une triangulation de $X$. Cela se produit typiquement lorsqu'on enlève un point $x$ de type $4$ à une courbe $k$-analytique $Y$ quasi-lisse de squelette analytique non vide. Si $r$ désigne la rétraction canonique de $Y$ sur $S^\an(Y)$, enlever $x$ revient à rajouter l'arête ouverte $]x, r(x)]$ : $S^\an(Y\setminus \lbrace x\rbrace)=S^\an(Y)\cup \; ]x, r(x)]$. Or $r^{-1}(]x, r(x)])$ n'a aucune raison d'être une couronne. Cela nous amène à définir la notion de \emph{triangulation généralisée} : 
 \end{rem}\
 
 \begin{defi}[Triangulation généralisée et squelettes associés] \label{définition triangulation généralisée et des squelettes}\
 
\begin{enumerate}
\item Soit $X$ une courbe $k$-analytique quasi-lisse. Nous appellerons \emph{triangulation généralisée de $X$} tout sous-ensemble $S\subset X_{[2,3]}$ discret et fermé dans $X$ tel que toute composante connexe relativement compacte de $X\setminus S$ soit un disque ou une couronne et tel que toute composante connexe non relativement compacte $\mathcal{C}$ de $X\setminus S$ s'écrive $r_X^{-1}(S^\an(\mathcal{C}))$ où le squelette analytique $S^\an(\mathcal{C})$ est un intervalle tracé sur $X_{[2,3]}$ constitué de points intérieurs de genre nul. 
\item Si $S$ est une triangulation généralisée de $X$, la réunion $\Gamma^\natural$ de $S$ et des squelettes des composantes connexes de $X\setminus S$ qui sont des couronnes \emph{relativement compactes} est un sous-graphe localement fini et analytiquement admissible de $X$ tracé sur $X_{[2,3]}$. Le graphe $\Gamma^\natural$ sera alors appelé le \emph{squelette tronqué} de la triangulation $S$.
\item Si $S$ est une triangulation généralisée de $X$, la réunion $\Gamma$ de $S$ et des squelettes des composantes connexes de $X\setminus S$ qui ne sont pas des disques est un semi-graphe localement fini tracé sur $X_{[2,3]}$. Le semi-graphe $\Gamma$ sera alors appelé le \emph{squelette} de la triangulation généralisée $S$. Il est obtenu à partir du squelette tronqué $\Gamma^\natural$, qui est un graphe, en ajoutant des arêtes ouvertes qui sont exactement les squelettes analytiques des composantes connexes non relativement compactes de $X\setminus S$.  
\end{enumerate}
 \end{defi}\

 \begin{rem}
 Toute triangulation de $X$ est un cas particulier de triangulation généralisée, et ces deux notions coïncident dans le cas où $X$ est compacte. Avec les notations précédentes, $\Sigma_X$ constitue toujours une triangulation généralisée de $X$, de squelette $S^\an(X)$. Réciproquement, si $S$ est une triangulation généralisée de $X$ et que $\Gamma$ et $\Gamma^\natural$ désignent les squelette et squelette tronqué, alors $S$ contient tous les noeuds de $\Gamma$ (\cite{Duc}, Théorème $5.1.14$). En revanche il se peut qu'un nœud $s$ de $\Gamma$ ne soit pas un  nœud de $\Gamma^\natural$, cela se produit lorsque $s$ est dans l'intérieur de $X$, de genre $0$, que $s$ est de valence $2$ dans $\Gamma^\natural$ tout en étant un sommet de valence $>2$ dans $\Gamma$, autrement dit si parmi les composantes connexes de $X\setminus S$ issues de $s$ il y en a exactement deux qui sont des couronnes relativement compactes, et au moins une qui est non relativement compacte. 
 \end{rem}

\begin{rem}
Un point de $\Gamma^\natural$ qui est unibranche dans $\Gamma$ est un nœud de $\Gamma$ si et seulement s'il est soit de genre strictement positif, soit appartient au bord analytique $\partial^\an X$ de $X$. Cette définition diffère légèrement de celle de \cite{Duc}, cet écart nous permettant simplement d'alléger par la suite certains énoncés et d'affirmer qu'une triangulation généralisée non réduite à un singleton d'une courbe $k$-analytique est minimale si et seulement si elle est constituée exactement des nœuds de son squelette, comme nous le montre la proposition suivante.
\end{rem}\

 \begin{prop}[Points superflus d'une triangulation généralisée, cf. \cite{Duc}, $5.1.16$ pour le cas compact]\label{points superflus}  Soit $X$ une courbe $k$-analytique quasi-lisse et $S$ une triangulation généralisée de $X$ dont on note $\Gamma$ le squelette. Soit $s\in S$ et $\Sigma= S\setminus\lbrace s\rbrace$, notons $\Gamma_0$ la composante connexe de $\Gamma\setminus \Sigma$ contenant $s$. Le sous-ensemble $S\setminus \lbrace s\rbrace$ est une triangulation généralisée si et seulement si les trois conditions suivantes sont satisfaites : 
\begin{itemize}
\item[i)] $s\not\in \partial^\an X$
\item[ii)]$g(s)=0$
\item[iii)] $\Gamma_0$ est un intervalle ouvert ou bien un intervalle semi-ouvert d'extrémité $s$.
\end{itemize}
\end{prop}\

Si $S$ est une triangulation généralisée de $X$ de squelette $\Gamma$, nous savons d'une part que $S$ contient tous les nœuds de $\Gamma$, d'autre part que $\Gamma$ contient le squelette analytique $S^\an(X)$ de $X$. Par conséquent toute triangulation généralisée de $X$ contient les nœuds de $S^\an(X)$. Voici un résultat de minimalité de triangulation pour les courbes compactes : 
 
 \begin{prop}[\cite{Duc}, $5.4.12$]\label{triangulation minimale dans le cas compact} 
 Soit $X$ une courbe $k$-analytique quasi-lisse, connexe et compacte. Soit $\Sigma_X$ l'ensemble (potentiellement vide) des nœuds de $S^\an (X)$. Trois cas sont possibles : \bigskip
 \begin{enumerate}
 \item \emph{$\Sigma_X$ est non vide} : ce cas se présente si et seulement si $X$ possède une plus petite triangulation, auquel cas cette triangulation coïncide avec $\Sigma_X$. Si $\vert k^\times \vert \neq \lbrace 1\rbrace$ et si $X$ est strictement $k$-analytique, alors $\Sigma_X$ est une triangulation de type $2$.
 \item \emph{$S^\an (X)$ est non vide mais $\Sigma_X$ est vide} : cela se produit si et seulement si $X$ est une courbe de Tate, i.e. isomorphe à un quotient de la forme $\Gm^\an/{q^\Z}$ pour un certain $q\in k^*$ vérifiant $\vert q \vert < 1$ (cela n'est possible que lorsque $k$ est non trivialement valué). Dans ce cas $S^\an (X)$ est un cercle, et les triangulations minimales de $X$ sont les singletons $\{x\}$ avec $x\in S^\an (X)$. $X$ possède donc une infinité de triangulations minimales, dont une infinité sont de type $2$.
 \item \emph{$S^\an (X)$ est vide} : cela se présente si et seulement si $X\simeq \P_k^{1, \an}$. Dans ce cas les triangulations minimales sont les singletons $\{x\}$ avec $x\in X_{[2,3]}$. La courbe $X$ possède donc une infinité de triangulations minimales, dont une infinité sont de type $2$ dès que $\vert k^\times \vert \neq \{1\}$.
 \end{enumerate}

 \end{prop}

\subsubsection{Semi-graphe d'anabélioïdes d'une triangulation généralisée}

Nous allons maintenant associer à une triangulation généralisée un graphe d'anabélioïdes. Fixons pour la suite de cette partie une courbe $k$-analytique $X$ connexe, quasi-lisse et non vide munie d'une triangulation généralisée $S\subset X_{[2, 3]}$. Commençons par lui associer un graphe (resp. semi-graphe) qui n'est autre que l'analogue du squelette tronqué (resp. squelette) de la triangulation  généralisée S.\\

Notons $\mathcal{C}^\natural(X,S)$ (respectivement $\mathcal{C}^\infty(X,S)$) l'ensemble des composantes connexes de $X\setminus S$ qui sont des couronnes relativement compactes (respectivement des composantes non relativement compactes) et $\DX$ l'ensemble des composantes connexes de $X\setminus S$ qui sont des disques. Si $\CC\in \mathcal{C}^\natural(X,S)$, l'inclusion $S^\an(\CC)\hookrightarrow X$ s'étend à la compactification naturelle $\widehat{S}^\an(\CC)$ (cf. \ref{chapeautée d'une couronne}) en une application $\zeta_\CC : \widehat{S}^\an(\CC)\rightarrow X$ (mais celle-ci n'est pas forcément injective, en particulier si la couronne $\CC$  \og boucle sur elle-même \fg, i.e. si $\widehat{S}^\an(\CC)$ est homéomorphe à un cercle).\\

\begin{defi}\

\begin{enumerate}
 \item Le graphe localement fini associé à la triangulation généralisée $S$, noté $\G^\natural(X,S)$, est défini par : 
\begin{itemize}
\item[•] ses sommets correspondent aux éléments de $S$;
\item[•] ses arêtes correspondent aux éléments de $ \mathcal{C}^\natural(X,S)$ : si $e$ est une arête nous désignerons par $\CC_e$ la couronne associée, $e$ étant alors définie comme l'ensemble constitué des deux bouts de $\CC_e$;
\item[•] les branches (au sens des semi-graphes) d'une arête $e$ correspondent aux bouts (au sens $k$-analytique) de $\CC_e$, et une branche $b$ de $e$ aboutit à $\zeta_{\CC_e}(b)$.
\end{itemize}
\item Le semi-graphe localement fini associé à la triangulation généralisée $S$, noté $\G(X,S)$, est défini par : 
\begin{itemize}
\item[•] ses sommets et ses arêtes fermées (c'est-à-dire contenant deux branches) correspondent exactement  à ceux de $\G^\natural(X,S)$;
\item[•] ses arêtes ouvertes (contenant une seule branche) correspondent aux éléments de $\mathcal{C}^\infty(X,S)$, i.e. aux arêtes non relativement compactes : si $e$ est une telle arête nous désignerons par $\CC_e$ la composante associée, $e$ étant alors définie comme le singleton constitué de l'unique bout de $\CC_e$ correspondant à son bord topologique dans $X$;
\item[•] l'unique branche d'une arête ouverte $e$ correspondant à une composante $\CC_e$ aboutit au sommet de $S$ qui correspond au bord topologique de $\CC_e$ dans $X$.
\end{itemize}
 \end{enumerate} 

\end{defi}

Si $s\in S$, définissons $\mathsf{Comp}(X,S,s)$ comme l'image de l'application naturelle : $\br(X,s)\rightarrow \pi_0(X\setminus S)$. Soit $(\widetilde{X},\widetilde{s})$ le revêtement topologique universel de $(X,s)$, et $\widetilde{S}\subset \widetilde{X}$ l'ensemble des préimages des éléments de $S$. Définissons : 
\begin{align*}
&\CC(X,S,s) =\left( \mathcal{C}^\natural(X,S) \cup \mathcal{C}^\infty(X,S) \right) \cap \mathsf{Comp}(X,S,s) \\
&\mathcal{D}(X,S,s) =\DX \cap \mathsf{Comp}(X,S,s) \\
&\St(X,S,s) = \lbrace \widetilde{s}\rbrace\cup \bigsqcup_{T\in \mathsf{Comp}(\widetilde{X},\widetilde{S},\widetilde{s})} T \;\;\;\; \text{: \;«\;l'étoile en $s$ de la triangulation généralisée\;$S$ »}
\end{align*}

\begin{defi}[Graphe d'anabélioïdes associé à une triangulation généralisée]\label{définition du semi-graphe d'anabélioide associé à une triangulation} \

Le graphe $\G^\natural(X,S)$ (resp. le semi-graphe $\G(X,S)$) peut être enrichi en un graphe d'anabélioïdes (resp. semi-graphe d'anabélioïdes) que l'on notera $$\GG^\natural(X,S) \;\; (\text{resp.} \;\;\GG(X,S))$$ et qui sera défini par : 
\begin{itemize}
\item[•]si $s\in S$ est un sommet, $\GG_s = \Covt(\mathrm{St}(X,S,s))$
\item[•]si $e\in \mathcal{C}^\natural(X,S) \cup \mathcal{C}^\infty(X,S)$ est une arête, $\GG_e = \Covt(\CC_e)$
\item[•]si $b$ est une branche de $e$ aboutissant à $s$, alors $$b : \Covt(\CC_e)\rightarrow \Covt(\mathrm{St}(X,S,s))$$ est défini par le foncteur exact $b^* : \Covt(\mathrm{St}(X,S,s))\rightarrow \Covt(\CC_e)$ induit naturellement par l'unique plongement $\CC_e\hookrightarrow \mathrm{St}(X,S,s)$ qui relève à $\widetilde{X}$ le plongement naturel $\CC_e\hookrightarrow X$ et selon lequel le relevé de la branche $b$ aboutit à $\widetilde{s}$.
\end{itemize}

\end{defi}

\begin{rem}
Le graphe $\G^\natural(X, S)$ prend en compte uniquement les couronnes \emph{relativement compactes} (i.e. les éléments de $\mathcal{C}^\natural(X,S)$) parmi les arêtes qui apparaissent comme composantes connexes de $X\setminus S$. Néanmoins la situation est différente pour le graphe d'anabélioïdes $\GG^\natural(X,S)$, dont $\G^\natural(X,S)$ est le graphe sous-jacent : la définition même de l'anabélioïde connexe associé à un sommet permet de garder une trace des composantes non relativement compactes (i.e. des éléments de $\mathcal{C}^\infty(X,S)$) issues de ce sommet. 
\end{rem}

\begin{rem}
Le passage par le revêtement topologique universel dans la définition de l'étoile $\St(X,S,s)$ implique que si une couronne $\CC_e\in \CC(X,S,s)$ \og boucle sur $e$ \fg{} (ses deux bouts sont constitués du sommet $s$), elle apparaît deux fois dans $\St(X,S,s)$ via ses deux branches. Sans un tel recours au revêtement universel le graphe d'anabélioïdes $\GG^\natural(X,S)$ ne saurait être totalement détaché. La compatibilité de $\GG^\natural(X,S)$ avec les hypothèses mochizukiennes sera investi dans la section \ref{Compatibilité}.
\end{rem}

\subsubsection{Revêtements de Kummer, $\mu_\ell$-torseurs et cochaînes harmoniques}

Si $X$ est un espace $k$-analytique et si $\ell\in\N^*$ est un entier inversible dans $k$, la \emph{suite exacte de Kummer} sur $X_{\et}$ 

$$ 1\longrightarrow \mu_l \longrightarrow \G_m\overset{z\mapsto z^\ell}{\longrightarrow }\G_m\longrightarrow1 $$

induit une injection $$\mathscr{O}_X(X)^\times / (\mathscr{O}_X(X)^\times)^\ell\overset{\iota}{\hookrightarrow} \HH^1(X_\et, \mu_\ell)$$ dont on notera on notera $\Kum_\ell(X)$ l'image. Nous savons d'après \cite{Ber2} que tout faisceau étale localement constant sur $X_\et$ est représentable, ce qui est en particulier le cas de tout $\mu_\ell$-torseur. Par conséquent $\HH^1(X_\et, \mu_\ell)$ classe tous les $\mu_\ell$-torseurs étales \emph{analytiques} sur $X$ à isomorphisme près. Si $f \in \mathscr{O}_X(X)^\times$ son image $(f)$ dans $\HH^1(X_\et, \mu_\ell)$ par $\iota$ correspond à $ \mathscr{M}(\mathscr{O}_X[T]/(T^\ell-f)) $. Les éléments de $\Kum_\ell(X)$ vus comme $\mu_\ell$-torseurs étales analytiques seront appelés les $\mu_\ell$-torseurs \emph{de Kummer}.\

\begin{ex}
Si $I$ est un intervalle non vide de $\R_+^*$, la fonction $T^\ell$ induit un revêtement fini et plat (en fait étale modéré) $\mathcal{C}(I)\rightarrow \mathcal{C}(I^\ell)$ qui permet d'identifier $\mathcal{C}(I)$ au $\mu_\ell$-torseur de Kummer $\mathscr{M}(\mathscr{O}_{\mathcal{C}(I^\ell)}[T]/(T^\ell-S))$, où $S$ est la coordonnée standard de $\mathcal{C}(I^\ell)$.
\end{ex}

\begin{prop}[\cite{Duc}, $3.6.30$ et $3.6.31$]\label{propriétés kummériennes des revêtements de couronnes} Soit $X$ une $k$-couronne :
\begin{enumerate}
\item Si $Y$ est une sous-couronne de $X$ (voir définition \ref{définition d'une sous-couronne}), l'injection de $Y$ dans $X$ induit un isomorphisme de groupes :  $\Kum_\ell(X)\rightarrow \Kum_\ell(Y)$
\item Comme $k^\times$ est $\ell$-divisible (puisque $k$ est supposé algébriquement clos), $\Kum_\ell(X)$ est isomorphe à $\Z/\ell\Z$. Cet isomorphisme est non canonique mais le devient dès que l'on fixe une orientation de $X$.
\item Toute composante connexe d'un $\mu_\ell$-torseur de Kummer de $X$ est une couronne. 
\end{enumerate}
\end{prop}\

\begin{defi}[Cochaînes harmoniques sur un semi-graphe localement fini]\label{cochaînes harmoniques} 
Soit $\Gamma$ un semi-graphe localement fini et $A$ un groupe abélien. Une \emph{$A$-cochaîne harmonique sur $\Gamma$} est une application : $c : \lbrace \text{arêtes orientées de } \Gamma\rbrace \rightarrow A$ satisfaisant les deux propriétés suivantes : 
\begin{enumerate}
\item si $e$ et $e'$ correspondent à la même arête munie de ses deux orientations distinctes : $c(e')=-c(e)$. 
\item si $x$ est un sommet de $\Gamma$ : $$\sum_{\text{arêtes orientées vers }e}c(e)=0_A. $$
\end{enumerate}
L'ensemble des $A$-cochaînes harmoniques de $\Gamma$ forme un groupe abélien noté $\mathrm{Harm}(\Gamma, A)$.

\end{defi}\

\begin{lem}\label{lemme sur les cochaines harmoniques} 
Soit $\Gamma$ un semi-graphe connexe, et $e$, $e'$ et $e''$ trois arêtes ouvertes distinctes de $\Gamma$, orientées vers le sommet auquel elles aboutissent. Soit $A$ un groupe commutatif et $(a, a')\in A^2$. Alors il existe au moins un élément $c\in \mathrm{Harm}(\Gamma, A)$ tel que $c(e)=a, c(e')=a'$ et $c(e'')=-(a+a')$. On peut par ailleurs imposer $c(f)=0$ pour toute autre arête ouverte $f$ de $\Gamma$. 
\end{lem}

\begin{proof}
Soient $x, x'$ et $y$ les sommets auxquels aboutissent $e, e'$ et $e''$ ($x, x'$ et $y$ ne sont pas forcément distincts). Il existe des suites finies $(x=x_1, \ldots, x_r=y)$ et $(x'=x'_1, \ldots, x'_s=y)$ de sommets de $\Gamma$ distincts et consécutifs, de telle sorte qu'il existe des arêtes orientées $e_i$ (resp. $e'_j$) de $\Gamma$ reliant $x_i$ à $x_{i+1}$ (resp. $x'_j$ à $x'_{j+1}$) et orientées vers $x_{i+1}$ (resp. $x'_{j+1}$). Soit $c$ (resp. $c'$) $\in \mathrm{Harm}(\Gamma, A)$ défini par $c(e)=c(e_i)=-c(e'')=a$ pour tout $i=1\ldots r-1$ et nul sur les autres arêtes (resp. $c'(e')=c'(e'_j)=-c'(e'')=a'$ pour tout $j=1\ldots s-1$ et nul sur les autres arêtes). Alors la cochaîne harmonique $c+c'\in \mathrm{Harm}(\Gamma, A)$ vérifie les propriétés requises. 
\end{proof}\

Soit $X$ une courbe $k$-analytique, $\pi$ le morphisme du site étale $X_{\et}$ vers le site topologique $X_{\mathrm{top}}$ et $\ell$ un entier premier à $p$. La suite spectrale $$\HH^p(X_{\mathrm{top}}, \mathrm{R}^q\pi_*\mu_\ell)\Rightarrow \HH^{p+q}(X_{\et}, \mu_\ell) $$ induit la suite exacte : \begin{equation}\label{suite exacte 1} 
0\rightarrow \HH^1(X_{\mathrm{top}}, \mu_\ell )\rightarrow   \HH^1(X_{\et}, \mu_\ell )\rightarrow \HH^0(X_{\mathrm{top}}, \mathrm{R}^1\pi_* \mu_\ell)\rightarrow 0.
\end{equation}

Par ailleurs, si la courbe $X$ est supposée \emph{lisse} (i.e. quasi-lisse et sans bord), si $S$ est une triangulation généralisée de $X$ de squelette $\Gamma$, comme $k$ est algébriquement clos on a la suite exacte suivante (cf. \cite{Duc}, $5.2.4$) : 
\begin{equation}\label{suite exacte 2} 
0\rightarrow \prod_{x\in S}\mathbb{J}_x\rightarrow  \HH^0(X_{\mathrm{top}}, \mathrm{R}^1\pi_* \mu_\ell)\rightarrow \mathrm{Harm}(\Gamma, \Z/\ell \Z)\rightarrow 0,
\end{equation}
où le groupe $\mathbb{J}_x$ est défini comme suit : 
\begin{itemize}
\item[•]si $x$ est de type $3$ alors $\mathbb{J}_x$ est trivial;
\item[•]si $x$ est de type $2$ alors $\mathbb{J}_x=_\ell\mathscr{J}(\widetilde{k})$ est le groupe de $\ell$-torsion des $\widetilde{k}$-points de la jacobienne de la courbe résiduelle $\mathscr{C}_x$ (cf. \ref{courbe résiduelle}).
\end{itemize}\

\begin{rem}
Ducros obtient la suite exacte \ref{suite exacte 2} dans \cite{Duc}, $5.2.4$ en se plaçant dans le cadre d'une vraie triangulation. Néanmoins, son raisonnement reste valable dans notre cadre des triangulations généralisées. En effet, on remarque qu'une telle suite exacte appliquée à une rabotée $X^\circ$ de $X$ (voir la définition \ref{définition d'une rabotée}, $S$ est alors une vraie triangulation de $X^\circ$) ne dépend pas du choix de $X^\circ$ : deux rabotées différentes donneront deux suites exactes isomorphes, ce qui permet d'étendre de telles suites exactes à $X$ tout entier.
\end{rem}

Quelque soit le type du point $x$, le groupe $\mathbb{J}_x$ est isomorphe (non canoniquement en général) à $(\Z/\ell\Z)^{2g(x)}$. 
La combinaison des deux suites exactes \ref{suite exacte 1} et  \ref{suite exacte 2} et le fait que $\Gamma$ soit homotope à $X$ mène au résultat suivant : 

\begin{prop}\label{diagramme de suites exactes} 
On a le diagramme suivant, dont la ligne et la colonne sont exactes :\

 \xymatrix{
    & & 0 \ar[d] & & \\
    & & \HH^1(\Gamma_{\mathrm{top}}, \mu_\ell)\ar[d] & & \\
    & & \HH^1(X_{\et}, \mu_\ell)\ar[d] & & \\
    0 \ar[r] & \prod_{x\in S}\mathbb{J}_x \ar[r] & \HH^0(X_{\mathrm{top}}, \mathrm{R}^1\pi_* \mu_\ell) \ar[r] \ar[d] & \mathrm{Harm}(\Gamma, \Z/\ell \Z) \ar[r] & 0 \\
    & & 0 & &
  }

\end{prop}

\subsubsection{Revêtements sauvages et condition de déploiement}\label{sous-section p-revêtements} \

Supposons dans cette sous-section que $k$ est de caractéristique nulle et $p>0$, où $p$ désigne la caractéristique résiduelle de $k$. Supposons la norme $\vert \cdot \vert$ sur $k$ normalisée de telle sorte que  $\vert p \vert=p^{-1}$.

\begin{lem}\label{lemme sur la distance entre deux racines de l'unité} 
Soient $\xi$ et $\xi'$ deux racines $p$-ièmes de l'unité distinctes dans $k$ (supposé algébriquement clos). Alors $\vert \xi - \xi'\vert= p^{-\frac{1}{p-1}}$.
\end{lem}

\begin{proof}
Soit $\displaystyle{\Phi_p=\frac{X^p-1}{X-1}=\sum_{i=0}^{p-1}X^i=\prod_{\xi\in \mu'_p}X-\xi\in \Q[X]}$ le $p$-ième polynôme cyclotomique, où $\mu_p'$ désigne l'ensemble des $p-1$ racines primitives $p$-ièmes de l'unité dans $k$. L'évaluation en $1$ donne : $p=\prod_{\xi\in \mu'_p}1-\xi.$ Pour $\xi$ décrivant $\mu'_p$ tous les $1-\xi$ ont la même norme en tant que conjuguées par $\Gal(k/\Q)$. De là on déduit que $\vert 1-\xi \vert=p^{-\frac{1}{p-1}}$, puis on obtient le résultat en remarquant que la multiplication par n'importe quelle racine $p$-ième de l'unité est une isométrie de $k$. 
\end{proof}\

La proposition suivante, donnant les ensembles de déploiement du $\mu_{p^h}$-torseur donné par la fonction $\sqrt[p^h]{1+T}$, sera importante dans la suite du texte.

\begin{prop}\label{décomposition du torseur sauvage} 
Si $h\in \N^*$, le revêtement étale $\G_m^\an \xrightarrow{z\mapsto z^{p^h}} \G_m^\an$ est totalement déployé au-dessus du point $\eta_{z_0, r}$ vérifiant $r<\vert z_0\vert=:\alpha$ si et seulement si : $r<\alpha p^{-h-\frac{p}{p-1}}$. Plus précisément, l'image réciproque de $\eta_{z_0, r}$ contient :
\begin{itemize}
\item[•]un seul élément quand $r\in [\alpha p^{-\frac{p}{p-1}}, \alpha];$
\item[•]$p^i$ éléments quand $r\in[\alpha p^{-i-\frac{p}{p-1}}, \alpha p^{-i-\frac{1}{p-1}}[,$ avec $1\leqslant i\leqslant h-1$;
\item[•]$p^h$ éléments quand $r\in[0, \alpha p^{-h-\frac{1}{p-1}}[.$
\end{itemize}
\end{prop}

\begin{proof}
Soit $f: \G_m^\an\rightarrow \G_m^\an$ le revêtement donné par $f(z)=z^{p}$. Soit $z_1\in k^*$ et $\rho\in \R_+$ vérifiant $\rho<\vert z_1 \vert$ (de telle sorte que $\eta_{z_1, \rho} \notin ]0, \infty[$).

\item Afin de calculer $f(\eta_{z_1,\rho})$, remarquons que pour tout polynôme $P\in k[T]$ :

\begin{align*}
\vert P\left(f(\eta_{z_1, \rho})\right) \vert &= \vert (P\circ f)(\eta_{z_1, \rho})\vert= \vert P(T^p)(\eta_{z_1, \rho})\vert\\
&=\sup_{x\in \mathrm{B}(z_1, \rho)}\vert P\circ f(x)\vert\\
&=\sup_{y\in f(\mathrm{B}(z_1, \rho))} \vert P(y)\vert
\end{align*}

Comme $k$ est algébriquement clos, il existe $\widehat{\rho}>0$ tel que $f(\mathrm{B}(z_1, \rho))=\mathrm{B(z_1^p,\widehat{\rho}\,)}$, d'où l'on obtient $f(\eta_{z_1, \rho})=\eta_{z_1^p,\widehat{\rho}}$.\\

Afin de calculer $\widehat{\rho}$, remarquons que $\widehat{\rho}=\vert (T-z_1^p)(f(\eta_{z_1, \rho}))\vert=\vert (T^p-z_1^p)(\eta_{z_1, \rho})\vert$, et : 
$$T^p-z_1^p=\sum_{i=1}^p \binom{p}{i}z_1^{p-i}(T-z_1)^i=\sum_{i=1}^p \gamma_i (T-z_1)^i ,$$
où l'on pose $\gamma_i= \binom{p}{i}z_1^{p-i}$, avec : 
$$
\vert \gamma_i \vert = \left\{
    \begin{array}{ll}
        1& \hbox{si } i=p \\
        p^{-1}\vert z_1 \vert ^{p-i} & \hbox{si } 1\leqslant i\leqslant p-1
    \end{array}
\right.
$$
Ainsi $\widehat{\rho}=\vert (T-z_1^p)(f(\eta_{z_1, \rho}))\vert=\max_{1\leqslant i \leqslant p} \lbrace\vert \gamma_i \vert \rho^i \rbrace=\max \lbrace \rho^p, \left( p^{-1}\rho^i \vert z_1 \vert^{p-i}\right)_{1 \leqslant i \leqslant p-1}  \rbrace $. Puisque l'on a supposé $\rho< \vert z_1 \vert$, alors on obtient $\widehat{\rho}= \max \lbrace \rho^p, p^{-1}\rho\vert z_1 \vert^{p-1}\rbrace $, soit : 
$$
\widehat{\rho} = \left\{
    \begin{array}{ll}
         p^{-1}\rho \vert z_1 \vert^{p-1}& \hbox{si } \rho \leqslant \vert z_1 \vert p^{-\frac{1}{p-1}} \\
        \rho^p & \hbox{si } \rho \geqslant \vert z_1 \vert p^{-\frac{1}{p-1}}
    \end{array}
\right.
$$

De là : 

$$
f(\eta_{z_1, \rho})= \left\{
    \begin{array}{ll}
         \eta_{z_1^p,p^{-1}\rho \vert z_1 \vert^{p-1}}& \hbox{si } \rho \leqslant \vert z_1 \vert p^{-\frac{1}{p-1}} \\
        \eta_{z_1^p, \rho^p} & \hbox{si } \rho \geqslant \vert z_1 \vert p^{-\frac{1}{p-1}}
    \end{array}
\right.
$$

\bigskip
\item Cherchons maintenant les préimages par $f$ de $\eta_{z_0, r}$, où $0\leqslant r<\alpha:=\vert z_0 \vert$. Posons : 

$$
\widetilde{r} = \left\{
    \begin{array}{ll}
         rp \alpha^{-\frac{p-1}{p}}& \hbox{si } r \leqslant \alpha p^{-\frac{p}{p-1}} \\
        r^{\frac{1}{p}} & \hbox{si } r \geqslant \alpha p^{-\frac{p}{p-1}} 
    \end{array}
\right.
$$
D'après ce qui précède, si $\widetilde{z_0}$ désigne une racine $p$-ième de $z_0$, alors : $$\eta_{\widetilde{z_0}, \widetilde{r}}\in f^{-1}\left(\lbrace \eta_{z_0,r}\rbrace \right),$$ et $f^{-1}\left(\lbrace \eta_{z_0,r}\rbrace \right)$ est constituée de tous les conjugués $\eta_{\xi\widetilde{z_0}, \widetilde{r}}$ de $\eta_{\widetilde{z_0}, \widetilde{r}}$ pour $\xi\in \mu_p$. De là : 
$$
f^{-1}\left(\lbrace \eta_{z_0,r}\rbrace \right) = \left\{
    \begin{array}{ll}
         \lbrace \eta_{\xi \widetilde{z_0},rp\alpha^{-\frac{p-1}{p}}}\rbrace_{\xi\in \mu_p} & \hbox{si } r \leqslant \alpha p^{-\frac{p}{p-1}} \\
         \lbrace \eta_{\xi \widetilde{z_0},r^{\frac{1}{p}}}\rbrace_{\xi\in \mu_p} & \hbox{si } r \geqslant \alpha p^{-\frac{p}{p-1}} 
    \end{array}
\right.
$$
Puisque $\vert \widetilde{z_0} \vert= \alpha^{\frac{1}{p}}$, on a $\vert \xi \widetilde{z_0}- \xi' \widetilde{z_0}\vert=\alpha^{\frac{1}{p}} p^{-\frac{1}{p-1}}$ dès que $\xi\neq \xi' \in \mu_p$, d'après le lemme \ref{lemme sur la distance entre deux racines de l'unité}. On en déduit que $f^{-1}\left(\lbrace \eta_{z_0,r}\rbrace \right)$ a un unique élément si $r \geqslant \alpha p^{-\frac{p}{p-1}}$, $p$ sinon. 
\item Dans le cas général, avec $h\geqslant 1$, un raisonnement par récurrence sur $h$ permet de conclure.

\end{proof}

\begin{figure}
\includegraphics[width=1\textwidth]{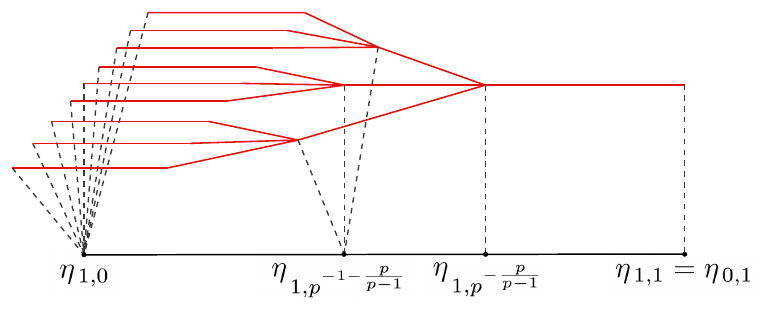}
\caption{Revêtement $\G_m^\an \xrightarrow{z\mapsto z^{p^h}} \G_m^\an$ avec $p=3$, $h=2$ et $z_0=1$.}
\end{figure}

\begin{coro}\label{corollaire permettant de trouver un revêtement d'une couronne totalement décomposé en un de ses bouts} 
Soit $\CC$ une couronne $k$-analytique ouverte de module fini, que l'on peut supposer être un domaine analytique de $\P_k^{1,\an}$. Soit $S^\an(\CC)=]x,y[$ son squelette analytique, avec $x, y \in \P_k^{1,\an}$, et $h\in \N^*$. Si $\ell(\CC)>h-1$ et $y_0\in S^\an(\CC)_{[2]}$ est tel que $$\varepsilon:= d(y_0,y)<\min\left(\frac{p}{p-1}, \ell(\CC)-(h-1)\right),$$ il existe un $\mu_{p^h}$-torseur analytique $f\in H^1(\CC,\mu_{p^h})$ tel que : 
\begin{itemize}
\item[•] l'ensemble des points de $S^\an(X)$ en lesquels $f$ n'a qu'un seul antécédent est $[y_0,y[$, en particulier $f$ est non déployé au dessus de $[y_0,y[$,
\item[•]$f$ est totalement déployé sur un voisinage ouvert de l'extrémité donnée par $x$; plus précisément $f$ est totalement déployée au-dessus de $\ds{r^{-1}\left(\,]x, y_1[   \,\right)}$ où $y_1\in ]x,y_0]_{[2]}$ est tel que $d(y_1,y_0)=h-1$.
\end{itemize}

\end{coro}

\begin{rem}
Avec les notations ci-dessus, si $f$ correspond au torseur $Y\to X$, l'unique antécédent de $y_0$ est un \og point de branchement \fg{} de $S^\an(Y)$, c'est-à-dire que c'est un nœud de $S^\an(Y)$. La condition $d(y_0,y)<\frac{p}{p-1}$ implique que l'on peut prendre $y_0$ arbitrairement proche de $y$. 
\end{rem}

\begin{proof}
Si $I$ désigne un intervalle de $\R_+^*$, notons $\mathcal{C}_{(I)}$ le domaine analytique de $\A_k^{1, \an}$ défini par la condition $\mid T-1\mid \in I$, c'est une couronne. Par ailleurs si $J$ est un autre intervalle de $\R_+^*$, $\mathcal{C}_{(I)}$ et $\mathcal{C}_{(J)}$ sont isomorphes si et seulement si $J\in \vert k^\times\vert^{\pm 1}$. Par conséquent la couronne $\mathcal{C}$ est isomorphe à $\mathcal{C}_{(I)}$ pour un intervalle ouvert $I\varsubsetneq ]0,1[$ vérifiant $\sup I\neq 1$, $\inf I\neq 0$ (car $\mathcal{C}$ est supposée de longueur finie), et tel que via cet isomorphisme le bout $x$ corresponde à $\inf I$. 

Posons :$$\left\{
\begin{array}{l}
  r_1=p^{\;\varepsilon-\frac{p}{p-1}} \\
  r_2=\frac{r_1}{\mathrm{Mod}(\CC)}
\end{array}
\right.$$  

On a alors
$$r_1\in \left] p^{-\frac{p}{p-1}}, \min\lbrace 1,\Mod(\mathcal{C})p^{-(h-1)-\frac{p}{p-1}}\rbrace \right[\cap (\sup I)\cdot \vert k^\times\vert,$$ 
\item et $$0<r_2<p^{-(h-1)-\frac{p}{p-1}}<p^{-\frac{p}{p-1}}<r_1<1.$$ Ainsi, par construction, $\mathcal{C}$ est isomorphe à la couronne $\mathcal{C}_{(\, ]r_2,r_1[\,)}$, d'où un isomorphisme $$\psi : \mathcal{C}\iso \mathcal{C}_{ (\,]r_2,r_1[\,)}$$ qui identifie $y_0, x$ et $y$  respectivement à $\eta_{1,p^{-\frac{p}{p-1}}}, \eta_{1,r_2}$ et $\eta_{1,r_1}$. Notons $r$ la rétraction de $\mathcal{C}_{(\, ]r_2,r_1[\,)}$ sur son squelette analytique $]\eta_{1,r_2}, \eta_{1,r_1}[$. 
\item La proposition \ref{décomposition du torseur sauvage} nous indique que le $\mu_{p^h}$-torseur $f: \G_m^\an \xrightarrow{z\mapsto z^{p^h}} \G_m^\an$ vu comme revêtement de $\mathcal{C}_{ ]r_2,r_1[}$ (on le restreint à $f^{-1}(\mathcal{C}_{ ]r_2,r_1[})$) est :
\begin{itemize}
\item[•]totalement déployé au-dessus de $r^{-1}\left(]\eta_{1,r_2}, \eta_{1,p^{-(h-1)-\frac{p}{p-1}}} [\right)$
\item[•]non déployé avec un seul antécédent au-dessus de $[\eta_{1,p^{-\frac{p}{p-1}}}, \eta_{1,r_1}[ ,$ 

\end{itemize}\
ce qui permet de conclure. 
\end{proof}

\subsection{Groupe fondamental tempéré d'une courbe analytique}\label{groupe fondamental tempéré} \

Soit $K$ un corps valué non archimédien complet et $p$ l'exposant caractéristique de son corps résiduel. Soit $X$ un espace strictement $K$-analytique quasi-lisse et connexe muni d'un point géométrique $x\in X$ (i.e. un point défini sur une extension complète et algébriquement close de $K$). Nous allons définir les \emph{revêtements tempérés} de $X$ ainsi que le \emph{groupe fondamental tempéré} $\pi_1^\tp(X,x)$ de $X$ en $x$, notions introduites par Yves André dans \cite{And} . Ce dernier se place dans le cadre d'un sous-corps complet $K$ de $\C_p$ (en désignant par $\C_p$ la complétion d'une clôture algébrique fixée $\overline{\Q}_p$ de $\Q_p$). Or ces définitions s'étendent \emph{verbatim} à tout corps non archimédien, ce que nous nous permettrons de faire librement dans la suite.\\

Rappelons que De Jong définit dans \cite{DJg} un \emph{revêtement étale d'espaces $K$-analytiques} $Y\rightarrow X$ comme étant localement sur la base $X$ une union disjointe de revêtements étales finis. Les \emph{revêtements topologiques} (lorsque ces revêtements étales finis sont des isomorphismes) sont en particulier des revêtements étales, mais de nombreux revêtements étales ne donnent lieu à aucun revêtement topologique.

Les revêtements topologiques de l'espace $K$-analytique pointé $(X,x)$ sont classés par le \emph{groupe fondamental topologique} classique $\pi_1^\mathrm{top}(X,x)$, tandis que De Jong définit dans \cite{DJg} le \emph{groupe fondamental de $X$ en $x$}, $\pi_1^\et(X,x)$, qui classe les revêtements étales de $X$, et dont  $\pi_1^\mathrm{top}(X,x)$ est un quotient discret. Son complété profini, le \emph{groupe fondamental algébrique} $\pi_1^\mathrm{alg}(X,x)$, classe les revêtements finis étales. Ce dernier est l'analogue du groupe fondamental classique de la géométrie algébrique, coïncidant avec celui-ci pour l'analytifiée d'une courbe algébrique. En revanche, sur les corps algébriquement clos de caractéristique nulle, ce groupe ne porte aucun phénomène anabélien pour les courbes analytiques de nature algébrique.\\

Les autres groupes ont en pratique un comportement très éloigné de ce qui advient en géométrie complexe :  

\begin{itemize}
\item[•]$\pi_1^\mathrm{top}(X,x)$ est \og trop petit \fg{} , par exemple $\pi_1^\mathrm{top}(\P_{\C_p}^{1,\an}\setminus \lbrace 0,1,\infty\rbrace), x)=\lbrace 1\rbrace$;
\item[•]$\pi_1^\et(X,x)$ est \og trop grand \fg{}, par exemple $\pi_1^\et(\P_{\C_p}^{1,\an},x)$ admet des groupes linéaires de la forme $\mathrm{SL}_n(\Q_p)$ comme quotients (\cite{DJg}, Proposition $7.4$).
\end{itemize}\

C'est pour remédier à ce \og défaut \fg \,qu'Yves André introduit dans \cite{And} la notion de \emph{revêtement tempéré}, qui capture à la fois la théorie des revêtements topologiques et finis étales tout en offrant des comportements raisonnables qui s'apparentent à ce que l'on connaît en géométrie complexe : 

\begin{defi}\

\begin{enumerate}
\item Un \emph{revêtement tempéré} d'un espace strictement $K$-analytique $X$ connexe non vide sera défini comme un revêtement étale \emph{(au sens de \cite{DJg})}\, $\varphi: Y\rightarrow X$ qui devient topologique après un changement de base selon un revêtement étale fini, autrement dit qui est facteur du revêtement topologique d'un revêtement étale fini de $X$, i.e. s'il existe un diagramme commutatif suivant 
$$ \xymatrix{ & Z \ar[ld]_\psi \ar[rd]&\\ W \ar[rd]_\chi & & Y \ar[ld]^\varphi\\ & X &}$$ où $\chi$ est un revêtement étale fini et $\psi$ un revêtement topologique.
\item Soit $\Cov^\tp(X)$ la catégorie des revêtements tempérés de $X$. Si $x\in X$ est un point géométrique, soit $$F_x : \Cov^\tp(X)\to \mathrm{Set}$$ le foncteur fibre en $x$ qui associe à un revêtement $Y\to X$ la fibre $Y_x$ en $x$. Le \emph{groupe fondamental tempéré pointé en $x$} est défini comme le groupe d'automorphismes du fonceur fibre en $x$ : 
$$\pi_1^\tp(X,x):=\Aut(F_x).$$ 
Le groupe $\pi_1^\tp(X,x)$ devient un groupe topologique en considérant l'ensemble des sous-groupes stabilisateurs $(\mathrm{Stab}_{F_x(Y)}(y))_{Y\in \Cov^\tp(X),\; y\in F_x(Y)}$ comme base d'ouverts. C'est alors un groupe topologique prodiscret. 

\bigskip
Si $x$ et $x'$ sont deux points géométriques distincts, les foncteurs fibre $F_x$ et $F_{x'}$ sont (non canoniquement) isomorphes, et tout automorphisme de $F_x$ induit un automorphisme intérieur de $\pi_1^\tp(X,x)$. Ainsi, on peut considérer le \emph{groupe fondamental tempéré $\pi_1^\tp(X)$ de $X$}, défini à unique isomorphisme extérieur près.
\end{enumerate}
\end{defi}

\begin{rem}Il y a équivalence de catégories entre la catégorie des sommes directes de revêtements tempérés de $X$ et la catégorie des $\pi_1^\tp(X,x)$-ensembles discrets sur lesquels $\pi_1^\tp(X,x)$ agit à gauche continûment. 
\end{rem}\

La classe d'isomorphisme de $\pi_1^\tp(X,x)$ ne dépend pas du point géométrique $x$ considéré, de sorte que l'on pourra parler du \emph{groupe fondamental tempéré $\pi_1^\tp(X)$ de $X$}, défini à unique isomorphisme extérieur près.\

\begin{rem}
Une somme directe de revêtements tempérés n'est pas toujours un revêtement tempéré. Supposons que $(H_i)_{i\in\N}$ soit une suite décroissante de sous-groupes ouverts de $\pi_1^\tp(X,x)$ telle que l'intersection $\bigcap_{i\in\N} H_i$ ne soit pas ouverte dans $\pi_1^\tp(X,x)$. Soit $X_i\to X$ le revêtement tempéré de $X$ correspondant au $\pi_1^\tp(X,x)$-ensemble $\pi_1^\tp(X,x)/H_i$. Le revêtement étale $\varphi : \bigsqcup_{i\in \N} X_i\to X$, correspondant au $\pi_1^\tp(X,x)$-ensemble discret $\bigsqcup_{i\in \N} \pi_1^\tp(X,x)/H_i$, n'est en revanche pas tempéré puisque les degrés $[\mathscr{H}(y):\mathscr{H}(\varphi(y))]$ ne sont pas bornés lorsque $y$ parcourt $\bigsqcup_{i\in \N} X_i$.
\end{rem}

\begin{rem}
Comme toute composante connexe d'un bon espace $K$-analytique est dénombrable à l'infini, un revêtement étale $\varphi : Y\to X$ avec $Y$ connexe est tempéré si et seulement s'il s'insère dans un diagramme commutatif similaire à celui ci-dessus avec $\psi$ un revêtement topologique \emph{à fibres dénombrables}. C'est le cas de tout revêtement tempéré galoisien. 
\end{rem}\

Si $(Y,y)$ est un revêtement étale fini galoisien géométriquement pointé de $(X,x)$, alors le revêtement universel pointé $\psi : (\widetilde{Y}, \widetilde{y}) \rightarrow(Y, y)$ reste un revêtement étale galoisien de $X$ de groupe de Galois $\Delta_Y:= \Gal(\widetilde{Y}/X)$ s'insérant dans une suite exacte : 
$$\{1\} \longrightarrow \pi_1^\mathrm{top}(Y,{y})\longrightarrow \Delta_Y\longrightarrow \Gal(Y/X)\longrightarrow \{1\} $$

Tout morphisme $(Y_2,y_2)\rightarrow (Y_1,y_1)$ de revêtements finis étales galoisiens géométriquement pointés de $(X,x)$ induit un morphismes de groupes : $\Delta_{Y_2}\rightarrow \Delta_{Y_1}$.\\

\begin{prop} [\cite{And}]\label{description du pi1 tp de l'espace analytique en terme de limite} 
Soit $X$ un espace strictement $K$-analytique connexe et $x\in X$ un point géométrique. Alors on a l'égalité : 
$$\pi_1^\tp(X,x)= \varprojlim_{i\in I} \Delta_{Y_i} , $$
dès que $(Y_i,y_i)_{i\in I}$ est un système projectif cofinal de revêtements finis étales galoisiens géométriquement pointés au-dessus de $(X,x)$. L'ensemble d'indexation $I$ n'est \textit{a priori} pas dénombrable en général (mais il l'est par exemple dès que $X$ est de nature algébrique et $K$ de caractéristique nulle).
\end{prop}

\begin{rem}
Lorsque que les espaces $K$-analytiques considérés ci-dessus sont des courbes, le groupe  $\pi_1^\mathrm{top}(Y,{y})$ est un groupe libre discret. Cela découle du fait qu'une courbe $K$-analytique a le type d'homotopie d'un graphe. 
\end{rem}

Le groupe fondamental tempéré $\pi_1^\tp(X)$ est par ailleurs un groupe \emph{tempéré} (au sens de la section \ref{section le monde des tempéroïdes}). La catégorie $\Cov^\tp(X)$ des revêtements tempérés de $X$ munie de la classe d'isomorphisme $\Upsilon_X$ des foncteurs $\Cov^\tp(X)\to \mathfrak{Ens}^\mathrm{dd}$ contenant les foncteurs fibres $F_x$ au-dessus des points géométriques est un \emph{tempéroïde connexe} isomorphe à $(\BB^\tp(\pi_1^\tp(X)),\Upsilon_{\pi_1^\tp(X)}^{\mathrm{oubli}})$.\\

Décrivons quelques propriétés du groupe fondamental tempéré : 

\begin{prop}[\cite{And}, $2.1.6$]\label{proprétés topologiques du groupe tempéré}Soit $X$ un bon espace $K$-analytique connexe quasi-lisse, et $x$ un point géométrique : 
\begin{enumerate}
\item Le morphisme naturel $\pi_1^\tp(X,x)\to \pi_1^\mathrm{top}(X,x)$ est surjectif (le groupe tempéré peut donc avoir une infinité de quotients discrets). 
\item Le morphisme naturel $\pi_1^\tp(X,x)\to \pi_1^\mathrm{alg}(X,x)$ est d'image dense, $\pi_1^\mathrm{alg}(X,x)$ peut être identifié à la complétion profinie de $\pi_1^\tp(X,x)$.
\item Si $X$ est une courbe $k$-analytique, le morphisme $\pi_1^\tp(X,x)\to \pi_1^\mathrm{alg}(X,x)$ est injectif.
\end{enumerate}

\end{prop}

Le deuxième point montre que le groupe tempéré n'est pas \og trop petit \fg{}, tandis que le troisième point indique qu'il n'est pas non plus \og trop grand \fg{} lorsque $X$ est une courbe.

\begin{rem}\label{remarque}
Nous savons que la droite affine $\C_p$-analytique est algébriquement simplement connexe : $\pi_1^\mathrm{alg}(\A_{\C_p}^{1,\an})=\lbrace 0\rbrace$ (cela découle du caractère simplement connexe de droite affine algébrique $\A_{\C_p}^1$ et de l'équivalence entre les revêtements étales finis de $\A_{\C_p}^1$ et ceux de $\A_{\C_p}^{1,\an}$). On déduit du troisième point de la proposition \ref{proprétés topologiques du groupe tempéré} que le groupe tempéré de la droite affine est trivial : $\pi_1^\mathrm{temp}(\A_{\C_p}^{1,\an})=\lbrace 0\rbrace$. Cela prouve le caractère non tempéré du revêtement étale donné par le logarithme $\mathrm{log : }\; D_{\C_p}(0, 1)^\circ \to \A_{\C_p}^{1,\an}$. 
\end{rem}\

Dans le cas d'un revêtement étale galoisien, la propriété d'être de nature topologique ou tempérée peut se lire directement sur son groupe de Galois, comme le montre la proposition suivante : 

\begin{prop}[\cite{And}, $2.1.9$]
Soit $X$ un espace strictement $K$-analytique connexe et $Y$ un revêtement étale galoisien de $X$ de groupe de Galois $G$.
\begin{enumerate}
\item Si $G$ est \emph{sans torsion}, $Y$ est un revêtement topologique de  $X$;
\item si $G$ est \emph{virtuellement sans torsion} (i.e. s'il contient au moins un sous-groupe d'indice fini qui est sans torsion), alors $Y$ est un revêtement tempéré de $X$, et plus précisément un revêtement topologique d'un revêtement fini étale galoisien de $X$. 
\end{enumerate}
\end{prop}

\begin{defi}[Revêtements tempérés modérés]\label{p' version du groupe fondamental}
Notons $\Cov^{\tp, \pp}(X)$ la sous-catégorie pleine de $\Cov^\tp(X)$ constituée des revêtements tempérés qui sont quotients d'un revêtement topologique d'un revêtement galoisien fini \emph{de degré premier à $p$}. Si $x\in X$ est un point géométrique, notons $\pi_1^{\tp, \pp}(X,x)$ le groupe d'automorphismes du foncteur fibre en $x$ :  $$\mathrm{F}_{x}^{\pp} : \Cov^{\tp, \pp}(X)\to \mathrm{Ens}$$ qui à un revêtement $Y\to X$ associe la fibre $Y_x$ au-dessus de $x$. 
La classe d'isomorphisme de $\pi_1^{\tp(X,x), \pp}$ ne dépend pas du point géométrique $x$ considéré, de sorte que l'on pourra considérer le groupe $\pi_1^{\tp(X),\pp}$, défini à unique isomorphisme extérieur près, appelé le \emph{groupe fondamental tempéré modéré de $X$}.
\end{defi}

\begin{rem}
Munissons $\Cov^{\tp, \pp}(X)$ de la classe d'isomorphisme $\Upsilon_X^{\pp}$ des foncteurs $\Cov^{\tp, \pp}(X)\to \mathfrak{Ens}^\mathrm{dd}$ contenant les foncteurs fibres $F_x^{\pp}$ au dessus des points géométriques (ils sont tous isomorphes). Alors $(\Cov^{\tp, \pp}(X),\Upsilon_X^{\pp})$ est un tempéroïde connexe dont $\pi_1^{\tp, \pp}(X)$ est le groupe tempéré, de sorte que l'on a un isomorphisme de tempéroïdes connexe : $$\left(\Cov^{\tp, \pp}(X),\Upsilon_X^{\pp}\right)\simeq \left(\mathcal{B}^{\tp}(\pi_1^{\tp, \pp}(X)),\Upsilon^{\mathrm{oubli}}_{\pi_1^{\tp, \pp}(X)}\right).$$
\end{rem}

\begin{rem}\label{pi1t}
Lorsque $X$ est une courbe $K$-analytique, le groupe $\pi_1^{\tp, \pp}(X)$ peut se construire directement à partir de $\pi_1^\tp(X)$ : en effet, $\pi_1^{\tp, \pp}(X)$ n'est autre la limite projective des quotients de $\pi_1^\tp(X)$ qui ont un sous-groupe distingué sans torsion d'indice fini premier à $p$. Ce résultat est en revanche faux en dimension supérieure, lorsque les groupes topologiques fondamentaux des espaces $k$-analytiques ne sont plus libres. 
\end{rem}

\section{Géométrie anabélienne dans les courbes analytiques}

Dans toute la suite du texte, fixons un corps $k$ ultramétrique complet non trivialement valué, algébriquement clos, et $p$ l'exposant caractéristique de $\widetilde{k}$. Fixons un système de racines primitives l'unité, de sorte que l'on ait des isomorphismes naturels : $\mu_n(k)\simeq \Z/n\Z$ et $\ds{\varprojlim_{n\in\N^*}\mu_n(k)=\widehat{\Z}}$ (plutôt que de parler de $\widehat{\Z}(1)$).

\subsection{Description du graphe d'anabélioïdes d'une triangulation généralisée} \label{section description}

\subsubsection{Comportement des composantes coronaires dans le cas d'une vraie triangulation}\label{partie comportement des composantes coronaires}

Nous appellerons \emph{kummérien} tout revêtement étale fini $\varphi : \CC'\rightarrow \CC$ entre deux couronnes $k$-analytiques tel qu'il existe un certain entier $\ell$ premier à $p$ permettant de considérer $\CC'$ comme un $\mu_\ell$-torseur de Kummer de $\mathcal{C}$ via $\varphi$. Rappelons que le choix d'une orientation de $\mathcal{C}$ permet d'exhiber un isomorphisme canonique : $\Kum_\ell(\mathcal{C})\iso \Z/\ell\Z$.\\

Soit $\widehat{\Z}^{\pp}$ le complété pro-$p'$-fini de $\Z$, égal à $\displaystyle{\prod_{q\in \P\setminus \lbrace p\rbrace}\Z_q}$.\\

\begin{prop}[Description du groupe fondamental modéré des disques et des couronnes, \cite{Ber2}, Théorèmes $6.3.2$ et $6.3.5$]\label{proposition caractere kummerien et trivial}\

\begin{enumerate}
\item Tout revêtement fini modéré d'un disque est trivial. De là, si $\mathcal{D}$ est un disque $k$-analytique, alors $\pi_1^\ttt(\mathcal{D}) = \lbrace 1\rbrace$;
\item Tout revêtement fini galoisien modéré d'une couronne de type $\CC(I)$ avec $I$ intervalle compact de $\R_+^*$ est kummérien. De là, si $\CC$ est une telle couronne $k$-analytique, alors ${\pi_1^\ttt(\CC) = \widehat{\Z}^{\pp}}$.
\end{enumerate}
\end{prop}\

Si $\CC'\rightarrow \CC$ est un plongement de couronnes $k$-analytiques compactes permettant de considérer $\CC'$ comme une sous-couronne de $\CC$, alors comme on a d'une part pour tout entier $\ell$ premier à $p$ un isomorphisme $\Kum_\ell(\CC')\simeq\Kum_\ell(\CC)$ (voir Proposition \ref{propriétés kummériennes des revêtements de couronnes}) et que de l'autre tout revêtement galoisien modéré de $\CC$ ou $\CC'$ est kummérien, alors le foncteur naturel associé $\Covt(\CC')\rightarrow\Covt(\CC)$ est une équivalence de catégories. Par conséquent, pour n'importe quelle couronne $k$-analytique $\CC$ (en particulier pour toute couronne ouverte) l'on a : $$\pi_1^\ttt(\CC) = \widehat{\Z}^{\pp}.$$ 

Si $S$ est une triangulation généralisée de $X$, les composantes $\GG_e$ du graphe d'anabélioïdes $\GG(X,S)$ données par les arêtes correspondant à des couronnes, sont donc isomorphes (en tant qu'anabélioïdes connexes) à $\mathcal{B}(\widehat{\Z}^{\pp}).$\\

\begin{rem}
En reprenant le raisonnement précédent restreint aux $\mu_\ell$-torseurs on obtient pour toute couronne $\mathcal{C}$ l'isomorphisme : 
$$ \HH^1(\mathcal{C}_\et, \mu_\ell) \simeq \Kum_\ell (\CC)\simeq \Z/\ell\Z . $$ 
\end{rem}

\subsubsection{Triangulation généralisée et rabotée d'une courbe}\label{sous-section rabotée d'une courbe}

Expliquons pourquoi la description du semi-graphe d'anabélioïdes $\GG(X, S)$ associé à une triangulation généralisée $S$ d'une courbe $k$-analytique $X$ peut se ramener, quitte à restreindre $X$ en une courbe que l'on appellera une \emph{rabotée} de $X$, à supposer que $S$ est une vraie triangulation de $X$. Reprenons l'exemple esquissé à la remarque \ref{remarque sur l'ensemble des noeuds qui peut ne pas être une triangulation} : celui d'une courbe $k$-analytique $Y$ quasi-lisse de squelette analytique non vide, à laquelle on enlève un point $x$ de type $4$. Enlever $x$ revient à rajouter l'arête ouverte $]x, r(x)]$ : $S^\an(Y\setminus \lbrace x\rbrace)=S^\an(Y)\cup \; ]x, r(x)]$. Or $r^{-1}(]x, r(x)])$ n'a aucune raison d'être une couronne. Néanmoins, tout point $\widetilde{x}\in ]x, r(x)]$ est de type $2$ ou $3$, et $r^{-1}(]\widetilde{x}, r(x)])$ a le bon goût d'être une couronne $k$-analytique d'après \cite{Duc} ($5.1.12.3$). Cela nous amène à définir la notion de \emph{rabotée d'une courbe}.

 \begin{defi}[Rabotée d'une courbe]\label{définition d'une rabotée} 
 Soit $X$ une courbe $k$-analytique connexe quasi-lisse dont $S$ est une triangulation généralisée non vide de squelette $\Gamma$. Notons $r$ la rétraction canonique de $X$ sur $\Gamma$. Soit $\mathscr{I}$ l'ensemble (éventuellement vide) des composantes connexes de $\Gamma\setminus S$ non relativement compactes dans $\Gamma$ dont la préimage par $r$ n'est pas une couronne $k$-analytique. Il existe une famille $(x_I)_{I\in \mathscr{I}}$ de $S$ telle que $I$ est un intervalle ouvert non relativement compact aboutissant à $x_I$, pour tout $I\in \mathscr{I}$.
Une \emph{rabotée de $X$ relativement à $S$} sera définie comme une courbe $k$-analytique $X^\circ$ obtenue à partir de $X$ en ne gardant de chaque $r^{-1}\left( I\right)$ que les couronnes $r^{-1}\left( \;]x_I, \widetilde{x_I}[\; \right)$ : 
\begin{align*}
X^\circ&=\left( X\setminus \bigcup_{I\in \mathscr{I}}r^{-1}\left( I\right)\right) \cup \bigcup_{I\in \mathscr{I}}r^{-1}\left( \;]x_I, \widetilde{x_I}[\; \right)\\
&=X\setminus \bigcup_{I\in \mathscr{I}}r^{-1}\left( I\;\setminus \;]x_I, \widetilde{x_I}[\; \right),
\end{align*}
 
 où $(\widetilde{x_I})_{I\in \mathscr{I}}$ est une collection de points de $\Gamma$ vérifiant $\widetilde{x_I}\in I$ pour tout $I\in \mathscr{I}$. 
 \item Si l'ensemble $\Sigma_X$ des nœuds du squelette analytique $S^\an(X)$ est non vide, on désignera par \emph{rabotée de $X$} toute rabotée de $X$ relativement à la triangulation généralisée $\Sigma_X$.
  \end{defi}
  
  \begin{rem}
  La courbe $X$ est une rabotée d'elle-même si et seulement si $S$ constitue une triangulation de $X$. La notion de \emph{rabotée} n'est donc utile que dans le cas contraire. Remarquons qu'une rabotée $X^\circ$ de $X$ relativement à $S$ contient toujours $S$ en tant qu'ensemble des sommets du squelette de la triangulation de $X^\circ$ définie par $S$, et que la donnée d'une rabotée de $X$ est équivalente au choix de la famille $\left(\widetilde{x_I}\right)_{I\in \mathscr{I}}$. 
  
  \end{rem}

  \begin{prop}\label{équivalence des semi-graphe de la rabotée et de la non rabotée} 
  Soit $X$ une courbe $k$-analytique connexe quasi-lisse munie d'une triangulation généralisée $S$ non vide de squelette $\Gamma$. Notons $X^\circ$ une rabotée de $X$ relativement à $S$. Alors les deux semi-graphes d'anabélioïdes $\GG(X^\circ,S)$ et $\GG(X,S)$ sont isomorphes. Ainsi, tous les groupes des arêtes (que celles-ci soient ouvertes ou fermées) sont isomorphes à $\widehat{\Z}^{\pp}$.
   \end{prop}
   
   \begin{proof}
   L'application $\mathcal{I}\mapsto r(\mathcal{I})$ établit une bijection entre les composantes connexes de $X\setminus S$ qui ne sont ni des disques ni des couronnes, et l'ensemble $\mathscr{I}$ (cf. définition \ref{définition d'une rabotée}) des composantes connexes de $\Gamma\setminus S$ non relativement compactes dans $\Gamma$ dont la préimage par $r$ n'est pas une couronne. Soit $(\widetilde{x_I})_{I\in \mathscr{I}}$ la famille de $\Gamma$ permettant de définir $X^\circ$. Eu égard à la construction de $\GG(X^\circ,S)$ et $\GG(X,S)$, il suffit de montrer que pour toute arête ouverte $e$ de $\GG(X,S)$ (telle que la composante $\mathcal{I}_e$ corresponde à $I_e\in \mathscr{I}$), le foncteur $\Covt(\mathcal{I}_e)\rightarrow \Covt \left(r^{-1}\left(\;]x_{I_e},\widetilde{x_{I_e}} [ \;\right)\right)$ induit naturellement par l'injection $r^{-1}\left(\;]x_{I_e},\widetilde{x_{I_e}} [ \;\right)\hookrightarrow \mathcal{I}_e$ est une équivalence de catégories. 
   Comme $r^{-1}\left(\;]x_{I_e},\widetilde{x_{I_e}} [ \;\right)$ est une couronne et que tout revêtement galoisien modéré d'une couronne est kummérien, on a $$\pi_1^\ttt\left(r^{-1}(\;]x_{I_e},\widetilde{x_{I_e}} [ \;)\right)= \widehat{\Z}^{\pp}.$$ Soit $\left(J_k=]x_{I_e}, y_k[\right)_{k\in \N}$ une suite exhaustive de sous-intervalles emboités de $I_e$, avec $J_k\subseteq J_{k+1}$. Pour tout $k\in \N$, notons $\CC_k=r^{-1}(J_k)$ :  c'est une couronne $k$-analytique. Remarquons que $\CC_k$ est une sous-couronne de $\CC_{k'}$ dès que $k<k'$, d'où l'équivalence de catégories $\Covt(\CC_k)\iso \Covt(\CC_k')$. Comme chaque $\Covt(\CC_k)$ est équivalente à $\BB\left(\widehat{\Z}^{\pp}\right)$ et que l'on a : $$\mathcal{I}_e = \varinjlim_{k\in \N} \CC_k, $$ alors on obtient en passant à la limite inductive $\Covt(\mathcal{I}_e)\simeq \BB\left(\widehat{\Z}^{\pp}\right)$, d'où l'équivalence souhaitée $\Covt(\mathcal{I}_e)\iso \Covt \left(r^{-1}\left(\;]x_{I_e},\widetilde{x_{I_e}} [ \;\right)\right)$. 
   
   \end{proof}
   
   La proposition précédente nous indique qu'étudier le semi-graphe d'anabélioïdes d'une courbe munie d'une triangulation généralisée revient à étudier le semi-graphe d'anabélioïde d'une courbe munie d'une vraie triangulation. Par conséquent, toutes les arêtes de $\GG(X,S)$ se comportent comme si elles étaient associées à de vraies couronnes en ce qu'elles ont chacun pour groupe associé $\widehat{\Z}^{\pp}$. Il est ainsi justifiable d'appeler \emph{composantes coronaires} d'une triangulation généralisée toutes les arêtes (ou les composantes associées aux arêtes) de son semi-graphe d'anabélioïdes.  
   
   \begin{rem}
   On peut montrer que le morphisme $\pi_1^\tp(X^\circ)\rightarrow \pi_1^\tp(X)$ induit par le plongement $X^\circ\hookrightarrow X$ induit un isomorphisme de groupes : $$\pi_1^{\tp, \pp}(X^\circ)\iso \pi_1^{\tp, \pp}(X).$$ Ainsi, une courbe et sa rabotée sont indistinguables sous le prisme des revêtements (tempérés) modérés. En revanche ils n'auront plus du tout les mêmes $p$-revêtements, et le morphisme  $\pi_1^\tp(X^\circ)\rightarrow \pi_1^\tp(X)$ n'a aucune chance d'être un isomorphisme dès que ${X}^\circ\varsubsetneq X$. 
   
  \end{rem}

\subsubsection{Germes d'espaces $k$-analytiques}

Commençons par rappeler la définition de la catégorie des \emph{germes d'espaces $k$-analytiques} :  \

Un \emph{germe $k$-analytique} est une paire $(X,S)$ où $X$ est un espace $k$-analytique et $S$ un sous-espace de l'espace topologique sous-jacent $\vert X\vert$. Si $x\in X$, le germe $(X,\lbrace x\rbrace)$ est simplement noté $(X,x)$. Soit $k$-$\mathscr{G}$ la catégorie dont les objets sont les germes d'espaces $k$-analytiques et les morphismes entre deux germes $(Y,T)$ et $(X,S)$ sont les morphismes d'espaces analytiques $\varphi : Y\to X$ tels que $\varphi(T)\subseteq S$. La catégorie des germes $k$-analytiques, notée $k$-$\mathit{Germ}$, est définie comme la localisée de la catégorie $k$-$\mathscr{G}$ par la famille des morphismes $\varphi : (Y,T)\to (X,S)$ induisant un isomorphisme de $Y$ vers un voisinage ouvert de $S$ dans $X$. Cette famille admet un \og calcul des fractions à droite \fg{}, se sorte que dans cette catégorie l'ensemble des morphismes $\Hom_{k-\mathit{Germ}}((Y,T),(X,S))$ correspond à la limite inductive de l'ensemble des morphismes $k$-analytiques $\varphi : \mathscr{V}\to X$ où $\mathscr{V}$ décrit les voisinages ouverts de $T$ dans $Y$, et $\varphi(T)\subseteq S$. Un tel morphisme $k$-analytique $\varphi : \mathscr{V}\to X$ correspondant à un morphisme de germes $\psi : (Y,T)\to (X,S)$ est un \emph{représentant} de $\psi$. Tout morphisme de germes $k$-analytiques admet un représentant, et un morphisme $\psi : (Y,T)\to (X,S)$ est un isomorphisme dans $k$-$\mathit{Germ}$ si et seulement s'il admet un représentant induisant un isomorphisme $k$-analytique entre des voisinages ouverts de $T$ et $S$.

\bigskip

\begin{defi} Un morphisme de germes $k$-analytiques $\psi: (Y,T)\rightarrow (X,x)$ sera appelé un \emph{revêtement modéré fini} (resp. \emph{revêtement fini étale}) \emph{de $(X,x)$} s'il admet au moins un représentant $\varphi : \mathscr{V}\to X$ tel que $\mathscr{V}\to\varphi(\mathscr{V})$ est un revêtement modéré fini ((resp. fini étale) d'espaces $k$-analytiques, et $T=\varphi^{-1}(x)$. Nous noterons par la suite $\Covt((X,x))$ (resp. $\Cov^\fet((X,x))$) les catégories respectives de ces revêtements du germe $(X,x)$.\

Si $X$ désigne notre courbe $k$-analytique quasi-lisse de départ, $S$ une triangulation généralisée et $s\in S$, remarquons que l'on a, au sens $2$-catégorique :

$$\Covt((X,s)):=\varinjlim_U \Covt(U)$$ où $U$ décrit l'ensemble des voisinages ouverts de $s$.

Si $\CC\in \CX$ et $\mathcal{D}\in \DX$, définissons de même : 
\begin{align*}
\Covt(\CC\cap (X,s)) & :=\varinjlim_U \Covt(\CC\cap U) \\
\Covt(\mathcal{D}\cap (X,s)) & :=\varinjlim_U \Covt(\mathcal{D}\cap U).
\end{align*}
\end{defi}

La proposition suivante signifie que les comportements locaux des revêtements modérés autour de $s$ sont équivalents aux produits d'extensions finies du corps résiduel complété $\HC(s)$ dominées par une extension galoisienne de degré premier à $p$.
\begin{prop}
Le foncteur naturel $\Covt((X,s))\rightarrow \Covt(\HC(s))$ est une équivalence de catégories, où $\Covt(\HC(s))$ désigne la catégorie des schémas finis étales au-dessus  de $\Spec(\HC(s))$ dominés par un revêtement galoisien de degré premier à $p$.
\end{prop}

\begin{proof}
Cela découle directement du théorème $3.4.1$ de \cite{Ber2} qui nous indique que le foncteur naturel $\Cov^\fet((X,s))\rightarrow \Cov^\fet(\HC(s))$ est une équivalence de catégories, où $\Cov^\fet(\HC(s))$ désigne la catégorie des schémas finis étales au-dessus  de $\Spec(\HC(s))$.
\end{proof}\

\begin{lem}\label{lemme sur l'équivalence de catégories avec une sous-couronne} 
Considérons deux domaines analytiques $U_1$ et $U_2$ de $X\setminus S$ qui sont des couronnes $k$-analytiques aboutissant à $s$, avec $U_1 \subseteq U_2$. Alors le plongement ${U_1\hookrightarrow U_2}$ induit une équivalence de catégories : $\Covt(U_2) \iso \Covt(U_1)$.
\end{lem}

\begin{proof}
De manière générale, il se peut qu'une couronne analytique $U_1$ contenue dans une deuxième couronne analytique $U_2$ ne soit pas une sous-couronne de cette dernière (voir définition \ref{définition d'une sous-couronne})). Les domaines analytiques de $\P_k^{1,\an}$ donnés par les conditions ${\vert T-p \vert \in ]0,\frac{1}{p}[}$ et ${\vert T\vert \in ]0,1[}$ en sont des exemples. Néanmoins, $U_1$ est une sous-couronne de $U_2$ dès que l'intersection des squelettes $S^\an(U_1)\cap S^\an(U_2)$ est non vide. Or, dans le cadre de ce lemme, l'inclusion $U_1 \subseteq U_2$ et le fait que $U_1$ et $U_2$ aient un bout en commun impliquent qu'il existe $x_1, x_2$ et $y\in X$ tels que $S^\an(U_1)=]s, x_1[$, $S^\an(U_2)=]s, x_2[$ et $S^\an(U_1)\cap S^\an (U_2)=]s,y[\neq \emptyset$. Ainsi $U_1$ est une sous-couronne de $U_2$. D'après la proposition \ref{propriétés kummériennes des revêtements de couronnes} combinée au caractère kummérien des revêtements modérés d'une couronne, on a l'équivalence ${\Covt(U_1)\iso \Covt(U_2)}$.

\end{proof}

\begin{lem}\label{lem1} 
Pour toute composante $T\in \mathsf{Comp}(X,S,s)$ la catégorie $\Covt(T\cap (X,s))$ est équivalente à la catégorie $\BB( \widehat{\Z}^{\pp})$ des $\widehat{\Z}^{\pp} $-ensembles finis, et si $T\in \CC(X,S,s)$ le foncteur de restriction $\Covt(T)\rightarrow \Covt(T \cap (X,s))$ est une équivalence de catégorie.
\end{lem}

\begin{proof}
D'après la sous-section \ref{sous-section rabotée d'une courbe} on peut supposer que $S$ est une vraie triangulation de $X$. Soit $T\in \mathsf{Comp}(X,S,s)$, et soit $(U_i)_i$ une base de voisinages ouverts de $s$ telle que chaque $U_i\cap T$ soit une couronne (une telle base d'ouverts existe toujours). Si $U_i\subset U_j$ alors $U_i\cap T\rightarrow U_j \cap T$ est un plongement de couronnes, d'où l'on déduit d'après le lemme \ref{lemme sur l'équivalence de catégories avec une sous-couronne} précédent que le foncteur de restriction $\Covt(U_j \cap T)\rightarrow \Covt(U_i \cap T)$ est une équivalence de catégories. Par ailleurs ces équivalences sont compatibles avec l'identification des $\Covt(U_i \cap T)$ à  $\BB( \widehat{\Z}^{\pp})$. On obtient donc en passant à la limite inductive que $\Covt((X,s)\cap T)$ est elle-même équivalente à  $\BB( \widehat{\Z}^{\pp})$, ainsi que l'équivalence lorsque $T$ est une couronne. 
\end{proof}\

\begin{prop}\label{allure locale} 
Soit $(X, S=\lbrace s\rbrace)$ une courbe $k$-analytique connexe quasi-lisse munie d'une triangulation possédant un unique sommet, et dont le squelette ne contient aucune arête qui boucle sur $s$ :\

\begin{enumerate}
\item Soit $U$ un voisinage ouvert de $s$ totalement découpé autour de $s$. Alors le foncteur naturel $\Covt(X)\rightarrow \Covt(U)$ induit par l'injection $U\hookrightarrow X$ est pleinement fidèle et son image essentielle est constituée des objets de $\Covt(U)$ qui sont triviaux sur tout disque de $\mathsf{Comp}(X,S,s)$, i.e. dont l'image par le foncteur de restriction $\Covt(U)\rightarrow \Covt(\mathcal{D}\cap U)$ est triviale pour tout $\mathcal{D}\in \DX$.
\item Le foncteur de restriction $\Covt(X)\rightarrow \Covt((X,s))$ est pleinement fidèle et son image essentielle est constituée des objets de $\Covt((X,s))$ qui sont triviaux sur tout disque de $\mathsf{Comp}(X,S,s)$, i.e. dont l'image par le foncteur de restriction $\Covt((X,s))\rightarrow \Covt( \mathcal{D}\cap (X,s))$ est triviale pour tout $\mathcal{D}\in \DX$.
\end{enumerate}
\end{prop}

\begin{proof}\

\begin{enumerate}
\item Considérons le recouvrement ouvert de $X$ suivant :  $\mathcal{U}:= U \cup (T)_{T\in \mathsf{Comp}(X,S,s)}$.\

Une \emph{donnée de descente sur $\UU$} consiste en la donnée de $(Y_U, (Y_T, \varphi_T)_T)$ où $Y_U\in \Covt(U)$, $Y_T\in \Covt(T)$ et pour toute composante $T\in \mathsf{Comp}(X,S,s)$ un isomorphisme $$\varphi_T  : \restriction{Y_U}{U\cap T}\iso \restriction{Y_T}{U\cap T}.$$ Les conditions usuelles sur les $1$-cocycles sont ici tautologiques puisque l'intersection de trois des composantes ouvertes de $\UU$ est automatiquement vide. Par descente topologique nous savons que la catégorie $\Covt(X)$ est équivalente à la catégorie des données de descente sur $\UU$. Comme par ailleurs le foncteur de restriction $\Covt(T)\rightarrow \Covt(U\cap T)$ est toujours pleinement fidèle (si $T$ est un disque cela découle de la trivialité de tout revêtement modéré sur $T$, et si $T$ est une couronne cela découle du lemme \ref{lemme sur l'équivalence de catégories avec une sous-couronne}), alors le foncteur $\Covt(X)\rightarrow \Covt(U)$ est pleinement fidèle. 

Son image essentielle consiste en la donnée des éléments de $\Covt(U)$ dont l'image par le foncteur de restriction $\Covt(U)\rightarrow \Covt(U\cap T)$ est dans l'image essentielle de $\Covt(T)\rightarrow \Covt(U\cap T)$ pour toute composante $T\in \mathsf{Comp}(X,S,s)$. Or nous savons d'après \ref{lemme sur l'équivalence de catégories avec une sous-couronne} que si $T$ est une couronne ce dernier foncteur est une équivalence de catégories et que si $T$ est un disque tout revêtement modéré de $T$ est trivial, de là l'on déduit le résultat.
\item Soit $(U_i)_i$ la famille de voisinages ouverts de $s$ totalement découpés autour de $s$. Nous savons que cette famille forme une base de voisinages ouverts de $s$. Par le premier point démontré ci-dessus l'on sait que si $U_i\subset U_j$ le foncteur $\Covt(U_j)\rightarrow \Covt(U_i)$ est pleinement fidèle. En passant à la limite inductive on obtient que tous les foncteurs $\Covt(U_j)\rightarrow \Covt((X,s))$ sont pleinement fidèles. Comme les $\Covt(X)\rightarrow \Covt(U_i)$ sont eux-mêmes pleinement fidèles d'après le premier point, alors par transitivité le foncteur $\Covt(X)\rightarrow \Covt((X,s))$ est lui-même pleinement fidèle.

Cherchons son image essentielle. Si $Z\in \Covt((X,s))$ est dans l'image essentielle, alors il découle de la trivialité de tout revêtement modéré d'un disque que sa restriction à $\DD\cap (X,s)$ est triviale pour tout $\DD\in \DX$. Soit réciproquement $Y\in \Covt((X,s))$ tel que sa restriction à $\DD\cap (X,s)$ soit triviale pour tout $\DD\in \DX$. Il en existe alors un représentant $Y_i\in \Covt(U_i)$ pour un certain voisinage $U_i$ totalement découpé autour de $s$. Si $\DD\cap U_i$ est un disque, alors la restriction de $Y_i$ à  $\DD\cap U_i$ est automatiquement triviale, et si $\DD\cap U_i$ est une couronne nous savons d'après le lemme \ref{lem1} que le foncteur de restriction $\Covt(\DD\cap U_i)\rightarrow \Covt(\DD\cap (X,s))$ est une équivalence de catégorie, donc la restriction de $Y_i$ à  $\DD\cap U_i$ est également triviale. Par le premier point de la proposition cela nous indique que $Y_i$ est dans l'image essentielle de $\Covt(X)\rightarrow \Covt(U_i)$, ce qui nous permet de conclure que $Y$ est bien dans l'image essentielle de $\Covt(X)\rightarrow \Covt((X,s))$, ce qui termine la preuve.
\end{enumerate}
\end{proof}

\subsubsection{Courbe résiduelle,  valuation associée à une branche et groupe d'inertie modérée autour d'un point de type 2}\label{courbe résiduelle}

Supposons que $x$ est un point de type $2$ de $X$. Alors le corps résiduel $\widetilde{\HC(x)}$ est une extension de $\widetilde{k}$ de type fini et de degré de transcendance $1$. Ainsi $\widetilde{\HC(x)}$ est de la forme $\widetilde{k}(\mathscr{C}_x)$ où $\mathscr{C}_x$ est une $\widetilde{k}$-courbe projective, lisse et intègre qui est uniquement déterminée par $x$ et que l'on appellera la \emph{courbe résiduelle en $x$}. Rappelons que le \emph{genre de $x$}, noté $g(x)$, est défini comme le genre de la courbe résiduelle $\mathscr{C}_x$ lorsque $x$ est de type $2$ (il est égal à $0$ sinon). L'espace (topologique) de Zariski-Riemann $\PP_{\widetilde{\HC(x)}/\widetilde{k}}$ des valuations sur $\widetilde{\HC(x)}$ triviales sur $\widetilde{k}$ est alors homéomorphe à l'espace topologique sous-jacent au schéma $\mathscr{C}_x$ : au point générique correspond la valuation triviale et à chaque point fermé correspond une valuation discrète.\\

Dans \cite{Duc}, $(4.2.6)$, l'auteur associe à chaque branche $b\in \br(X,x)$ une valuation $\langle.\rangle_b$ sur $\widetilde{\HC(x)}$ triviale sur $\widetilde{k}$, ie. un élément de $\PP_{\widetilde{\HC(x)}/\widetilde{k}}$. En composant $\langle.\rangle_b$ avec la valuation $\vert.\vert_x$ définie naturellement sur $\HC(x)$ l'on obtient une valuation $\vert.\vert_b$ sur $\HC(x)$ qui raffine par construction $\vert.\vert_x$ et coïncide avec celle-ci sur $k$.\

En notant $\langle.\rangle_0$ la valuation triviale sur $\widetilde{\HC(x)}$, Ducros déduit de la théorie de la réduction de Temkin (\cite{Tem}) que l'application $b\mapsto\langle.\rangle_b$ établit une injection $$\br(X,x)\hookrightarrow \PP_{\widetilde{\HC(x)}/\widetilde{k}}\setminus\lbrace \langle.\rangle_0\rbrace$$ dont l'image est de complémentaire fini, et qui est bijective si et seulement si $x$ appartient à l'intérieur analytique de $X$. La réunion de l'image et de la valuation triviale correspond à la réduction de Temkin du germe $(X,x)$. \\

\begin{rem}
Ducros établit en fait un résultat un peu plus général en faisant usage de valuations graduées. Néanmoins, comme $x$ est de type $2$, nous n'en avons pas besoin dans ce texte, et toutes les valuations considérées sont des valuations classiques. 
\end{rem}\

Puisque $x$ est de type $2$, ce dernier résultat se reformule comme suit : l'application qui à une branche $b$ associe la valuation $\langle.\rangle_b$ sur $\widetilde{\HC(x)}=\widetilde{k}(\mathscr{C}_x)$ induit une injection dont l'image est de complémentaire fini de $\br(X,x)$ dans l'ensemble $\PP_{\widetilde{k}(\mathscr{C}_x)/\widetilde{k}}\setminus\lbrace \langle.\rangle_0\rbrace$. Comme ce dernier s'identifie à l'ensemble des points fermés de $\mathscr{C}_x$, cela revient à une bijection entre $\br(X,x)$ et l'ensemble des points fermés d'un ouvert de Zariski $\mathscr{U}_x$ de $\mathscr{C}_x$, égal à $\mathscr{C}_x$ tout entier si et seulement si $x$ est dans l'intérieur analytique de $X$.\\

Soit $b\in \br(X,x)$ une branche issue de $x$. Dans \cite{Duc} l'auteur définit l'anneau associé à $b$, noté $\mathscr{O}_X(b)$, par la limite inductive : 
$$\mathscr{O}_X(b)=  \varinjlim_V \OO_X(V)  $$ où $V$ décrit l'ensemble des sections de $b$. Il existe un morphisme naturel de $\OO_{X,x}$ dans $\OO_X(b)$. En tant que courbe quasi-lisse notre courbe $X$ est normale. En conséquence l'anneau local $\OO_{X,x}$ est artinien et intègre, c'est donc un corps qui coïncide alors avec son corps résiduel $\kappa (x)$. Par ailleurs toute section $Z$ de $b$ est un espace non vide, connexe et normal, par conséquent $\OO_X(Z)$ est intègre, et par passage à la limite on en déduit que $\OO_X(b)$ lui-même est intègre. Le corps $\kappa (x)$ coïncidant avec $\OO_{X,x}$, il se plonge dans $\OO_X(b)$.\\

\begin{defi}Notons $\OO_X(b)_\sep$ la fermeture intégrale séparable de $\kappa (x)$ dans l'anneau $\OO_X(b)$; c'est un corps naturellement muni d'une valuation prolongeant la valuation $\vert . \vert_b$ sur $\kappa (x)$ (qui est elle-même la restriction à $\kappa (x)$ de la valuation $\vert . \vert_b$ définie sur $\HC(x)$, $\kappa(x)$ étant dense dans $\HC(x)$ car $X$ est un bon espace) que nous continuerons à noter $\vert . \vert_b$ et qui jouit de la propriété suivante : \\

\end{defi} 

\begin{prop}[\cite{Duc08}, Proposition $3.6$ ; \cite{Duc}, Théorème $4.2.19$]\label{corps hensélisé par rapport à la valuation d'une branche} \

Le corps valué $(\OO_X(b)_\sep, \vert . \vert_b)$ est le hensélisé de $(\kappa(x), \vert . \vert_b)$.
\end{prop}

En tant que hensélisé $(\OO_X(b)_\sep, \vert . \vert_b)$ est une extension \emph{immédiate} de $(\kappa(x), \vert . \vert_b)$, i.e. ces deux corps valués ont même groupe de valeurs et même corps résiduel : \\

\begin{itemize}
\item[•] \emph{groupe des valeurs : }
Si l'on note $\CC_e$ la couronne correspondant à l'arête $e$ associée à la branche $b$, le Théorème $4.3.5$ de \cite{Duc} nous indique qu'il existe une fonction $\lambda \in \kappa (x)=\OO_{X,x}$ définie sur $\CC_e$, qui en est une fonction coordonnée, telle que $\vert \lambda \vert_b$ soit infiniment proche de $1$ inférieurement, et telle que  : 
$$ \vert \kappa (x)^\times \vert_b=\vert \OO_X(b)_\sep^\times \vert _b= \vert k^\times \vert\oplus \vert \lambda \vert_b^\Z $$
\item[•] \emph{corps résiduel : }
Les corps résiduels associés $\widetilde{\kappa (x)}_b$ et $\widetilde{\OO_X(b)_\sep}$ sont égaux à $\widetilde{k}$, et sont en particulier algébriquement clos.
\end{itemize}\

Notons $\OO_X(b)_\sep^\s$ une clôture séparable du corps $\OO_X(b)_\sep$. Comme $\OO_X(b)_\sep$ est hensélien la valuation sur $\OO_X(b)_\sep$ se prolonge en une unique valuation sur $\OO_X(b)_\sep^\s$ encore notée $\vert . \vert_b$. Notons $\mathcal{I}_b$ et $\mathcal{W}_b$ les groupes d'inertie et de ramification de l'extension de corps valués $ \OO_X(b)_\sep^\s/\OO_X(b)_\sep$, et $\mathcal{I}_b^\ttt= \mathcal{I}_b/\mathcal{W}_b$ son \emph{groupe d'inertie modérée}. Nous savons alors  que l'on dispose de l'isomorphisme suivant : 
$$\mathcal{I}_b^\ttt = \Hom(\vert \OO_X(b)_\sep^{\s \; \times} \vert _b/ \vert \OO_X(b)_\sep^\times \vert _b , \widetilde{\OO_X(b)_\sep}^\times) $$ que l'on peut récrire de la manière suivante sachant que  $\vert \OO_X(b)_\sep^\times \vert _b= \vert k^\times \vert\oplus \vert \lambda \vert_b^\Z $ et que $\vert k^\times \vert $ est divisible : 
$$\mathcal{I}_b^\ttt = \Hom(\vert\lambda \vert _b^\Q/ \vert\lambda \vert _b^\Z,  \widetilde{k}^\times  ) = \Hom(\Q/ \Z, \widetilde{k}^\times), $$ et comme $\widetilde{k}$ est algébriquement clos de caractéristique $p$ on obtient  :  $\mathcal{I}_b^\ttt \simeq \widehat{\Z}^{\pp}  $.\\

Si ${\kappa(x)_b^\s}$ désigne une clôture séparable de ${\kappa(x)_b}$ munie d'une valuation (encore notée $\vert . \vert_b$) prolongeant $\vert . \vert _b$ sur $\kappa(x)_b$), notons $\mathcal{I}_{\kappa(x)_b}^\ttt$ le groupe d'inertie modérée de l'extension ${\kappa(x)_b^\s}/{\kappa(x)_b}$. Comme $\OO_X(b)_\sep$ est le hensélisé de $\kappa(x)_b$ nous avons l'égalité des groupes d'inertie modérée : $$\mathcal{I}_b^\ttt=\mathcal{I}_{\kappa(x)_b}^\ttt\simeq  \widehat{\Z}^{\pp} $$\

\begin{rem}Pour tout entier $n$ premier à $p$, le seul quotient d'ordre $n$ de $\widehat{\Z}^{\pp} $ correspond à l'extension de $\OO_X(b)_\sep$ obtenue en extrayant une racine $n$-ième d'un élément de norme $\vert . \vert_b$ égale à  $\vert \lambda \vert_b$.

\end{rem}

\subsubsection{Comportement des composantes sommitales}\label{cmportement des composantes sommitales} 

\begin{prop}[\cite{Duc}, Théorème $4.3.13$]\label{allureD} 
Supposons que $\varphi : Y\rightarrow X$ est un morphisme étale entre deux courbes $k$-analytiques quasi-lisses, $y\in Y_{[2]}$ et $x\in X_{[2]}$ son image. Soit $b\in \br(Y,y)$ d'image $a\in \br(X,x)$. Soient $\mathscr{C}_y$ et $\mathscr{C}_x$ les courbes résiduelles de $y$ et $x$, et $\mathscr{Q}$ (resp. $\mathscr{P}$) le point fermé de $\mathscr{C}_y$ (resp. $\mathscr{C}_x$) correspondant à $b$ (resp. $a$) par la bijection précédente. Alors $\mathscr{P}$ est l'image de $\mathscr{Q}$ par $\mathscr{C}_y\rightarrow \mathscr{C}_x$, et en notant $e$ l'indice de ramification en $\mathscr{Q}$, on a  : $e=\deg(b\rightarrow a)$.

\end{prop}

\begin{rem}
En reprenant les notations de la proposition précédente, comme $x\in X_{[2]}$ nous savons que le groupe $\vert \HC(x)^\times\vert/\vert k^\times\vert$ est fini, et comme $k$ est algébriquement clos $\vert k^\times\vert$ est divisible, d'où $\vert k^\times\vert=\vert\HHH(x)^\times\vert$ et $\vert\HHH(x)^\times\vert$ est divisible. Par ailleurs comme le corps $\HC(x)$ est \emph{stable} (cf \cite{Duc}, Théorème $4.3.14$) nous avons l'égalité : 
$$[\HC(y): \HC(x)]=[\widetilde{\HC(y)}: \widetilde{\HC(x)}]\cdot[\vert\HC(y)^\times\vert :\vert\HC(x)^\times\vert] $$
Ainsi $$[\HC(y): \HC(x)]=[\widetilde{\HC(y)}: \widetilde{\HC(x)}]=[\widetilde{k}(\mathscr{C}_y) : \widetilde{k}(\mathscr{C}_x)]=\deg(\mathscr{C}_y\rightarrow \mathscr{C}_x),$$ en particulier lorsque $\varphi$ est galoisien alors $\varphi$ est modéré en $y$ si et seulement si le morphisme $\mathscr{C}_y\rightarrow \mathscr{C}_x$ induit par $\varphi$ en $y$ est de degré premier à $p$.
\end{rem}\

En combinant la remarque précédente avec les propositions \ref{allure locale} et \ref{allureD} nous en déduisons naturellement un résultat dont nous nous servirons par la suite pour étudier le comportement local de $\GG(X,S)$ :

\begin{coro}[Allure locale de $\Covt(X)$]\label{allure} 
Soit $X$ une courbe $k$-analytique connexe quasi-lisse munie d'une triangulation généralisée ayant un unique sommet $s$, de type $2$, et dont le squelette ne contient aucune arête qui boucle sur $s$. Soit $\mathscr{U}_s$ le sous-schéma ouvert de $\mathscr{C}_s$ associé à $s$, tel que l'on ait une bijection entre $\br(X,s)$ et les points fermés de $\mathscr{U}_s$. Soit $\mathscr{D}_s \subset \mathscr{U}_s$ le sous-schéma ouvert correspondant aux éléments de $\DD(X,s)$. Alors la catégorie $\Covt(X)$ est équivalente à la catégorie des revêtements finis de $\mathscr{U}_s$ dominés par des revêtements galoisiens modérément ramifiés de degré premier à $p$ et étales au-dessus de $\mathscr{D}_s$.
\end{coro}\

\begin{defi}\label{tame group} 
Soit $\mathscr{C}$ une courbe propre, normale et intègre au-dessus de $\widetilde{k}$. Soit $\mathscr{D}\subset \mathscr{C}$ une sous-courbe ouverte (éventuellement égale à $\mathscr{C}$) et $\mathscr{Y}\rightarrow  \mathscr{D}$ un revêtement fini étale connexe de $\mathscr{D}$. On dit que $\mathscr{Y}$ est \emph{modérément ramifié le long de $\mathscr{C}\setminus \mathscr{D}$} lorsque les valuations définies par les points fermés de $\mathscr{C}\setminus \mathscr{D}$ sont toutes \emph{modérément ramifiées} dans l'extension $\widetilde{k}(\mathscr{C})/\widetilde{k}(\mathscr{D})$.\

 Notons $\Cov^\ta(\mathscr{D})$ la catégorie des revêtements finis étales de $\mathscr{D}$ modérément ramifiés le long de $\mathscr{C}\setminus \mathscr{D}$. C'est une catégorie galoisienne dont nous noterons $\pi_1^\ta(\mathscr{D})$ le groupe fondamental, défini à isomorphismes intérieurs près. Notons qu'il ne dépend que de $\mathscr{D}$, car $\mathscr{C}$ se retrouve à partir de $\mathscr{D}$ comme égal à sa compactification normale. 
 \end{defi}\
 
Rappelons que si $G$ est un groupe et $p$ un nombre premier, la \emph{$p'$-complétion profinie} de $G$, notée $G^{\pp}$, est le groupe profini défini comme la limite inverse du système projectif des quotients finis $G/N$ où $N$ décrit l'ensemble des sous-groupes distingués de $G$ d'indice premier à $p$. Lorsque $G$ est lui-même profini $G^{\pp}$ sera appelé le \emph{$p'$-quotient maximal} de $G$.\\

\begin{rem}\label{remarque sur le tame groupe} Reprenons les notations de la définition \ref{tame group}. Notons $\pi_1(\mathscr{D})$ le groupe fondamental étale de $\mathscr{D}$. En remarquant que toute extension de corps valués dominée par une extension galoisienne de degré premier à p est automatiquement modérément ramifiée on obtient l'isomorphisme : $$\pi_1(\mathscr{D})^{\ta, {\pp}}=\pi_1(\mathscr{D})^{\pp}$$
\end{rem}

 \begin{prop}[Grothendieck]\label{Grothendieck} 
 Soit $\mathscr{C}$ une courbe propre, normale et intègre de genre $g$ au-dessus d'un corps algébriquement clos de caractéristique $p > 0$. Soit $\mathscr{D}\subset \mathscr{C}$ une sous-courbe ouverte (éventuellement égale à $\mathscr{C}$) et $n\geqslant 0$ le nombre de points fermés de $\mathscr{C}\setminus \mathscr{D}$. Alors $\pi_1(\mathscr{D})^{\pp}$ est isomorphe à la $p'$-complétion profinie de $\Pi_{g,n}$ où $\Pi_{g,n}$ est le groupe fondamental (topologique) du complémentaire de $n$ points sur une surface de Riemann sur $\C$ de genre $g$, défini par la formule : 
 
$$\Pi_{g,n}:= \langle a_1, b_1, \ldots, a_g, b_g, \gamma_1, \ldots, \gamma_n \,\,\vert \,\, [a_1, b_1] \ldots [a_g, b_g] \gamma_1 \ldots \gamma_n=1\rangle$$
 
 \end{prop}\
 
 \begin{defi}
 Soit $X$ une courbe $k$-analytique quasi-lisse munie d'une triangulation généralisée $S$ de squelette $\Gamma$. Si $s\in S$, définissons la \emph{valence généralisée en $s$} comme l'entier $$n_s = \mathrm{val}_{\Gamma}(s) + \widetilde{n}_s \in \N,$$ où \;$\mathrm{val}_{\Gamma}(s)$ désigne la valence de $s$ dans le semi-graphe $\Gamma$, et $\widetilde{n}_s$ est le cardinal du complémentaire de l'image de $\br(X,s)\hookrightarrow \mathscr{C}_s(\widetilde{k}) $ (cf. \ref{courbe résiduelle}) dans $\mathscr{C}_s(\widetilde{k})$ : $$ \widetilde{n}_s= \sharp (\mathscr{C}_s(\widetilde{k}) \setminus \br(X,s)). $$
 \end{defi}\
 
 \begin{coro}[Allure locale de $\GG(X,S)$]\label{allure locale groupe} 
  Soit $X$ une courbe $k$-analytique quasi-lisse munie d'une triangulation généralisée $S$. Si $s\in S$ est un point de type $2$ de valence généralisée $n_s$, notons $\mathscr{C}_s$ la courbe résiduelle en $s$, $g(s)$ son genre (qui est par définition le genre du point $s$), et $\mathscr{D}_s$ le sous-schéma ouvert correspondant aux disques, i.e. aux éléments de $\DD(X,S,s)$. Alors $n_s$ n'est autre que le cardinal de $\mathscr{C}_s(\widetilde{k})\setminus \mathscr{D}_s(\widetilde{k})$, et le groupe fondamental de la composante sommitale en $s$ de $\GG(X,S)$ est : 
  $$\pi_s= \pi_1(\mathscr{D}_s)^{\pp} \simeq \Pi_{g(s),n_s}^{\pp}$$
 \end{coro}
 
 \begin{proof}
Puisque l'étoile $\mathrm{St}(X,S,s)$ est définie sur un revêtement topologique universel pointé $(\widetilde{X},\widetilde{s})$ de $(X,s)$, $\mathrm{St}(X,S,s)$ est une courbe $k$-analytique quasi-lisse dont $\lbrace \widetilde{s}\rbrace$ constitue une triangulation généralisée et dont le squelette ne contient aucune arête qui boucle sur $ \widetilde{s}$. Par ailleurs, on a $g(s)=g(\widetilde{s})$ et $n_s=n_{\widetilde{s}}$. Cela permet d'appliquer le corollaire \ref{allure} à l'étoile $\mathrm{St}(X,S,s)$, puis de conclure par la proposition \ref{Grothendieck} couplée à la remarque \ref{remarque sur le tame groupe}.
 \end{proof}

\subsection{Compatibilité avec les hypothèses mochizukiennes}\label{Compatibilité}

\subsubsection{Simplification des hypothèses}

Soit $X$ une courbe $k$-analytique quasi-lisse connexe munie d'une triangulation généralisée $S$ de squelette $\Gamma$. Nous voudrions savoir à quelles conditions sur $S$ le graphe d'anabélioïdes $\GG^\natural(X,S)$ construit plus haut vérifie l'ensemble des hypothèses que nous appellerons \emph{mochizukiennes} permettant d'appliquer le théorème de reconstruction \ref{reconsmochi}, à savoir : de type injectif, quasi-cohérent, totalement élevé, totalement détaché et sommitalement mince. D'après \ref{sous-section rabotée d'une courbe} on peut sans perte de généralité supposer que $S$ est une vraie triangulation de $X$. Supposons $\GG^\natural(X,S)$ connexe et non réduit à un sommet.\\

Vérifions tout d'abord que les hypothèses mochizukiennes sur $\GG^\natural(X,S)$ excluent l'appartenance à $S$ de \og points superflus \fg{} (cf. \ref{points superflus}). Supposons que $s\in S$ est tel que $g(s)=0$ et $s\not\in \partial^\an X$. Notons alors $\Gamma_0$ la composante connexe de $\Gamma\setminus \lbrace S\setminus\lbrace s\rbrace \rbrace $ contenant $s$, et distinguons les deux cas donnés par la proposition \ref{points superflus}: \\

 \begin{itemize}
 \item[•] \textit{Cas où $\Gamma_0$ est un intervalle ouvert} : \

 Dans ce cas $\St(X,S,s)$ est une couronne de squelette analytique $\Gamma_0$. Si $b$ (associée à une arête $e$) est une branche aboutissant à $s$, alors le morphisme $b_* : \pi_b\rightarrow \pi_s$ correspond au foncteur $\iota_b : \Covt(\mathrm{St}(X,S,s))\rightarrow \Covt(\CC_e)$ induit naturellement par l'e plongement naturel $\CC_e\hookrightarrow \mathrm{St}(X,S,s)$. Comme il s'agit d'un \emph{bon plongement} alors $\iota_b$ est une équivalence de catégorie, par conséquent $b_*$ correspond au morphisme identité $\widehat{\Z}^{\pp}\rightarrow \widehat{\Z}^{\pp}$, ce qui exclut la possibilité pour $\GG^\natural(X,S)$ d'être totalement détaché.\\
 
 \item[•] \textit{Cas où $\Gamma_0$ est un intervalle semi-ouvert d'extrémité $s$} : \
 
 Dans ce cas $\St(X,S,s)$ est un disque $k$-analytique (cf. \cite{Duc}, 5.1.16.2). Ainsi le morphisme $b_* : \pi_b\rightarrow \pi_s$ correspondant au foncteur $\iota_b : \Covt(\mathrm{St}(X,S,s))\rightarrow \Covt(\CC_e)$ induit naturellement par le plongement naturel $\CC_e\hookrightarrow \mathrm{St}(X,S,s)$ n'est autre que le morphisme trivial $\widehat{\Z}^{\pp}\rightarrow \lbrace 1\rbrace$, ce qui interdit à $\GG^\natural(X,S)$ d'être de type injectif.
  \end{itemize}\
  
 De là nous en déduisons qu'une triangulation généralisée $S$ sur $X$ vérifiant les hypothèses mochizukiennes sera automatiquement \og minimale \fg{}  (i.e. ne contenant aucun point superflu) et constituée des nœuds de son squelette. Par ailleurs, une triangulation généralisée $S$ de $X$ vérifiant les hypothèses mochizukiennes ne contient nul point de type $3$. En effet, supposons que $s\in S$ est un nœud de type $3$ de son squelette. Alors $s\in \partial^\an X$ et nous savons que $s$ admet un voisinage ouvert qui est une couronne $Z$ de type $]*,*]$. Ainsi $\St(X,S,s)$ est elle-même une couronne de type $]*,*]$. Si $b$ (associée à une arête $e$) est la branche aboutissant à $s$, le plongement $\CC_e\hookrightarrow \mathrm{St}(X,S,s)$ est un \emph{bon plongement}, ainsi le foncteur induit $\iota_b : \Covt(\mathrm{St}(X,S,s))\rightarrow \Covt(\CC_e)$ est une équivalence de catégorie, donc le morphisme $b_* : \pi_b\rightarrow \pi_s$ correspond à l'isomorphisme $\widehat{\Z}^{\pp}\rightarrow \widehat{\Z}^{\pp}$, ce qui exclut la possibilité pour $\GG^\natural(X,S)$ d'être totalement détaché. Nous allons donc seulement considérer dans la suite des triangulations généralisées constituées de points de type $2$ qui sont les nœuds du  squelette qu'elles constituent.\\

 Il reste à montrer que si une triangulation $S$ de squelette $\Gamma$ est exactement constituée des nœuds de $\Gamma$, alors le graphe d'anabélioïdes associé, $\GG^\natural(X,S)$, vérifie sous certaines hypothèses à préciser l'ensemble des conditions mochizukiennes. Eu égard à l'ingratitude des définitions de ces hypothèses la tâche semble malaisée, mais elle se simplifie en combinant la proposition \ref{allure} à la remarque suivante de Mochizuki : \\
 
 \begin{prop}[\cite{M3}, Example 2.10] \label{remarque de Mochizuki} Soit $\Sigma$ un ensemble non vide de nombres premiers contenant au moins un entier différent de $p$. Soit $\GG$ un semi-graphe d'anabélioïdes tel que, pour chaque sommet $v$, $\pi_v$ est le pro-$\Sigma$ quotient maximal du groupe fondamental d'une surface de Riemann hyperbolique, et tel que les différents morphismes $\pi_b\hookrightarrow \pi_v$ correspondent aux inclusions des groupes d'inertie des «cusps» de la courbe hyperbolique. Alors $\GG$ est cohérent, totalement élevé, totalement détaché et sommitalement mince. 
 \end{prop}\
 
 \begin{rem}
 Nous n'allons pas prouver ici la proposition \ref{remarque de Mochizuki}, mais observer que la minceur sommitale est une conséquence du corollaire $1.3.4$ de \cite{Nak}. En effet si $n>0$ le groupe $\Pi_{g,n}^{\pp}$  associé à un point $s\in S$ est un pro-$p'$ groupe libre de rang $2g+n-1$. Or Nakamura nous dit exactement que les $p'$-complétions des groupes libres de rang $\geqslant 2$ ou des groupes du type $\Pi_{g,0}$ avec $g\geqslant 2$ ont un centre trivial. Par ailleurs si $H$ est un sous-groupe ouvert de $\Pi_{g,n}^{\pp}$ nous savons que $H$ est le groupe $\pi_1 ^{\ta, \pp}(\mathscr{H}_s)$ correspondant à un revêtement étale de degré premier à $p$ (et donc modérément ramifié le long de $\mathscr{C}_s\setminus \mathscr{U}_s$) de $\mathscr{U}_s$. Or la formule de Riemann-Hurwitz nous assure que le genre de $\mathscr{H}_s$ ne peut pas être strictement plus petit que celui de $\mathscr{U}_s$, la courbe $\mathscr{H}_s$ reste donc hyperbolique et ainsi $H$ est du type $\Pi_{g',n'}^{\pp}$ avec $2g'+n'>2$. Par conséquent $H$ est également de centre trivial, et cela est vrai pour tout sous-groupe $H$ ouvert dans $\Pi_{g,n}^{\pp}=\pi_s$, ce qui signifie exactement que $\pi_s$ est mince. 
 \end{rem}\
 
  En reprenant les notations de la proposition \ref{Grothendieck}, rappelons que la courbe $\mathscr{D}$ est dite \emph{hyperbolique} lorsque : $2g+n>2$. De la même manière une \emph{surface de Riemann hyperbolique} est définie comme le complémentaire de $n\geqslant 0$ points d'une surface de Riemann compacte de genre $g$ telle que $2g+n>2$, i.e. lorsqu'on est dans l'un des cas suivants : 
 \begin{itemize}
 \item[$\mathit{i)}$]$g=0$ et $n\geqslant 3$
 \item[$\mathit{ii)}$]$g=1$ et $n\geqslant 1$
 \item[$\mathit{iii)}$]$g\geqslant 2$
 \end{itemize}\

 Cela nous mène naturellement aux définitions suivantes : 
 
 \begin{defi}[Nœuds hyperboliques]\label{noeuds hyperboliques}
 Soit $X$ une courbe $k$-analytique connexe quasi-lisse munie d'une triangulation généralisée $S$ uniquement constituée des nœuds de son squelette $\Gamma$. Soit $s\in S$, rappelons que la \emph{valence généralisée en $s$} est définie comme l'entier $$n_s = \mathrm{val}_{\Gamma}(s) + \widetilde{n}_s \in \N$$ où \;$\mathrm{val}_{\Gamma}(s)$ désigne la valence de $s$ dans le semi-graphe $\Gamma$, et $\widetilde{n}_s$ est le cardinal du complémentaire de l'image de $\br(X,s)\hookrightarrow \mathscr{C}_s(\widetilde{k}) $ (cf. \ref{courbe résiduelle}) dans $\mathscr{C}_s(\widetilde{k})$ : $$ \widetilde{n}_s= \sharp (\mathscr{C}_s(\widetilde{k}) \setminus \br(X,s)). $$
 Le nœud $s$ est dit \emph{nœud hyperbolique} lorsque : 
 $\;\;\;2g(s)+n_s>2. $

 \end{defi}\
 
 \begin{rem}
 Soit $\mathscr{U}_s$ l'ouvert de Zariski de $\mathscr{C}_s$ tel que $\br(X,s)\hookrightarrow \mathscr{C}_s(\widetilde{k}) $ établisse une bijection entre $\br(X,s)$ et $\mathscr{U}_s(\widetilde{k})$, et $\mathscr{D}_s\subset \mathscr{U}_s$ le sous-ouvert de Zariski tel que les branches discales (i.e. associées aux éléments de $\DD(X,S,s)$) correspondent exactement à $\mathscr{D}_s(\widetilde{k})$. Alors on a  : \begin{align*}
&\mathrm{val}_{\Gamma}(s)=\sharp(\mathscr{U}_s(\widetilde{k})\setminus \mathscr{D}_s(\widetilde{k})) \\
& \widetilde{n}_s=  \sharp(\mathscr{C}_s(\widetilde{k})\setminus \mathscr{U}_s(\widetilde{k})) \\
& n_s= \sharp(\mathscr{C}_s(\widetilde{k})\setminus \mathscr{D}_s(\widetilde{k}))
\end{align*}
Par ailleurs, $\widetilde{n}_s>0$ si et seulement si $s\in \partial^\an X$.

 \end{rem}\
 
 \begin{rem}\label{exemple de noeud non hyperbolique} 
 Avec les définitions précédentes, un nœud $s$ de $X$ non isolé dans $\Gamma$ (i.e. $\mathrm{val}_{\Gamma}(s)\neq 0$) est non hyperbolique seulement dans la situation où $s$ est unibranche dans $\Gamma$, de genre $0$ et sur le bord analytique de $X$ mais qu'il ne « manque qu'une branche », c'est-à-dire $\widetilde{n}_s=1$. C'est le cas par exemple d'une couronne compacte : son bord analytique est constitué de deux nœuds qui forment une triangulation mais qui sont deux nœuds non hyperboliques. Par ailleurs un nœud de type $3$ ne peut jamais être hyperbolique. En effet, si $s$ est un tel nœud, nous savons que $s$ admet un voisinage ouvert qui est une couronne $Z$ de type $]*,*]$ avec $\lbrace s\rbrace=\partial^\an Z$, ce qui donne $n_s=2$, et ainsi $2g(s)+n_s=2$, d'où l'obstruction à l'hyperbolicité.   
 \end{rem}\
 
 \subsubsection{Courbes $k$-analytiquement hyperboliques}
 
  \begin{defi}[Courbes $k$-analytiquement hyperboliques]\label{définition des courbes marquées anabéliennes} \
Soit $X$ une courbe $k$-analytique connexe quasi-lisse. Nous dirons que $X$ est une courbe \emph{$k$-analytiquement hyperbolique} lorsque l'ensemble des nœuds de son squelette analytique $S^\an(X)$ est non vide et uniquement constitué de nœuds hyperboliques. 
 \end{defi}\
 
 \begin{prop}\label{corollaire indiquant que le granphe vérifie les hypothèses Mochizukiennes} 
 Soit $X$ une courbe $k$-analytiquement hyperbolique dont on note $\Sigma_X \subset S^\an(X)$ l'ensemble non vide des nœuds du squelette analytique (automatiquement de squelette $S^\an(X)$). Alors le graphe d'anabélioïdes $\GG^\natural(X,\Sigma_X)$ vérifie l'ensemble des hypothèses \emph{mochizukiennes}, à savoir : de type injectif, quasi-cohérent, totalement élevé, totalement détaché et sommitalement mince.\
 \end{prop}
 
\begin{proof}

D'après la proposition \ref{remarque de Mochizuki}, afin de montrer ce résultat il suffit de montrer que notre graphe d'anabélioïdes $\GG^\natural(X,\Sigma_X)$ vérifie les hypothèses prescrites dans la proposition sus-citée, en prenant $\Sigma=\P\setminus \lbrace p\rbrace$ (auquel cas le \emph{pro-$\Sigma$ quotient maximal} ne sera rien d'autre que le \emph{$p'$-quotient maximal}) :\\
 
 \begin{itemize}

 \item[•]\textit{Pour tout $s\in \Sigma_X$, $\pi_s$ est le $p'$-quotient maximal du groupe fondamental d'une surface de Riemann hyperbolique :}\\
 
 D'après le corollaire \ref{allure locale groupe} nous savons que $\pi_s$ est le $p'$-quotient maximal de $\Pi_{g(s),n_s}$. Puisque $s$ est par hypothèse un nœud hyperbolique nous savons que ${2g(s)+n_s>2}$, c'est-à-dire que $\Pi_{g(s),n_s}$ est le groupe fondamental d'une surface de Riemann hyperbolique (de genre $g(s)$ avec $n_s$ points manquants).
 
 \

  \item[•] \textit{Les morphismes $\pi_b\xrightarrow{b_*} \pi_s$ correspondent aux inclusions des groupes d'inertie des «cusps» :}\\
  
  Reprenons ici les notations $\OO_X(b)_\sep$ et $\mathcal{I}_b^\ttt$ de la partie \ref{courbe résiduelle}. Soit $\mathcal{A}$ une $\OO_X(b)_\sep$-algèbre finie modérément ramifiée, i.e. $\mathcal{A}$ définit un $\mathcal{I}_b^\ttt$-ensemble fini, autrement dit un objet de $\BB(\mathcal{I}_b^\ttt)$. Via la tensorisation par $\OO_X(b)$ au-dessus de $\OO_X(b)_\sep$, $\mathcal{A}$ définit un revêtement étale fini de $\CC_e$ qui est modéré. Réciproquement, tout revêtement modéré $Y$ de $\CC_e$ provient d'un revêtement modéré de $X$ tout entier. Cela découle du caractère kummérien des revêtements modérés des couronnes (\ref{proposition caractere kummerien et trivial}) couplé à la proposition \ref{diagramme de suites exactes} pour un choix idoine de $\ell\in \N$ premier à $p$. Un tel revêtement fourni une $\OO_X(b)_\sep$-algèbre modérément ramifiée : cela découle par exemple de la Proposition $4.2.20$ de \cite{Duc} selon laquelle une branche $a\in \br(Y,y)$  au-dessus de $b$  induit une extension de corps $\OO_Y(a)_\sep/\OO_X(b)_\sep$ finie de degré $\deg (a\rightarrow b)$.\
  
 Ceci définit une équivalence de catégories $\BB(\mathcal{I}_b^\ttt)\xrightarrow{\tau^*} \Cov^\ttt(\CC_e)$ qui induit un  isomorphisme de groupes : $\pi_b\xrightarrow{\tau} \mathcal{I}_b^\ttt  $ (nous savions déjà d'après les parties \ref{partie comportement des composantes coronaires} et \ref{courbe résiduelle} que ces groupes sont tous deux isomorphes à $\widehat{\Z}^{\pp}$). Comme $\mathcal{I}_b^\ttt=\mathcal{I}_{\kappa(s)_b}^\ttt$, $\tau$ correspond à un isomorphisme $\pi_b \iso \mathcal{I}_{\kappa(s)_b}^\ttt$.\
  
  Le corps résiduel de $\widetilde{\kappa(s)}=\widetilde{\HC(s)}=\widetilde{k}(\mathscr{C}_s)$ muni de $\langle .\rangle_b$ n'est autre que $\widetilde{k}$ qui est algébriquement clos. Par conséquent, en reprenant les raisonnements précédents on obtient l'égalité : $$\mathcal{I}_{\widetilde{\kappa(s)}_b}^\ttt= \widehat{\Z}^{\pp},$$ 
  
où $\mathcal{I}_{\widetilde{\kappa(s)}_b}^\ttt$ désigne le groupe d'inertie modérée de $\widetilde{\kappa(s)}$ munie de $\langle . \rangle_b$ (où l'on a préalablement fixé une clôture séparable de $\widetilde{\kappa(s)}$ munie d'une valuation prolongeant $\langle . \rangle_b$ sur $\widetilde{\kappa(s)}$). On en déduit que le morphisme naturel $\mathcal{I}_{\kappa(s)_b}^\ttt \xrightarrow{\eta} \mathcal{I}_{\widetilde{\kappa(s)}_b}^\ttt$ est un isomorphisme. Notons $\mathscr{P}$ le point de $\mathscr{C}_s$ correspondant à $b$. Le groupe $\mathcal{I}_{\widetilde{\kappa(s)}_b}^\ttt$ correspond alors au groupe d'inertie modérée $\mathcal{I}_\mathscr{P}^\ttt$ du point $\mathscr{P}$, qui se plonge dans $\pi_s$. Par ailleurs on a le diagramme commutatif suivant : 

$$\begin{tikzcd}[row sep=2.5em]
 & \mathcal{I}_{\kappa(s)_b}^\ttt \arrow{dr}{\eta} \\
\pi_b \arrow{ur}{\tau} \arrow{rr}{b_*} && \mathcal{I}_\mathscr{P}^\ttt\subset \pi_s
\end{tikzcd}$$

où $\tau$ et $\eta$ sont les isomorphismes naturels invoqués plus haut. En conséquence $b_*$ correspond bien au plongement du groupe d'inertie du \emph{cusp} $\mathscr{P}$ dans $\pi_s$.

 \end{itemize}\
 
 \end{proof}
 
 \begin{rem}
 Soit $S=\{s\}$ une triangulation réduite à un singleton d'une courbe $k$-analytique $X$ quasi-lisse connexe dont le squelette n'admet aucune arête qui boucle sur $s$. Dans ce cas $\GG^\natural(X,S)$ n'est rien d'autre que l'anabélioïde connexe $\BB(\pi_s)$. Comme le graphe sous-jacent n'est constitué que d'un sommet et qu'il n'y a aucune branche, alors $\GG^\natural(X,S)$ est automatiquement de type injectif, totalement élevé, totalement détaché et quasi-cohérent. D'après \cite{Nak}, $1.3.4$, le groupe $\pi_s$ n'est mince que si le comportement en $s$ est hyperbolique, ce qui correspond à $2g(s)+n_s>2$. On retrouve ainsi un cas particulier de la proposition \ref{remarque de Mochizuki}.
 \end{rem}\
 
 \begin{enumerate}

 \item\textbf{Exemples de courbes $k$-analytiquement hyperboliques compactes} \\
 
 Soit $X$ une courbe $k$-analytique connexe quasi-lisse, \emph{compacte} et non vide. Notons $\Sigma_X$ l'ensemble des nœuds de $S^\an(X)$. Nous savons que $S^\an(X)=\emptyset$ si et seulement si $X= \P_k^{1,\an}$ ($k$ est supposé algébriquement clos), dans ce cas \textit{a fortiori} $\Sigma_X=\emptyset$. Si $S^\an(X)\neq\emptyset$, alors $\Sigma_X=\emptyset$ si et seulement si $X$ est une courbe de Tate (cf. \ref{triangulation minimale dans le cas compact}), auquel cas $S^\an(X)$ est un cercle. Dans ces deux cas $X$ n'est à l'évidence pas $k$-analytiquement hyperbolique. Exclure ces deux cas revient à supposer que $\Sigma_X\neq \emptyset$. Alors deux cas se présentent : \bigskip
 \begin{enumerate}
 \item si $S^\an(X)=\lbrace s\rbrace$ est un singleton, alors $X$ est $k$-analytiquement hyperbolique si et seulement si $s$ est un nœud hyperbolique, ce qui est le cas dès que $s$ vérifie $g(s)\geqslant 2$ ou $\lbrace g(s)= 1$ et  $s\in \partial^\an X\rbrace$. Un disque n'est jamais $k$-analytiquement anabélien;
 \item si $S^\an(X)$ contient au moins une arête et $X$ n'est pas une courbe de Tate, alors $X$ est $k$-analytiquement hyperbolique dès que ses points unibranches sont de genre $>0$, ceci permet d'éviter les nœuds non hyperboliques décrits au \ref{exemple de noeud non hyperbolique}.
 \end{enumerate}\ 
 
 \item\textbf{Courbes marquées comme exemples de courbes $k$-analytiquement hyperboliques non compactes} \\

 \begin{defi}
 Une \emph{courbe $k$-analytique marquée} est définie comme une courbe $k$-analytique $Y$ de la forme $Y=X\setminus \E$ où $(X,\E)$ est la donnée d'une courbe $k$-analytique non vide $X$, quasi-lisse et connexe, et d'un sous-ensemble $\E\subset X$ fermé et discret de points \emph{rigides} de $X$. 
 
\end{defi}\

Soit $Y=X\setminus \E$ une courbe marquée. Comme les éléments de $\E$ en tant que points rigides sont des points unibranches de $X$ par quasi-lissité, la courbe $Y$ reste connexe, quasi-lisse, mais ne jouit pas de la compacité dès que $\E\neq \emptyset$, même lorsque $X$ est elle-même compacte. Nous allons étudier les nœuds de $S^\an(Y)$, et notamment les conditions d'hyperbolicité de ces derniers. \\

\textbf{Si $S^\an (X)$ est non vide :} Nous savons que c'est un sous-graphe analytiquement admissible de $X$, et l'on note $r$ la rétraction canonique de $X$ sur $S^\an (X)$. Remarquons que l'on a : $$S^\an(Y) = S^\an (X)\cup\bigcup_{x\in \E}]x, r(x)].$$ 

Notons $\Sigma_Y$ l'ensemble des nœuds de $S^\an(Y)$, et $\Sigma_{Y}^\bot$, l'ensemble des sommets pluribranches (de valence $\geqslant 3$) de $S^\an(Y)\setminus S^\an(X)$. De tels points apparaissent lorsque deux points de $\E$ sont sur le même disque de $X\setminus \Sigma_X$ : si $t\in \Sigma_{Y}^\bot$, il existe $x_1, x_2\in \E$ tels que $r(x_1)=r(x_2):=r_0$, avec $$]x_1,r(x_1)]\cap ]x_2,r(x_2)]=[t, r_0].$$ Tous les éléments de $\Sigma_{Y}^\bot$ sont des nœuds hyperboliques de $S^\an(Y)$ puisqu'ils sont pluribranches.

Par ailleurs $\Sigma_Y$ peut être décrit par :  $$\Sigma_Y=\Sigma_X \cup\bigcup_{x\in \E}\lbrace r(x)\rbrace\cup \Sigma_{Y}^\bot. $$

Soit $x\in \E$, alors $x\notin S^\an(X)$ puisque ce dernier n'est constitué que de points de type $2$ ou $3$ et que $x$ est rigide. Alors $]x, r(x)]$ définit une branche de $r(x)$ dans $S^\an(Y)$ qui n'est pas dans $S^\an(X)$. Par conséquent :  $${\mathrm{val}_{S^\an(Y)}(r(x))>\mathrm{val}_{S^\an(X)}(r(x))}.$$ 
Supposons que $S^\an(X)$ n'est pas réduit à un singleton. Si $r(x)$ est de valence $2$ en tant que point de $S^\an(X)$, cela implique que ${\mathrm{val}_{S^\an(Y)}(r(x))\geqslant 3}$, ce qui fait de $r(x)$ un nœud hyperbolique de $S^\an(Y)$. \

Si maintenant $r(x)$ est un point unibranche de $S^\an(X)$, nous savons que c'est un nœud de $S^\an(X)$ qui est de genre $\geqslant 1$ ou sur le bord analytique de $X$. Ainsi $r(x)$ reste un nœud de $S^\an(Y)$. On a au pire des cas (lorsque le nœud est non hyperbolique): $2g(r(x))+\mathrm{val}_{S^\an(X)}(r(x))+\widetilde{n_{r(x)}}\geqslant 2$. Par conséquent : $$2g(r(x))+\mathrm{val}_{S^\an(Y)}(r(x))+\widetilde{n_{r(x)}}> 2,$$ ce qui fait de $r(x)$ un nœud hyperbolique de $S^\an(Y)$.\bigskip

Si on suppose que $S^\an(X)$ est réduit à un singleton $\lbrace s\rbrace$, alors tous les points de $\E$ se rétractent sur $s$, c'est-à-dire que pour tout $x\in \E$, $r(x)=s$. Si $s$ est un nœud hyperbolique dans $S^\an(X)$, il le reste dans $S^\an(Y)$. Si $g(s)\geqslant 1$ et $\E\neq \emptyset$, alors $s$ est un nœud hyperbolique de $S^\an(Y)$ par l'inégalité ${\mathrm{val}_{S^\an(Y)}(s)>\mathrm{val}_{S^\an(X)}(s)}$. En revanche, lorsque $g(s)=0$, le nœud $s$ de $S^\an(Y)$ n'a \textit{a priori} plus aucune raison d'être hyperbolique. Il le sera dès que $\widetilde{n_s}\geqslant 2$, ou dès que $s\in \partial^\an(X)$, que $\mathrm{Card}(\E)\geqslant 2$ et que les points de $\E$ sont « suffisamment éloignés » au sens où il existe $x_1, x_2\in \E$ tels que $]x_1,s]\cap ]x_2,s]=\lbrace s\rbrace$ (i.e. $x_1$ et $x_2$ ne sont pas dans le même disque de $X\setminus \lbrace s\rbrace$). 

\bigskip

\textbf{Si $S^\an (X)$ est vide :}  Lorsque $X$ est compact cela correspond à $X=\P_k^{1,\an}$. Dans tous les cas $S^\an(Y)$ est l'enveloppe convexe de $\E$, et $\Sigma_Y$ correspond exactement aux sommets pluribranches de cette enveloppe convexe, c'est-à-dire que $\Sigma_Y$ est non vide si et seulement si $\mathrm{Card}(\E)\geqslant 3$, auquel cas tous les éléments de $\Sigma_Y$ sont des \emph{nœuds hyperboliques}.\\

\begin{rem}
Nous savons que $\Sigma_Y$ est une triangulation de $X$ tout autant que de $Y$, avec pour cette dernière des composantes non relativement compactes qui apparaissent. Soit $e\in \E$ et $V$ la composante connexe de $X\setminus \Sigma_Y$ contenant $e$. Nous savons que $V$ est un disque $k$-analytique ouvert (\cite{Duc}, $5.5.3$). Si l'on se place maintenant dans $Y$, la composante connexe de $Y\setminus \Sigma_Y$ correspondante est exactement $V\setminus \lbrace e \rbrace$. Ce n'est ainsi plus un disque, mais une couronne $k$-analytique $\mathcal{C}_e$ non relativement compacte dans $Y$, c'est-à-dire un élément de $\mathcal{C}^\infty(Y,\Sigma_Y)$. Lorsque $X$ est compacte il est facile de voir que l'application $e\mapsto \mathcal{C}_e$ établit une bijection entre $\E$ et $\mathcal{C}^\infty(Y,\Sigma_Y)$.\\
\end{rem}

\textbf{Conclusion :} Une courbe marquée $Y=X\setminus \E$ avec $X$ compacte sera alors $k$-analytiquement hyperbolique dès que l'on est dans l'un des cas suivants :   
 \begin{enumerate}
 \item $X$ est elle-même une courbe $k$-analytiquement hyperbolique;
 \item $X$ est une courbe de Tate et $\E\neq \emptyset$;
 \item $X=\P_k^{1,\an}$ et $\mathrm{Card}(\E)\geqslant 3$;
 \item $S^\an(X)$ est réduit à un nœud \emph{non hyperbolique} $\lbrace s\rbrace$, et l'une des conditions suivantes est vérifiée : \begin{itemize}
 
 \item[•]  $g(s)=1$ et $\E\neq \emptyset$,
 \item[•] $s\in \partial^\an(X)$, $\mathrm{Card}(\E)\geqslant 2$, $\exists \;x_1, x_2\in \E, \;]x_1,s]\cap ]x_2,s]=\lbrace s\rbrace$.
 \end{itemize} 
 \end{enumerate}\

\begin{ex}\label{exemple des courbes elliptiques épointées comme courbes hyperboliques} 
 \emph{Le cas d'une courbe elliptique.} Soit $\mathscr{X}$ une courbe elliptique sur $k$ munie de son élément neutre $e\in \mathscr{X}(k)$. Si $\mathscr{X}$ est à mauvaise réduction, $S^\an (\mathscr{X}^\an)$ est un cercle ne contenant pas le point rigide $e$. Si $\mathscr{X}$ est à bonne réduction, son squelette est constitué d'un unique nœud de genre $1$. Dans les deux cas, $\mathscr{Y}=\mathscr{X}^\an\setminus \lbrace e\rbrace$ définit une courbe $k$-analytiquement hyperbolique. Cela permet d'associer de manière canonique une courbe marquée $k$-analytiquement hyperbolique à toute courbe elliptique sur $k$.
 \end{ex}\

 \item\textbf{Le demi-plan de Drinfeld}\\
 
Voici un exemple important de courbe $k$-analytiquement hyperbolique n'ayant aucune nature algébrique. Soit $k_0$ un corps non-archimédien \emph{local} de caractéristique $p>0$ (i.e. $k_0$ est soit une extension finie de $\Q_p$, soit du type $\F_{q}((T))$ où $q$ est une puissance de $p$) contenu dans le corps $k$ non archimédien, complet et algébriquement clos dont la valuation prolonge celle de $k_0$.\
 
 Le \emph{demi-plan de Drinfeld}, noté $\mathfrak{D}^1_{k/k_0}$, est la courbe $k$-analytique analogue non archimédien du demi-plan de Poincaré, définie en tant qu'ouvert de $ \P_{k}^{1, \an}$ par :
 
 $$\mathfrak{D}^1_{k/k_0}= \P_{k}^{1, \an}\setminus  \P_{k}^{1, \an}(k_0).  $$
 
 Le demi-plan de Drinfeld $\mathfrak{D}^1_{k/k_0}$ est une courbe $k$-analytiquement anabélienne non compacte et sans bord. En effet, son squelette analytique ainsi que l'ensemble de ses nœuds peuvent être décrits par : $$\left\{
\begin{array}{l}
 S^\an(\mathfrak{D}^1_{k/k_0})=\lbrace \text{enveloppe convexe de } \P_{k}^{1, \an}(k_0) \rbrace \setminus \P_{k}^{1, \an}(k_0)=\lbrace \eta_{a,r} \rbrace_{\;a\in k_0,\; r\in \R_+^*}  \\
  \\
  \Sigma_{\mathfrak{D}^1_{k/k_0}}=\lbrace \eta_{a,r} \rbrace_{\;a\in k_0,\; r\in \vert k_0^* \vert}
\end{array}
\right.$$

Ainsi, dans $S^\an(\mathfrak{D}^1_{k/k_0})$, tous les nœuds sont de valence $q+1$ (donc hyperboliques) où $q$ désigne le cardinal du corps fini $\widetilde{k_0}$. Par ailleurs $S^\an(\mathfrak{D}^1_{k/k_0})$ n'est pas seulement un semi-graphe, mais un vrai graphe (toutes les arêtes sont \og fermées \fg{}) infini, de là l'égalité avec le squelette tronqué : $$S^\an(\mathfrak{D}^1_{k/k_0})=S^\an(\mathfrak{D}^1_{k/k_0})^\natural.$$ Berkovich  montre dans \cite{Ber4} que le groupe $\Aut_k(\mathfrak{D}^1_{k/k_0})$ des automorphismes $k$-analytiques de $\mathfrak{D}^1_{k/k_0}$ coïncide avec $\mathrm{PGL}_2(k_0)$, et Alon étend dans \cite{Al} ce résultat à un produit fini  : si $k_1, \ldots, k_r\subset k$ sont des sous-extensions de $k$ qui sont des corps non-archimédiens locaux pour la valuation induite par celle de $k$, et si $M=\lbrace \sigma\in \mathfrak{S}_r, k_{\sigma(i)}=k_i, 1\leqslant i \leqslant r \rbrace$, alors : $$\Aut_k\left( \prod_{i=1}^r \mathfrak{D}^1_{k/k_i}\right)=M\ltimes \prod_{i=1}^r \mathrm{PGL}_2(k_i).$$
\begin{rem}
Le squelette $S^\an(\mathfrak{D}^1_{k/k_0})$ n'est autre que l'\emph{immeuble de Bruhat-Tits} de $\mathrm{SL}_2/k_0$, noté $\mathfrak{B}^1_{k_0}$. On s'en rend facilement compte en considérant $\mathfrak{B}^1_{k_0}$ selon sa description en termes de normes : considérons le $k_0$-espace vectoriel $$ V=\bigoplus_{i=0}^1 k_0T_i \subset k[T_0, T_1],$$ alors $\mathfrak{B}^1_{k_0}$ est défini comme l'ensemble des normes sur $V$ modulo la relation qui identifie deux normes $\rho_1, \rho_2$ si et seulement si celles-ci sont \emph{homothétiques}, c'est-à-dire s'il existe $c\in \R_+^*$ tel que $\rho_1=c.\rho_2$.
\end{rem}

 \end{enumerate}

\subsection{Reconstruction du squelette}\label{section reconstruction du squelette}

\subsubsection{Tempéroïde connexe associé à $\mathbf{\GG(X,S)}$}\label{subsection tempéroide connexe équivalence} 

\begin{prop}\label{equimod} 
Soit $X$ une courbe $k$-analytique quasi-lisse connexe munie d'une triangulation généralisée $S$ non vide. On a une équivalence de catégories (donc un isomorphisme d'anabélioïdes connexes) :  $$\Covt(X)\iso \BB(\GG(X,S))$$
\end{prop}

\begin{proof}
Remarquons tout d'abord que $\mathcal{U(X,S)}  := (\mathrm{St}(X,S,s_i))_{s_i\in S}$ est un recouvrement de $X$ pour la topologie étale sur le site $X_{\et}$, puisqu'un revêtement (universel) topologique est étale. Notons $U_{i,j}:=\mathrm{St}(X,S,s_i)\times_X \mathrm{St}(X,S,s_j)$.\\

Remarquons que si $(s_i, s_j)$ est un couple d'éléments distincts de $S$, alors : 
\begin{itemize}
\item[•]$\ds{U_{i,j}=\mathrm{St}(X,S,s_i)\times_X \mathrm{St}(X,S,s_j)\simeq\bigsqcup_e \CC_e}$, où $e$ décrit l'ensemble des arêtes de $\G(X,S)$ ayant deux branches aboutissant en $s_i$ et en $s_j$.
\item[•]$\ds{U_{i,i}=\mathrm{St}(X,S,s_i)}\times_X\mathrm{St}(X,S,s_i)\simeq \mathrm{St}(X,S,s_i)\bigsqcup_e \CC_e$, où $e$ décrit l'ensemble des arêtes de $\G(X,S)$ ayant deux branches aboutissant toutes deux en $s_i$.
\end{itemize}\

Une \emph{donnée de descente sur} $\mathcal{U(X,S)}$ est la donnée pour tout $s_i \in S$ d'un $Y_{s_i}\in \Covt(\mathrm{St}(X,S,s_i))$ et d'une famille d'isomorphismes 
$(\varphi_{j,i})_{j, \;s_j\in S}$ avec 
\begin{equation} \label{eq1} 
\varphi_{i,j} : Y_{s_i}\times_{\mathrm{St}(X,S,s_i)} U_{i,j}\iso Y_{s_j}\times_{\mathrm{St}(X,S,s_j)} U_{i,j},
\end{equation}
telle que chaque $\varphi_{i,i}$ se restreigne à l'identité sur la diagonale de $U_{i,i}$. Cette dernière contrainte provient de la conditions sur les cocycles $\varphi_{i,i}=\varphi_{i,i}\circ \varphi_{i,i}$. Les autres conditions de cocycles sont ici tautologiques puisqu'avec des notations évidentes $U_{i,j,k}=\emptyset$ dès que les indices $i, j$ et $k$ ne sont pas tous égaux.

 \bigskip 
Il reste à remarquer que les conditions de l'équation (\ref{eq1}) et la contrainte sur $\varphi_{i,i}$ correspondent exactement aux \og conditions de recollement \fg{}  permettant de définir les objets de la catégorie $\BB(\GG(X,S))$, d'où l'on déduit une équivalence de catégories entre $\BB(\GG(X,S))$ et la catégorie des données de descente sur $\mathcal{U}(X,S)$. D'autre part, par \emph{descente étale}, la catégorie $\Covt(X)$ est elle-aussi équivalente à la catégorie des données de descente sur $\mathcal{U}(X,S)$, ce qui donne par transitivité le résultat escompté. 

\end{proof}

\begin{coro}\label{equitemp}
Sous les mêmes hypothèses on a une équivalence de catégories $$\Cov^{\tp, \pp}(X)\iso \BB^\tp(\GG(X,S)),$$ qui induit un isomorphisme de tempéroïdes connexes : $$\left(\Cov^{\tp, \pp}(X),\Upsilon_X^{\pp}\right)\iso \left(\BB^\tp(\GG(X,S)),\Upsilon_{\GG(X,S)}^{\mathrm{g\acute{e}o}}\right).$$
\end{coro}

\begin{proof}
D'après la proposition précédente on peut identifier les anabélioïdes connexes $\Covt(X)$ et  $\BB(\GG(X,S))$. Sous cette identification l'ensemble des systèmes projectifs cofinaux de revêtements finis étales galoisiens de $\GG(X,S)$ correspond à l'ensemble des systèmes projectifs cofinaux de revêtements galoisiens finis de $X$ de degré premier à $p$ (donc modérés). Par ailleurs, si $\GG_i$ est un revêtement fini étale galoisien de $\GG(X,S)$ correspondant au revêtement $Y_i$ de $X$ selon l'équivalence de la proposition précédente, on a par définition du semi-graphe $\G_i$ sous-jacent à $\GG_i$ un isomorphisme de semi-graphes entre $S^\an(Y_i)$ et $\G_i$. En reprenant les notations de la remarque \ref{description du pi1 temp d'un graphe d'anabélioïdes en terme de limite} et de la proposition \ref{description du pi1 tp de l'espace analytique en terme de limite}, on a un isomorphisme entre les groupes $T_i:=\Gal(\widetilde{\GG}_i/\GG_i)$ et $\pi_1^\mathrm{top}(Y_i)$ (qui sont des groupes libres discrets) ainsi qu'entre les groupes  $\Gal(\GG_i/\GG(X,S))$ et $\Gal(Y_i/X)$. De plus ces isomorphismes sont fonctoriels en $Y_i$. De là on obtient un isomorphisme : $$\Gal(\widetilde{\GG}_i/\GG(X,S)) \simeq \Delta_i,  $$
qui induit en vertu de la remarque \ref{description du pi1 temp d'un graphe d'anabélioïdes en terme de limite} et de la proposition \ref{description du pi1 tp de l'espace analytique en terme de limite} un isomorphisme : $$ \pi_1^\tp(\GG(X,S))\simeq \pi_1^{\tp, \pp}(X) .$$
\end{proof}

\begin{rem}\label{remarque sur l'équivalence des catégories liées aux graphe tronqué et au semi-graphe} 
Si $X$ est une courbe $k$-analytique quasi-lisse et connexe munie d'une triangulation généralisée $S$, alors on a une équivalence de catégories suivante : $\BB^\tp(\GG(X,S))\simeq\BB^\tp(\GG^\natural(X,S))$. En effet, selon le même raisonnement que ce qui précède on montre que $\BB^\tp(\GG^\natural(X,S))$ est elle aussi équivalente à $\Cov^{\tp, \pp}(X)$. 
\end{rem}

\subsubsection{Reconstruction du squelette tronqué}

\begin{thm}\label{théorème de reconstruction du squelette tronqué} 
Soient $X_1$ et $X_2$ deux courbes $k$-analytiques quasi-lisses et connexes munies respectivement de triangulations généralisées $S_1$ et $S_2$ non vides de squelettes tronqués $\Gamma_1^\natural$ et $\Gamma_2^\natural$ dont les graphes d'anabélioïdes vérifient tous deux l'ensemble des hypothèses mochizukiennes. Supposons qu'il existe un isomorphisme entre les groupes tempérés modérés : $$\varphi : \pi_1^{\tp, \pp}(X_1)\iso \pi_2^{\tp, \pp}(X_2).$$ 
Alors $\varphi$ induit un isomorphisme de graphes entre les squelettes tronqués associés : $$\Gamma_1^\natural\iso \Gamma_2^\natural.$$
\end{thm}\

\begin{rem}
Comme les groupes $\pi_1^{\tp, \pp}(X_i)$ peuvent se reconstruire directement à partir des $\pi_1^\tp(X_i)$ (cf. remarque \ref{pi1t}), tout isomorphisme entre $\pi_1^\tp(X_1)$ et $\pi_1^\tp(X_2)$ induit un isomorphisme entre $\pi_1^{\tp, \pp}(X_1)$ et $\pi_1^{\tp, \pp}(X_2)$. Ainsi, par le théorème ci-dessus, tout isomorphisme entre les $\pi_1^\tp(X_i)$ induit un isomorphisme de graphes entre les squelettes tronqués. 
\end{rem}

\begin{proof}
Puisque les équivalences de catégories $\Cov^{\tp, \pp}(X_i)\simeq\BB^\tp(\pi_1^{\tp, \pp}(X_i))$ et $\BB^\tp(\pi_1^{\tp, \pp}(X_2)) \simeq \BB^\tp(\pi_1^{\tp, \pp}(X_1)) $ induisent des isomorphismes entre les différents tempéroïdes connexes correspondant, alors $\varphi$ induit un isomorphisme de tempéroïdes connexes $$\left(\Cov^{\tp, \pp}(X_1),\Upsilon_{X_1}^{\pp}\right)\iso \left(\Cov^{\tp, \pp}(X_2),\Upsilon_{X_2}^{\pp}\right).$$ 

Par le corollaire \ref{equitemp} couplé à la remarque \ref{remarque sur l'équivalence des catégories liées aux graphe tronqué et au semi-graphe} on déduit l'isomorphisme de tempéroïdes connexes : 
$$\left(\BB^\tp(\GG^\natural(X_1,S_1)),\Upsilon^{\mathrm{g\acute{e}o}}_{\GG^\natural(X_1,S_1)}\right)\iso \left(\BB^\tp(\GG^\natural(X_2,S_2)),\Upsilon^{\mathrm{g\acute{e}o}}_{\GG^\natural(X_2,S_2)}\right).$$

En ayant montré au \ref{corollaire indiquant que le granphe vérifie les hypothèses Mochizukiennes} que les $\GG^\natural(X_i,S_i)$ sont quasi-cohérents, totalement élevés, totalement détachés et sommitalement minces,  alors on déduit du théorème \ref{reconsmochi} un isomorphisme de graphes d'anabélioïdes: 
$$\GG^\natural(X_1, S_1)\iso \GG^\natural(X_2, S_2),$$
qui induit en particulier un  isomorphisme de graphes : $\G^\natural(X_1, S_1)\iso \G^\natural(X_2, S_2),$ i.e. un isomorphisme de graphes entre les squelettes tronqués : 
$\Gamma_1^\natural\iso \Gamma_2^\natural.$

\end{proof}

\begin{coro}\label{corollaire de reconstruction du squelette tronqué des courbes analytiquement anabéliennes} 
Soient $X_1$ et $X_2$ deux courbes $k$-analytiquement hyperboliques dont on note $\Sigma_{X_1}$ et $\Sigma_{X_2}$ les ensembles de nœuds des squelettes respectifs. Supposons qu'il existe un isomorphisme entre les groupes tempérés modérés : $$\varphi : \pi_1^{\tp, \pp}(X_1)\iso \pi_2^{\tp, \pp}(X_2).$$ 
Alors $\varphi$ induit un isomorphisme de graphes entre les squelettes analytiques tronqués : $$S^\an(X_1)^\natural \iso S^\an(X_2)^\natural.$$
\end{coro}

\begin{proof}
Nous savons d'après le corollaire \ref{corollaire indiquant que le granphe vérifie les hypothèses Mochizukiennes} que les graphes d'anabélioïdes $\GG^\natural(X_1, \Sigma_{X_1})$ et $\GG^\natural(X_2, \Sigma_{X_2})$ vérifient tous deux l'ensemble des hypothèses mochizukiennes. Il suffit donc de leur appliquer le théorème \ref{théorème de reconstruction du squelette tronqué} ci-dessus pour obtenir le résultat escompté. 
\end{proof}\

\begin{defi} [Courbes $k$-analytiquement hyperboliques à arêtes relativement compactes] \label{définition courbes hyperboliques à aretes relativement compactes}
Une courbe $k$-analytiquement hyperbolique sera dite \emph{à arêtes relativement compactes} lorsque toutes les composantes connexes de $X\setminus \Sigma_X $ sont relativement compactes dans $X$, ce qui revient à dire que $\CC^\infty(X, \Sigma_X)=\emptyset$, ou bien que $S^\an(X)$ est un vrai graphe, ou que $S^\an(X)=S^\an(X)^\natural$. 
\end{defi}

\begin{ex}
Toute courbe $k$-analytiquement hyperbolique compacte est à arêtes relativement compactes. Si $k_0$ est un corps local contenu dans $k$, le demi-plan de Drinfeld $\mathfrak{D}^1_{k/k_0}$ est également à arêtes relativement compactes.
\end{ex}

Le résultat suivant met en évidence le comportement \og anabélien \fg{} de telles courbes : 
\bigskip 

\begin{thm}\label{résultat anabélien pour les courbes hyperboliques à arêtes relativement compactes} 
Soient $X_1$ et $X_2$ deux courbes $k$-analytiquement hyperboliques à arêtes relativement compactes. Supposons qu'il existe un isomorphisme entre les groupes tempérés modérés : $$\varphi : \pi_1^{\tp, \pp}(X_1)\iso \pi_2^{\tp, \pp}(X_2).$$ 
Alors $\varphi$ induit (fonctoriellement en $\varphi$) un isomorphisme de graphes entre les squelettes analytiques associés : $S^\an (X_1)\iso S^\an (X_2)$.

\end{thm}

\begin{proof}
Ce résultat se déduit directement du corollaire \ref{corollaire de reconstruction du squelette tronqué des courbes analytiquement anabéliennes} et de la caractérisation des courbes $k$-analytiquement hyperboliques en terme de l'égalité entre squelette et squelette tronqué. 
\end{proof}\

\begin{prop}

Pour $i\in \lbrace 1,2 \rbrace$, soit $p_i$ un nombre premier et $k_i$ un corps de nombres $p_i$-adique, c'est-à-dire une extension finie de $\Q_{p_i}$. Supposons les groupes tempérés modérés isomorphes : $${\pi_1^{\tp, (p_1)'}(\mathfrak{D}^{1}_{\C_{p_1}/k_1})\simeq \pi_1^{\tp, (p_2)'}(\mathfrak{D}^{1}_{\C_{p_2}/k_2})},$$
alors $p_1=p_2$, et les extensions $k_1/\Q_{p_1}$ et $k_2/\Q_{p_2}$ ont même degré d'inertie : $[\widetilde{k_1} : \F_{p_1}]=[\widetilde{k_2} : \F_{p_2}]$, d'où un isomorphisme $\widetilde{k_1}\simeq \widetilde{k_2}$.
\end{prop}

\begin{proof}
Nous savons que $\mathfrak{D}^{1}_{\C_{p_i}/k_i}$ est une courbe $\C_{p_i}$-analytiquement hyperbolique à arêtes relativement compactes. Par conséquent, le théorème \ref{résultat anabélien pour les courbes hyperboliques à arêtes relativement compactes} permet d'affirmer que le groupe $\pi_1^{\tp, \pp}(\mathfrak{D}^{1}_{\C_{p_i}/k_i})$ détermine la classe d'isomorphisme du semi-graphe $S^\an(\mathfrak{D}^{1}_{\C_{p_i}/k_i})=\mathfrak{B}^1_{k_i}$ (on n'utilise pas ici la fonctorialité de $\varphi\mapsto \overline{\varphi}$ du théorème \ref{résultat anabélien pour les courbes hyperboliques à arêtes relativement compactes}). Or tous les sommets de $\mathfrak{B}^1_{k_i}$ ont une valence égale à $q_i+1$, où $q_i=p_i^{f_i}$ est le cardinal de $\widetilde{k_i}$ et $f_i=[\widetilde{k_i} : \F_{p_i}]$ est le degré d'inertie de $k_i/\Q_{p_i}$. Ainsi on en déduit $q_1=q_2$, puis $p_1=p_2$ et enfin $f_1=f_2$ en calculant $\ln q_i/\ln p_i$.
\end{proof}\

\subsubsection{Reconstruction des arêtes ouvertes}

Nous venons de montrer que l'on peut retrouver le squelette analytique \emph{tronqué} d'une courbe $k$-analytiquement hyperbolique à partir de son groupe tempéré, ce qui correspond au squelette tout entier dans les cas où les courbes en question sont à arêtes relativement compactes. Peut-on retrouver de manière générale tout le squelette d'une courbe $k$-analytiquement hyperbolique à partir de son groupe tempéré ? En nous inspirant de ce que fait Mochizuki dans \cite{M3}, nous allons définir les courbes $k$-analytiquement anabéliennes comme l'ensemble des courbes $k$-analytiquement hyperboliques vérifiant une condition technique dite d'\emph{ascendance vicinale}, que nous énoncerons, puis nous montrerons comment ces techniques de reconstruction \emph{au niveau de tous les revêtement étales finis} (sans se restreindre aux revêtements modérés) d'une courbe $k$-analytiquement anabélienne permettent de retrouver à partir du groupe tempéré de la courbe les arêtes ouvertes du squelette qui faisaient défaut dans le squelette tronqué. \\

\begin{lem}\label{lemme sur l'image réciproque des noeuds} 
Soit $X$ une courbe $k$-analytiquement hyperbolique, et $X'\xrightarrow{f} X$ un revêtement fini étale connexe. Alors $X'$ est une courbe $k$-analytiquement hyperbolique. De plus $f^{-1}(\Sigma_{X})\subseteq \Sigma_{X'}$, et $\Sigma_{X'}$ est une triangulation de $X'$ dès que $\Sigma_{X}$ est une triangulation de $X$.
\end{lem}

\begin{proof}
Le revêtement $f : X'\rightarrow X$ est fini et plat et $\Sigma_X$ et $\Sigma_{X'}$ sont respectivement l'ensemble des nœuds de $S^\an(X)$ et de $S^\an(X')$, ainsi l'on déduit directement de \cite{Duc} (Proposition $6.2.5$) que $f^{-1}(\Sigma_{X})\subseteq \Sigma_{X'}$. En particulier $\Sigma_{X'}$ est non vide. 

Soit $s'\in\Sigma_{X'}$, supposons qu'il existe $s\in \Sigma_X$ tel que $s'\in f^{-1}(\lbrace s\rbrace)$. En reprenant les notations de la définition \ref{noeuds hyperboliques}, l'hyperbolicité de $s$ se traduit par : $2g(s)+\mathrm{val}_{S^\an(X)}(s)+\widetilde{n}_s>2$. Or $f$ est fini, ainsi $s\in \partial^\an X$ implique $s'\in \partial^\an X'$, c'est-à-dire $\widetilde{n}_{s'}>0$ dès que $\widetilde{n}_s>0$. Par ailleurs $\mathrm{val}_{S^\an(X')}(s')\geqslant \mathrm{val}_{S^\an(X)}(s)$, et la formule de Riemann-Hurwitz permet d'affirmer que $g(s')\geqslant g(s)$. On en déduit que $s'$ est un nœud hyperbolique. Supposons que $s'\notin f^{-1}(\Sigma_{X})$. Cela implique que $s' \notin  \partial^\an(X')$, car dans le cas contraire la finitude de $f$ entraînerait $f(s')\in  \partial^\an(X)$ et de là $f(s')\in \Sigma_X$, ce qui n'est pas le cas. Par conséquent $\widetilde{n}_{s'}=0$. On a aussi $\mathrm{val}_{S^\an(X')}(s')\geqslant 1$ car si cette dernière était nulle, par connexité de $S^\an(X')$, on aurait $\Sigma_{X'}=\lbrace s'\rbrace$, et comme $\Sigma_{X'}\cap f^{-1}(\Sigma_{X})\neq \emptyset$ cela donnerait $s'\in f^{-1}(\Sigma_{X})$. Si $s'$ est unibranché dans $S^\an(X')$ ou de valence égale à $2$, alors en tant que nœud de $S^\an(X')$ il doit vérifier $g(s')>0$, d'où l'hyperbolicité de $s$ (si $\mathrm{val}_{S^\an(X')}(s')\geqslant 3$ le nœud $s'$ est automatiquement hyperbolique). On en déduit que $X'$ est une courbe $k$-analytiquement hyperbolique. 
\end{proof}\

Soient $X_1$ et $X_2$ deux courbes $k$-analytiquement hyperboliques. Pour $i\in \lbrace 1,2\rbrace$, notons ${\Delta_i := \pi_1^\tp(X_i)}$, et supposons donné un isomorphisme $\varphi : \Delta_1 \iso \Delta_2$. Soient $\Delta_1'\subseteq \Delta_1$ et $\Delta_2'\subseteq \Delta_2$ deux sous-groupes ouverts d'indice fini se correspondant via $\varphi$, et $X'_i\xrightarrow{f_i} X_i$ les revêtements finis étales connexes associés (en particulier on a l'égalité $\Delta_i'=\pi_1^\tp(X_i')$). 

\begin{coro}\label{corollaire de reconstruction restreinte à des sous-groupes} 
Les techniques de reconstruction du corollaire \ref{corollaire de reconstruction du squelette tronqué des courbes analytiquement anabéliennes} s'appliquent aux sous-groupes $\Delta_1'$ et $\Delta_2'$ se correspondant via $\varphi$, c'est-à-dire que l'isomorphisme $\varphi : \pi_1^\tp(X_1)\iso \pi_2^\tp(X_2)$ induit un isomorphisme de graphes d'anabélioïdes $\GG^\natural(X_1', \Sigma_{X_1'}) \iso \GG^\natural(X_2', \Sigma_{X_2'})$, donc en particulier un isomorphisme de graphes : $$S^\an(X_1')^{\natural} \iso S^\an(X_2')^{\natural}.$$

\end{coro}\

\begin{defi}
Soit $X$ une courbe $k$-analytiquement hyperbolique et $\Sigma_X$ l'ensemble non vide des nœuds de son squelette analytique. Après l'identification donnée par l'isomorphisme $\pi_1^\tp(\GG(X,\Sigma_X))\simeq \pi_1^{\tp, \pp}(X)$, nous dirons qu'un sous-groupe compact de $\pi_1^{\tp, \pp}(X)$ est un : 
\begin{itemize}
\item[•]\emph{sous-groupe vicinal} s'il est de la forme $\pi_e$ pour une arête $e$ fermée, c'est-à-dire associée à une couronne relativement compacte,
\item[•]\emph{sous-groupe cuspidal} s'il est de la forme $\pi_e$ pour $e$ une arête ouverte, en particulier c'est le groupe associé à une couronne non relativement compacte dès que $\Sigma_X$ forme une vraie triangulation de $X$.
\end{itemize}\

\end{defi}

L'étude menée des composantes coronaires nous permet d'affirmer que les sous-groupes vicinaux tout autant que cuspidaux sont isomorphes à $\widehat{\Z}^{\pp}$.

\begin{rem}
D'après \cite{M3} (théorème $3.7$), nous savons que les sous-groupes compacts maximaux de $\pi_1^{\tp, \pp}(X)$ correspondent exactement aux \emph{sous-groupes sommitaux}, c'est-à-dire de la forme $\pi_v$ pour $v\in \Sigma_X$, et que les sous-groupes vicinaux sont exactement les intersections non triviales de deux sous-groupes compacts maximaux. Un sous-groupe cuspidal n'est en revanche contenu que dans un seul sous-groupe compact maximal de $\pi_1^{\tp, \pp}(X)$. Par conséquent, un sous-groupe ne peut être à la fois vicinal et cuspidal. 

\end{rem}\

\begin{defi}
Si $G$ est un groupe et $H\subseteq G$ un sous-groupe, le \emph{commensurateur} de $H$ dans $G$, noté $C_G(H)$, est le sous-groupe de $G$ défini par : 
$$C_G(H):=\lbrace g\in G \mid  (gHg^{-1}\cap H) \; \text{est d'indice fini dans}\;  H \;\text{et dans}\; gHg^{-1}\rbrace.   $$ On dira que $H$ est \emph{commensurément terminal} dans $G$ lorsque $C_G(H)=H$. 

\end{defi}

\begin{rem}\label{propriété de restriction de la terminalité commensurable} 
La terminalité commensurable est préservée par restriction du groupe total : si $ G'\xhookrightarrow{\iota} G$ est une injection de groupes et $H'\subseteq G'$ un sous-groupe, alors $H'$ est commensurément terminal dans $G'$ dès que $\iota(H')$ est commensurément terminal dans $G$. Cela découle directement des définitions. 
\end{rem}\

\begin{lem}\label{lemme sur les groupes comensuréments terminaux} 
Si $G$ est un groupe dont $H$ est un sous-groupe commensurément terminal, alors $H$ est uniquement déterminé par n'importe lequel de ses sous-groupes d'indice fini. 
\end{lem}

\begin{proof}
Deux sous-groupes $A$ et $B$ de $G$ sont dits \emph{commensurables} si les deux indices $[A:A\cap B]$ et $[B:A\cap B]$ sont finis. Il s'agit d'une relation d'équivalence sur l'ensemble des sous-groupes de $G$, notée $\sim $. La symétrie et la réflexivité étant évidentes, la transitivité se montre comme suit : si $A, B$ et $C$ sont trois sous-groupes de $G$ vérifiant $A\sim B$ et $B\sim C$, alors  :
\begin{align*}
[A:A\cap C]\leqslant[A: A\cap B\cap C]&=[A:A\cap B]\cdot[A\cap B: A\cap B\cap C]\\
&\leqslant [A:A\cap B]\cdot [B:B\cap C]<+\infty,
\end{align*}
et l'on a également $[C:A\cap C]$ par symétrie, d'où $A\sim C$. Par ailleurs cette relation d'équivalence est stable par automorphisme de $G$, et $g\in C_G(A)$ si et seulement si $A\sim gAg^{-1}$. De là, si $A\sim B$, alors : \begin{align*}
g\in C_G(A) \Leftrightarrow A\sim gAg^{-1}\Leftrightarrow B\sim A\sim gAg^{-1}\sim gBg^{-1} \Leftrightarrow B\sim gBg^{-1} \Leftrightarrow g\in C_G(B).
\end{align*}
Ainsi deux groupes commensurables ont même commensurateur.\

Revenons au sous-groupe commensurément  terminal $H$ de $G$, et soit $V$ un sous-groupe d'indice fini de $H$. On a $H\sim V$, par conséquent $H=C_G(H)=C_G(V)$. \\

\end{proof}\

\begin{lem}\label{terminalité commensurale des sous-groupes cuspidaux} 
Soit $X$ une courbe $k$-analytiquement hyperbolique. Notons $\Delta:=\pi_1^\tp(X)$. Alors tout sous-groupe cuspidal de $\Delta^{\pp}$ est commensurément terminal dans $\Delta^{\pp}$. 
Un sous-groupe cuspidal étant compact, il est dès lors (d'après le lemme \ref{lemme sur les groupes comensuréments terminaux} précédent) uniquement déterminé par n'importe lequel de ses sous-groupes ouverts.

\end{lem}

\begin{proof}

Soit $\pi_x\subset \Delta^{\pp}$ un sous-groupe cuspidal associé à un cusp $x$ aboutissant à $s\in \Sigma_X$. En notant $\widehat{\Delta}^{\pp}$ le complété pro-$p'$ fini de $\Delta$, on a $\widehat{\Delta}^{\pp}=\pi_1^\mathrm{t}(X)$ (groupe fondamental modéré profini de $X$) et une injection $\Delta^{\pp}\hookrightarrow\widehat{\Delta}^{\pp}$ donnée par la proposition \ref{groupe tempéré d'un semi-graphe d'anabélioïdes}. D'après la remarque \ref{propriété de restriction de la terminalité commensurable}, il suffit de montrer que $\pi_x$ est commensurément terminal dans $\widehat{\Delta}^{\pp}$, cela impliquera qu'il l'est également dans $\Delta^{\pp}$. Supposons par l'absurde que $\pi_x$ n'est pas commensurément terminal dans $\widehat{\Delta}^{\pp}$, c'est-à-dire qu'il existe un élément $\sigma\in \widehat{\Delta}^{\pp}$ n'appartenant pas à $\pi_x$ et tel que $\pi_x\cap \pi_x^\sigma$ soit d'indice fini et donc ouvert dans $\pi_x$ et $\pi_x^\sigma$ (où l'on note $\pi_x^\sigma:=\sigma \pi_x \sigma^{-1}$). Puisque $\pi_x$ et $\pi_x^\sigma$ sont compacts, leurs topologies profinies sont celles induite par la topologie profinie de $\widehat{\Delta}^{\pp}$, et comme les sous-groupes ouverts distingués de $\widehat{\Delta}^{\pp}$ forment une base d'ouverts, il existe un sous-groupe ouvert distingué $\mathcal{U}\triangleleft \widehat{\Delta}^{\pp}$ vérifiant : 
\begin{equation*}
\pi_x\cap \mathcal{U}\subseteq \pi_x\cap \pi_x^\sigma\,\,\,\,\,\,\,\text{et}\,\,\,\,\,\,\,\pi_x^\sigma\cap \mathcal{U}\subseteq \pi_x\cap \pi_x^\sigma
\end{equation*}
De là on obtient : $\;\;\;\;\pi_x\cap \mathcal{U}=\pi_x\cap\pi_x^\sigma\cap\mathcal{U}=\pi_x^\sigma\cap \mathcal{U}=(\pi_x\cap \mathcal{U})^\sigma$.\

Soit $Z\xrightarrow{f} X$ le revêtement fini modéré connexe et galoisien donné par $\mathcal{U}$, de sorte que $\mathcal{U}=\pi_1^\mathrm{t}(Z)$. Soit $z$ un cusp de $Z$ (i.e. une composante associée à une arête ouverte de $S^\an(Z)$) au-dessus de $x$ tel que $z^\sigma\neq z$ (cela est possible quitte à prendre $\mathcal{U}$ suffisamment petit, car $\sigma\notin \pi_x$). Les groupes $\pi_z:= \pi_x\cap \mathcal{U}$ et  $\pi_{z^\sigma}:= (\pi_x\cap \mathcal{U})^\sigma$ sont les groupes d'inertie de $z$ et $z^\sigma$, et l'on a $\pi_z=\pi_{z^\sigma}$. On peut supposer, qui à restreindre à nouveau $\mathcal{U}$, que $Z$ admet un autre cusp $z'$ au-dessus de $x$, distinct de $z$ et $z^\sigma$.
\bigskip

Soit $s\in \partial^\an(Z)$. Comme $Z$ est quasi-lisse en $s$, il existe une courbe $k$-analytique lisse $Y_s$ et un voisinage $V_s$ de $s$ dans $Z$ tel que $V_s$ s'identifie à un domaine analytique fermé de $Y_s$. Fixons un voisinage affinoïde $W_s$ de $s$ dans $Y_s$; l'intersection ${T_s:=V_s\cap W_s}$ s'identifie à un voisinage analytique compact de $s$ dans $Z$. Quitte à restreindre chacun des $W_s$, pour $s$ parcourant $\partial^\an(Z)$, on peut supposer que les $T_s$ sont disjoints dans $Z$. Soit $Z^+$ la courbe $k$-analytique obtenue en recollant $\coprod W_s$ et $Z$ le long de $\coprod T_s$. On peut par ailleurs restreindre $Z^+$ de manière à ce que sous l'injection $S^\an(Z)\hookrightarrow S^\an(Z^+)$ les arêtes ouvertes de $Z$ (i.e. les cusps de $Z$) restent des arêtes ouvertes de $S^\an(Z^+)$.

\bigskip
Fixons un entier $\ell\geqslant 2$ premier à $p$. Comme la courbe $Z^+$ est lisse, le diagramme \ref{diagramme de suites exactes} nous permet d'exhiber un morphisme de groupes surjectif :  $$\HH^1(Z^+_\et, \mu_\ell)\xrightarrow{\qquad\theta \qquad} \mathsf{Harm}(S^\an(Z^+), \Z/\ell\Z).$$
Soient $e_z, e_{z^\sigma}$ et $e_{z'}$ les arêtes ouvertes de $S^\an(Z^+)$ correspondant aux cusps $z, z^\sigma$ et $z'$ de $Z$, et orientées vers le sommet auquel elles aboutissent. Le lemme \ref{lemme sur les cochaines harmoniques} permet d'exhiber une cochaîne harmonique $c\in \mathrm{Harm}(S^\an(Z^+), \Z/\ell\Z)$ telle que $c(e_z)\neq c(e_{z^\sigma})$. Soit $h\in \HH^1(Z_\et, \mu_\ell)$ la restriction à $Z$ d'un élément de $\theta^{-1}(\lbrace c\rbrace)$ vu comme une classe d'isomorphisme de $\mu_\ell$-torseurs analytiques sur $Z^+$. Nous savons que $h$ définit (à isomorphisme près) un revêtement modéré de $Z$. Or $\theta(h)(e_z)\neq \theta(h)(e_{z^\sigma})$ contredit l'égalité des groupes d'inertie $\pi_z=\pi_{z^\sigma}$, d'où la commensurabilité terminale de tout sous-groupe cuspidal dans  $\widehat{\Delta}^{\pp}$ et donc dans $\Delta^{\pp}$.
\bigskip

Par ailleurs tout sous-groupe cuspidal étant compact, ses sous-groupes ouverts sont d'indice fini. Il suffit dès lors d'appliquer le lemme \ref{lemme sur les groupes comensuréments terminaux} précédent pour obtenir qu'un sous-groupe cuspidal est uniquement déterminé par n'importe lequel de ses sous-groupes ouverts. 
\end{proof}\

\begin{rem}\label{remarque sous-groupe ouvert}Soit $\Delta'\subseteq \Delta=\pi_1^\tp(X)$ un sous-groupe ouvert d'indice fini, et $X'\to X$ le revêtement fini défini par ce sous-groupe. Soit $\CC'$ une couronne de $X'$ correspondant à une arête fermée $e'$ de $S^\an (X')$ s'envoyant sur une couronne $\CC$ de $X$ correspondant à une arête $e$ fermée ou ouverte de $S^\an(X)$. Quitte à identifier $\CC$ (resp. $\CC'$) au domaine analytique de $\P_k^{1,\an}$ défini par la condition $\vert T \vert \in I$ (resp. $\vert T \vert \in J$) avec $I$ (resp. $J$) un intervalle de $\R^*_+$, le morphisme envoyant $\CC$ sur $\CC'$ est défini par une fonction analytique $f$ inversible sur $\CC'$ s'écrivant $f=\sum_{i\in \Z} a_i T^i$ et telle qu'il existe $j\in \Z\setminus \lbrace 0\rbrace$ vérifiant $\vert a_j \vert r^j>\vert a_i \vert r^i$ pour tout $r\in J$ et $i\in \Z\setminus \lbrace j\rbrace$. Le morphisme $\CC'\rightarrow \CC$ est alors fini et plat de degré $\vert j \vert$ (cf. \cite{Duc}, $3.6.8$). Si $\CC'\rightarrow \CC$ est modéré, le morphisme de groupes injectif $\pi_{e'}\hookrightarrow \pi_e$ qu'il induit identifie $\pi_{e'}$ à l'unique sous-groupe d'indice $\vert j\vert$ de $\pi_e=\widehat{\Z}^{\pp}$, à savoir $\vert j \vert \widehat{\Z}^{\pp}$, qui est ouvert. 
\end{rem}

\begin{lem}\label{propriété de l'image d'un sous-groupe vicinal} 
Soit $X$ une courbe $k$-analytiquement hyperbolique. Notons $\Delta:=\pi_1^\tp(X)$. Si $\Delta'\subseteq \Delta$ est un sous-groupe ouvert d'indice fini induisant un morphisme $\iota : \Delta'^{\pp}\rightarrow \Delta^{\pp}$, alors l'image selon $\iota$ de tout sous-groupe vicinal de $\Delta'^{\pp}$ est soit triviale, soit un sous-groupe ouvert d'un sous-groupe vicinal de $\Delta^{\pp}$, soit un sous-groupe ouvert d'un sous-groupe cuspidal de $\Delta^{\pp}$, ces cas étant mutuellement exclusifs. 
\end{lem}

\begin{proof}
Le sous-groupe $\Delta'\subseteq \Delta$ définit un revêtement fini étale connexe $X'\xrightarrow{f} X$ envoyant $S^\an(X')$ sur $S^\an(X)$. Comme $f^{-1}(\Sigma_{X})\subseteq \Sigma_{X'}$ (cf. lemme \ref{lemme sur l'image réciproque des noeuds}), l'image d'une arête fermée de $S^\an(X')$ est soit hors du squelette de $X$ (auquel cas l'image du sous-groupe vicinal associé est triviale), soit contenue dans une arête fermée, soit dans une arête ouverte de $S^\an(X)$ (d'après la remarque précédente \ref{remarque sous-groupe ouvert}), d'où les cas possibles. 

\item Lorsque l'image est non triviale, le fait qu'un sous-groupe vicinal de $\Delta'^{\pp}$ ne puisse s'envoyer à la fois sur des sous-groupes de sous-groupes vicinaux et cuspidaux de $\Delta^{\pp}$ résulte directement du caractère totalement détaché du semi-graphe d'anabélioïdes $\GG(X, \Sigma_X).$
\end{proof}\

\begin{defi}[Courbes $k$-analytiquement anabéliennes]\label{définition des courbes analytiquement anabéliennes} 
Une courbe $k$-analytique non vide $X$, connexe, quasi-lisse et de groupe tempéré noté $\Delta:=\pi_1^\tp(X)$ sera dite \emph{$k$-analytiquement anabélienne} lorsqu'elle est $k$-analytiquement hyperbolique et qu'elle vérifie la condition suivante que l'on appellera propriété d'\emph{ascendance vicinale} :
\begin{itemize}\bigskip
\item[] Tout sous-groupe cuspidal de $\Delta^{\pp}$ admet un sous-groupe ouvert qui n'est autre que l'image d'un sous-groupe vicinal de $\Delta'^{\pp}$ par le morphisme naturel $\iota : \Delta'^{\pp}\rightarrow \Delta^{\pp}$ induit par l'injection d'un certain sous-groupe ouvert $\Delta'\subseteq \Delta$ d'indice fini.

\end{itemize}
\end{defi}\

\begin{thm}\label{théorème final sur les courbes maruées anabéliennes} 
Soit $X_1$ et $X_2$ deux courbes $k$-analytiquement anabéliennes. Supposons qu'il existe un isomorphisme entre les groupes tempérés :  $${\varphi : \pi_1^\tp(X_1)\iso \pi_2^\tp(X_2)}.$$ Alors $\varphi$ induit un isomorphisme de semi-graphes $\bar{\varphi} : S^\an (X_1)\iso S^\an (X_2),$ fonctoriel en $\varphi$.

\end{thm}\

\begin{proof}
Pour $i\in \lbrace 1,2 \rbrace$, notons $\Sigma_{X_i}$ l'ensemble (non vide) des nœuds du squelette analytique de $X_i$, et $\Delta_i= \pi_i^\tp(X_i) $. Le corollaire \ref{corollaire de reconstruction du squelette tronqué des courbes analytiquement anabéliennes} fournit un isomorphisme de graphes entre les squelettes tronqués : $$\widetilde{\varphi} : S^\an (X_1)^\natural\iso S^\an (X_2)^\natural.$$
Afin d'obtenir l'isomorphisme voulu entre les squelettes (non tronqués), il reste à montrer que la donnée de $\varphi$ permet de retrouver les arêtes ouvertes, c'est-à-dire d'établir une bijection entre $\CC^\infty(X_1, \Sigma_{X_1})$ et $\CC^\infty(X_2, \Sigma_{X_2})$ compatible avec $\widetilde{\varphi}$.

 Considérons une arête ouverte $x_1\in \CC^\infty(X_1, \Sigma_{X_1})$ de sous-groupe cuspidal associé $\pi_{x_1}\subset \Delta_1^{\pp}$. D'après la propriété d'ascendance vicinale dont jouit $X_1$, il est loisible de considérer un sous-groupe ouvert d'indice fini $\Delta'_1\subseteq \Delta_1$ tel que $\pi_{x_1}$ admette un sous-groupe ouvert qui est l'image d'un sous-groupe vicinal $\pi_{e_1}'$ de ${\Delta'_1}^{\pp}$ par le morphisme naturel $\iota_1 : {\Delta'_1}^{\pp}\rightarrow \Delta_1^{\pp}$. Considérons alors les groupes $\Delta'_2$ et $\pi'_2\subseteq {\Delta'_2}^{\pp}$ associés par $\varphi$ à  $\Delta'_1$ et $\pi'_{e_1}\subseteq {\Delta'_1}^{\pp}$. Le lemme \ref{propriété de l'image d'un sous-groupe vicinal} nous permet alors d'affirmer que $\pi_2'$ est un sous-groupe vicinal de ${\Delta'_2}^{\pp}$, dont l'image par le morphisme $\iota_2 : {\Delta'_2}^{\pp}\rightarrow \Delta_2^{\pp}$ est un sous-groupe ouvert d'un sous-groupe cuspidal $\pi_{x_2}$ de $\Delta_2^{\pp}$ pour un certain $x_2\in \CC^\infty(X_2, \Sigma_{X_2})$. Or, en sa qualité de commensurément terminal (voir lemme \ref{terminalité commensurale des sous-groupes cuspidaux}), $\pi_{x_2}$ est entièrement déterminé par ce sous-groupe ouvert. D'autre part $x_2$ est entièrement déterminé par $\pi_{x_2}$ puisque la propriété dont jouit $\GG^\natural(X_2, \Sigma_{X_2})$ d'être totalement détaché implique qu'une arête ouverte est entièrement déterminée par la classe de conjugaison de son sous-groupe cuspidal. Par ailleurs cette manière d'associer $x_2$ à $x_1$ ne dépend pas du choix du sous-groupe ouvert $\Delta'_1\subseteq \Delta_1$ vérifiant les propriétés requises : si $\Delta_1''$ est un autre sous-groupe ouvert d'indice fini de $\Delta_1$ permettant d'associer $x_2''\in  \CC^\infty(X_2, \Sigma_{X_2})$ à $x_1$, alors l'intersection $\Delta'_1\cap \Delta''_1$ est un sous-groupe ouvert de $\Delta_1$ d'indice fini permettant d'associer à $x_1$ un élément $x_2'''\in  \CC^\infty(X_2, \Sigma_{X_2})$, qui par le corollaire \ref{corollaire de reconstruction restreinte à des sous-groupes} doit être égal à $x_2$ tout autant qu'à $x_2''$, d'où l'égalité $x_2=x_2''$. 
 
Par conséquent $x_1\mapsto x_2$  est bien défini et induit une bijection $ \CC^\infty(X_1, \Sigma_{X_1}) \iso  \CC^\infty(X_2,\Sigma_{X_2})$ compatible avec $\widetilde{\varphi}$, d'où le résultat escompté. 
 \end{proof}\
 
 \begin{rem}
 La propriété d'ascendance vicinale fait appel à un sous groupe ouvert d'indice fini $\Delta'\subseteq \Delta$ associé à un revêtement étale $X'\xrightarrow{f} X$. Or, puisqu'il s'agit d'un sous-groupe de $\Delta=\pi_1^\tp(X)$ et non de $\pi_1^\tp(X)^{\pp}$, le morphisme $f$ n'a aucune raison d'être modéré. De là, c'est bien le groupe tempéré total $\pi_1^\tp(X)$ qui détermine le squelette analytique $S^\an(X)$ (d'après le théorème précédent), tandis que le groupe tempéré premier à $p$, à savoir $\pi_1^\tp(X)^{\pp}$, ne détermine que le squelette tronqué $S^\an(X)^\natural$. Par ailleurs, il n'y aucun espoir que le groupe $\pi_1^\tp(X)^{\pp}$ premier à $p$ puisse déterminer tout le squelette analytique $S^\an(X)$ dès que la courbe $k$-analytiquement hyperbolique $X$ n'est pas à arêtes relativement compactes. En effet, considérons une composante connexe $\mathcal{E}$ de $X\setminus \Sigma_X$ associée à une arête non relativement compacte (arête ouverte), et $Y=X\setminus \mathcal{E}$ la courbe obtenue à partir de $X$ par ablation de $\mathcal{E}$. Une telle ablation rajoute du bord, de telle sorte que $Y$ reste une courbe $k$-analytiquement hyperbolique. Par ailleurs il y a égalité des groupes tempérés premier à $p$ : $\pi_1^\tp(X)^{\pp}=\pi_1^\tp(Y)^{\pp}$, cela découle du théorème de rigidité $4.3.4$ de \cite{Ber2} appliqué à l'immersion $Y \hookrightarrow X$ ainsi que de la description des revêtement modérés d'une couronne et de l'expression de $\mathcal{E}$ comme limite inductive de couronnes. En revanche les deux squelettes analytiques sont différents en tant que semi-graphes en ce que $S^\an(Y)$ est obtenu à partir de $S^\an(X)$ par ablation de l'arête ouverte associée à $\mathcal{E}$, d'où l'impossibilité de retrouver $S^\an(X)$ à partir de $\pi_1^\tp(X)^{\pp}$.

 \end{rem}\
 
\subsubsection{Exemples de courbes $k$-analytiquement anabéliennes :}

Les courbes $k$-analytiquement hyperboliques à arêtes relativement compactes (\ref{définition courbes hyperboliques à aretes relativement compactes}) sont des exemples de courbes $k$-analytiquement anabéliennes. De là, les courbes $k$-analytiquement hyperboliques compactes tout autant que le demi-plan de Drinfeld sont des courbes $k$-analytiquement anabéliennes. Ce ne sont cependant pas les seuls exemples de telles courbes, comme le montrent les résultats suivants. 
 
 \begin{defi}[Cusps coronaires finis et discaux épointés]\label{définition cusps coronaires finis et discaux épointés} 
Soit $X$ une courbe $k$-analytiquement hyperbolique. Un cusp tel que la composante connexe de $X\setminus \Sigma_X$ qu'il définit est une vraie couronne sera dit : 
\begin{itemize}
\item[•]\emph{ coronaire fini} dès que la couronne associée est de module fini;
\item[•]\emph{discal épointé} dès que la couronne associée est isomorphe à un disque $k$-analytique privé d'un point rigide (la couronne associée est alors de module infini).
\end{itemize}
\end{defi}\

Le corps $k$ sera dit \emph{de carctéristique mixte} lorsque $\car(k)=0$ et $\car(\widetilde{k})>0$.

\begin{prop}\label{analélianité des courbes analytiquement hyperboliques à cusps finis}
 Si $k$ est de caractéristique mixte, alors toute courbe $k$-analytiquement hyperbolique dont tous les cusps sont coronaires finis est $k$-analytiquement anabélienne. 
 \end{prop}
 
 \begin{proof}
 Soit $X$ une telle courbe, et $(e_i)_{i\in I}$ l'ensemble des cusps de $X$, vus comme arêtes ouvertes de $S^\an(X)$. Par hypothèse, pour chaque $i\in I$, la composante connexe de $X\setminus \Sigma_X$ déterminée par $e_i$ est une couronne de module fini que l'on note $\CC_i$. Soit $x_i$ le sommet de $S^\an(X)$ auquel $\CC_i$ aboutit. Le corollaire \ref{corollaire permettant de trouver un revêtement d'une couronne totalement décomposé en un de ses bouts} assure l'existence d'un $\mu_p$-torseur $Y_i$ de $\CC_i$ et d'un point $z_i\in S^\an(\CC_i)$ tel que $Y_i$ est totalement déployé au-dessus de l'ouvert $r^{-1}\left( \,]z_i, x_i[ \,\right)$ et non déployé avec un seul antécédent au-dessus du complémentaire de $]z_i, x_i[$ dans $S^\an(\CC_i)$  ($r$ désigne ici la rétraction canonique de $X$ sur son squelette analytique). Soit $X'$ l'ouvert de $X$ défini par :  $$X'=X \setminus \bigcup_{i\in I}r^{-1}\left( S^\an(\CC_i)\setminus\,]z_i, x_i[ \, \right) =r^{-1}\left( S^\an(X)^\natural \right) \cup \bigcup_{i\in I} r^{-1}\left(\, ]z_i, x_i[\, \right).  $$ Soit $Y$ le revêtement fini étale de $X$ obtenu en recollant chacun des $Y_i$ au revêtement trivial de $X'$ à $p$ feuillets le long des $r^{-1}\left(\,]z_i, x_i[  \,  \right)$. Notons $y_i\in Y$ l'unique point au-dessus de $z_i$, c'est un nœud de $X'$. Si $\widetilde{x}_i\in Y$ désigne l'un des $p$ antécédents de $x_i$, alors $]y_i, \widetilde{x}_i[$ définit une arête fermée de $S^\an(X')$, dont le sous-groupe vicinal associé s'envoie sur un sous-groupe ouvert du sous-groupe cuspidal $\pi_{e_i}$, et ce pour chaque $i\in I$. Par conséquent le revêtement $X'$ assure que $X$ vérifie la propriété d'ascendance vicinale, d'où le résultat. On vient même de montrer par cette preuve que $X$ vérifie une propriété plus forte d' \og ascendance vicinale simultanée  \fg{} en ce que l'on peut trouver un revêtement étale fini, $X'$, permettant l'ascendance vicinale de tous les cusps simultanément.
 \end{proof}\

 \begin{prop}\label{anabélianité des analytifiées des courbes hyperboliques}
 Si $\mathscr{X}$ est une $k$-courbe hyperbolique et $k$ est de caractéristique mixte, la courbe $k$-analytique $\mathscr{X}^\an$ est $k$-analytiquement anabélienne. 
 \end{prop}
 
 \begin{proof}
 Cela est une reformulation directe de la propriété \emph{(iv)} apparaissant dans la preuve du résultat $3.11$ de \cite{M3}.
 \end{proof}
 
Remarquons qu'en vertu de la proposition ci-dessus, toute courbe marquée   $Y=X\setminus \E$ $k$-analytiquement hyperbolique avec $X$ compacte est $k$-analytiquement anabélienne.\\
 
 \begin{ex}
 \emph{Le cas d'une courbe elliptique.} Nous avons vu dans l'exemple \ref{exemple des courbes elliptiques épointées comme courbes hyperboliques} une manière canonique d'associer à une $k$-courbe elliptique $\mathscr{X}$ une courbe $k$-analytiquement hyperbolique : il s'agit de $\mathscr{Y}=\mathscr{X}^\an\setminus \lbrace e\rbrace$, où $e\in \mathscr{X}(k)$ est l'élément neutre de $\mathscr{X}$. Or, par ce qui précède, une telle courbe $\mathscr{Y}$ est $k$-analytiquement anabélienne. Cela permet d'associer de manière canonique une courbe $k$-analytiquement anabélienne à toute courbe elliptique sur $k$.
 \end{ex}

 \begin{coro}\label{anabélianité des courbes à semi-graphe fini ayant des cusps uniquement coronaires, soit finis soit épointés} 
 Toute courbe $k$-analytiquement hyperbolique dont le squelette analytique est un semi-graphe fini et dont les cusps sont uniquement coronaires finis ou discaux épointés est $k$-analytiquement anabélienne dès que $k$ est de caractéristique mixte.
 \end{coro}
 
 \begin{proof}
 Soit $X$ une telle courbe. En reprenant \emph{verbatim} le raisonnement de la preuve de la proposition \ref{analélianité des courbes analytiquement hyperboliques à cusps finis}, il est possible de construire un revêtement fini étale $X'$ de $X$ assurant l'ascendance vicinale de tout cusp coronaire fini de $X$. De là, quitte à prolonger tous les cusps coronaires finis en des cusps discaux épointés, on peut supposer dans le reste de la preuve que tous les cusps de $X$ sont discaux épointés. Remarquons que même si des cusps coronaires finis apparaissent dans $X'$, cela ne pose aucun problème en ce que l'ascendance vicinale des cusps ne doit pas nécessairement être \og simultanée \fg{}.
 \item Soit $\left(\varpi_i\right)_{i\in I}$ l'ensemble (fini) des cusps discaux épointés, chacun associés à une couronne $\CC_i$ s'incarnant comme un disque $\DD_i$ privé d'un point rigide $x_i$. Soit $\overline{X}$ la courbe $k$-analytique obtenue à partir de $X$ : 
 \begin{itemize}
 \item[•] en prolongeant chaque $\CC_i$ à $\DD_i$, pour tout $i\in I$;
 \item[•] pour chaque $x\in \partial^\an X$, en prolongeant en des disques $\left(\DD_{x,j}\right)_{1\leqslant j \leqslant \widetilde{n}_x}$ les $\widetilde{n}_x$ branches issues de $x$ manquant à $X$ en raison de l'appartenance de $x$ à son bord.
 \end{itemize}
 En sa qualité de compacte et sans bord, $\overline{X}$ est algébrisable : il existe une $k$-courbe projective lisse $\overline{\mathscr{X}}$ telle que $\overline{X}=\overline{\mathscr{X}}^\an$. Notons $\mathscr{X}=\overline{\mathscr{X}}\setminus \lbrace \left( x_i \right)_{i\in I}\rbrace$ : on a $S^\an(X)=S^\an(\mathscr{X}^\an)$, $\mathscr{X}$ est une $k$-courbe hyperbolique, et $X$ s'identifie à un domaine analytique fermé de $\mathscr{X}^\an$. Si $i\in I$, en considérant $\varpi_i$ comme un cusp discal épointé de $\mathscr{X}^\an$ qui est une courbe $k$-analytiquement anabélienne d'après \ref{anabélianité des analytifiées des courbes hyperboliques}, il existe un revêtement fini étale $Y_i$ de $\mathscr{X}^\an$ et une arête fermée de $S^\an(Y)$ dont le sous-groupe vicinal associé est naturellement envoyé sur un sous-groupe ouvert du sous-groupe cuspidal $\pi_{\varpi_i}$. Ainsi, la considération des revêtements $Y_i\times_{\mathscr{X}^\an} X$ permet de conclure que $X$ vérifie la propriété d'ascendance vicinale, d'où son anabélianité. 
 \end{proof}

 \bigskip
 
 Récapitulons les résultats obtenus en donnant une liste des courbes $k$-analytiquement hyperboliques que l'on sait être anabéliennes : 
 
 \begin{thm}\label{récapitulation des courbes analytiquement anabéliennes que l'on connaît}
 Soit $k$ un corps non archimédien complet algébriquement clos et non trivialement valué. Alors les trois types de courbes ci-dessous sont $k$-analytiquement anabéliennes : 
 \item
 \begin{itemize}
 \item[•]les courbes $k$-analytiquement hyperboliques à arêtes relativement compactes (\emph{ex :} le demi-plan de Drinfeld); 
 \item[•]si $k$ est de caractéristique mixte : les courbes $k$-analytiquement hyperboliques dont tous les cusps sont coronaires finis;
 \item[•]si $k$ est de caractéristique mixte : les courbes $k$-analytiquement hyperboliques dont le squelette analytique est un semi-graphe fini et dont les cusps sont coronaires finis ou discaux épointés. 
 \end{itemize}
 \end{thm}\

\begin{rem}
Dans la preuve de l'ascendance vicinale des analytifiées des $k$-courbes hyperboliques donnée par Mochizuki dans $\cite{M3}$, le morphisme que l'auteur utilise (associé à un sous-groupe ouvert $\Delta'\subseteq \Delta$ permettant l'ascendance vicinale d'un sous-groupe cuspidal de $\Delta^{\pp}$), qui n'est pas du tout modéré, est donné par la \og multiplication par $p$ au niveau de la Jacobienne \fg{}. Dans le cas des cusps coronaires finis, nous montrons la propriété d'ascendance vicinale en prenant une \og racine $p$-ième \fg{} d'une fonction coordonnée sur une couronne.  Il serait tentant de généraliser des constructions similaires à toute courbe $k$-analytiquement hyperbolique $X$. Il faudrait typiquement pouvoir fabriquer un revêtement étale fini de $X$ ayant au-dessus d'un point de type $2$ fixé un antécédent tel que l'extension résiduelle correspondante soit inséparable, ou montrer qu'une courbe $k$-analytiquement hyperbolique admet une fonction globale inversible, puis prendre le revêtement donné par sa \og racine $p$-ième \fg{}. Or la possibilité d'une telle construction ne semble pour l'instant pas évidente à l'auteur, mais il n'est pas exclu que toute courbe $k$-analytiquement hyperbolique soit $k$-analytiquement anabélienne.  
\end{rem}

Les courbes $k$-analytiquement hyperboliques dont on ne connaît pour l'instant pas l'éventuel caractère $k$-analytiquement anabélien sont de deux types : 
\begin{itemize}
\item[•] celles qui ont des cusps autres que coronaires finis ou discaux épointés;
\item[•] celles dont les cusps sont uniquement coronaires finis ou discaux épointés, qui ont au moins un cusp discal épointé mais dont le squelette analytique est un semi-graphe infini.
\end{itemize}

\end{document}